\documentclass[a4paper,12pt,twoside]{report}
\usepackage[left=2.5cm,right=2.5cm,top=3cm,bottom=2.5cm]{geometry}
\usepackage{tikz-cd}
\usepackage{tikz}
\usepackage{amsmath}
\usepackage{amsthm}
\usepackage{amssymb}
\usepackage{geometry}
\usepackage{setspace}
\usepackage[OT2,T1]{fontenc}
\usepackage{titlesec}
\usepackage{mathtools,leftindex,tensor,mhchem}
\usepackage{mathrsfs}
\usepackage{indentfirst}
\usepackage{caption}
\usepackage{hyperref}
\hypersetup{colorlinks=true,linkcolor=teal,anchorcolor=teal,citecolor=teal}
\usepackage{cleveref}
\usepackage{color}
\usepackage{etoolbox}
\usepackage{MnSymbol}
\usepackage{bbding}
\usepackage{pifont}
\usepackage[toc,page]{appendix}

\titleformat{\chapter}[hang] 
{\color{teal}\fontsize{30}{30}\bfseries} 
{{\color{teal} \fontsize{20}{20} \thechapter}} 
{1ex} 
{}
[]

\titleformat{\section}[hang]
{\color{teal}\huge\bfseries}
{{\color{teal}\huge \thesection.}}{0.5em}{}

\newtheorem{theorem}{Theorem}[section]
\newtheorem{conj}[theorem]{Conjecture}
\newtheorem{definition}[theorem]{Definition}
\newtheorem{proposition}[theorem]{Proposition}
\newtheorem{lemma}[theorem]{Lemma}
\newtheorem{corollary}[theorem]{Corollary}
\newtheorem{remark}[theorem]{Remark}
\newtheorem{example}[theorem]{Example}
\newtheorem*{theorem*}{Theorem}
\newtheorem*{proposition*}{Proposition}
\newtheorem*{conj*}{Conjecture}
\newtheorem*{definition*}{Definition}
\newtheorem*{corollary*}{Corollary}

\newcommand{\Q}{\mathbb{Q}}
\newcommand{\Z}{\mathbb{Z}}

\newcommand{\F}{\mathbb{F}}
\newcommand{\Cm}{\mathbb{C}}
\newcommand{\R}{\mathbb{R}}
\newcommand{\G}{\mathbb{G}}
\newcommand{\Hc}{\mathbb{H}}
\newcommand{\uU}{\underline{\mathrm{U}}}
\newcommand{\uH}{\underline{\mathrm{H}}}
\newcommand{\uD}{\underline{\mathrm{D}}}
\newcommand{\cris}{\mathrm{cris}}
\newcommand{\fppf}{\mathrm{fppf}}
\newcommand{\et}{\mathrm{\acute{e}t}}
\newcommand{\Et}{\mathrm{\acute{E}t}}
\newcommand{\dR}{\mathrm{dR}}
\newcommand{\Art}{\mathrm{Art}}
\newcommand{\prim}{\mathrm{prim}}
\newcommand{\Pic}{\mathrm{Pic}}
\newcommand{\Br}{\mathrm{Br}}
\newcommand{\NS}{\mathrm{NS}}
\newcommand{\HH}{\mathrm{H}}
\newcommand{\Hom}{\mathrm{Hom}}
\newcommand{\TC}{\mathrm{TC}}
\newcommand{\DM}{\mathrm{DM}}
\newcommand{\Fil}{\mathrm{Fil}}
\newcommand{\red}{\mathrm{red}}
\newcommand{\art}{\mathrm{Art}}
\newcommand{\op}{\mathrm{op}}
\newcommand{\fl}{\mathrm{fl}}
\newcommand{\Ld}{\mathrm{L}}
\newcommand{\Tor}{\mathrm{Tor}}
\newcommand{\length}{\mathrm{length}}
\newcommand{\Alb}{\mathrm{Alb}}
\newcommand{\W}{\mathbb{W}}
\newcommand{\sym}{\mathrm{sym}}
\newcommand{\rk}{\mathrm{rk}}
\newcommand{\gr}{\mathrm{gr}}
\newcommand{\spec}{\mathrm{spec}}
\newcommand{\T}{\mathrm{T}}
\newcommand{\Ker}{\mathrm{Ker}}
\newcommand{\BT}{\mathrm{BT}}
\newcommand{\IM}{\mathrm{Im}}
\newcommand{\PP}{\mathcal{P}}
\newcommand{\SSS}{\mathcal{S}}
\newcommand{\perf}{\mathrm{perf}}

\begin{document}

\begin{center}
\vspace*{6em}
\begin{doublespace}
{\textbf{{\Large \textcolor{teal}{$p$-PRIMARY TORSION OF THE BRAUER GROUP IN CHARACTERISTIC $p$}}}}
\end{doublespace}
\vspace*{7em}
{\small
{A THESIS PRESENTED FOR THE DEGREE OF}\\
{DOCTOR OF PHILOSOPHY OF IMPERIAL COLLEGE LONDON}\\
{AND THE}\\
{DIPLOMA OF IMPERIAL COLLEGE LONDON}\\
{BY}}\\
\vspace*{1em}
{\Large{YUAN YANG}}\\
\vspace*{5em}
{Department of Mathematics}\\
{Imperial College}\\

\vspace*{12em}
\begin{onehalfspace}
{\copyright  May 2025}\\ 

\end{onehalfspace}

\end{center}
\newpage

\vspace*{20em}
\large{I certify that this thesis, and the research to which it refers, are the product of my own work, 
and that any ideas or quotations from the work of other people, published or otherwise, 
are fully acknowledged in accordance with the standard referencing practices of the discipline. }

\newpage

\vspace*{4em}
\begin{flushright}
{\Huge\textbf{COPYRIGHT}}
\end{flushright}
\vspace*{3em}
\large{
The copyright of this thesis rests with the author. Unless otherwise indicated, its contents are licensed under a Creative Commons Attribution-NonCommercial 4.0 International Licence (CC BY-NC).

\noindent Under this licence, you may copy and redistribute the material in any medium or format. You may also create and distribute modified versions of the work. This is on the condition that: you credit the author and do not use it, or any derivative works, for a commercial purpose.

\noindent When reusing or sharing this work, ensure you make the licence terms clear to others by naming the licence and linking to the licence text. Where a work has been adapted, you should indicate that the work has been changed and describe those changes.

\noindent Please seek permission from the copyright holder for uses of this work that are not included in this licence or permitted under UK Copyright Law.}

\chapter*{Abstract}

Let $X$ be a proper smooth variety over an algebraically closed field $k$ of characteristic $p>0$. 
This thesis studies the 
$p$-primary component of the Brauer group of $X$, denoted by $\Br(X)[p^\infty]$. 

In the first part, we start by recalling
that \[\Br(X)[p^\infty]\cong (\Q_p/\Z_p)^{r-\rho}\oplus \HH^3(X,\Z_p(1))[p^\infty],\] 
where $r$ denotes the multiplicity of slope $1$ in $\HH^2_{\cris}(X/W)$, and $\rho$ is the Picard number of $X$.  
The term $\HH^3(X,\Z_p(1))[p^\infty]$ denotes the torsion part of the projective limit $\varprojlim_{n}\HH^3_{\fppf}(X,\mu_{p^n})$. 
By applying Illusie's theory of logarithmic de Rham-Witt sheaves, we show that
$\HH^3(X,\Z_p(1))[p^\infty]$ fits into a short exact sequence
\[0\rightarrow U(k)\rightarrow \HH^3(X,\Z_p(1))[p^\infty]\rightarrow J\rightarrow 0,\]
where $U$ is a connected commutative unipotent group over $k$ and $J$ is a finite group. 
We then review Crew's formula and Ekedahl's inequality, 
which provide information about the 
dimension of $U$. 
As an application, we obtain results concerning $U$ and $J$ for various classes of varieties. 

In the second part, we employ a formula introduced by Skorobogatov 
to explicitly compute $\Br(A)[p^n]$ for certain abelian varieties. 
Using this approach, we determine the dimension of $U[p]$, the $p$-torsion subgroup of $U$, for 
an arbitrary principally polarized abelian variety in terms of the Ekedahl-Oort type of $A$, when $p\neq 2$. 
Combined with a calculation for non-polarized abelian varieties using Kraft's cycles, 
this allows us to identify the isogeny class of $U$ for all abelian threefolds. 
We explicitly compute the connected component of the group $\Hom(A[p^2],A^\vee[p^2])^{\sym}$ in the case where
$A$ is a supergeneral principally polarized abelian threefold, 
and we show that this unipotent group varies in a well-behaved manner across the moduli space of such abelian threefolds. 
Finally, we propose a conjectural description of the isogeny class of $U$
for an arbitrary supergeneral abelian variety. 

The third part contains further applications, some conjectures and partial results in the direction
of these conjectures.
We relate the fppf cohomology of $\mu_{p^n}$ to the crystalline cohomology at level $n$.
In particular, we obtain an injectivity criterion that generalizes a result of Ogus.
We establish smoothness of flat Artin--Mazur formal groups for all abelian threefolds 
and for a large class of abelian fourfolds, 
without relying on Ekedahl's diagonal $t$-structure. 

\newpage

\begin{flushright}
    \Huge{\bf ACKNOWLEDGEMENTS}
\end{flushright}

I would like to deeply thank my supervisor Alexei Skorobogatov 
for suggesting this fascinating problem and giving me the opportunity to work on it. 
His unwavering support and enthusiasm for mathematics have 
continually inspired me and brightened my path to mathematics.

I am immensely grateful to Oliver Gregory 
for introducing me to Ekedahl's remarkable work, 
and for helping me to understand work of Artin and Mazur. 
Without his help my PhD would be much more difficult.

I would like to thank George Boxer and Ambrus P\'al for examining my early 
and late stage assessments. 
They have offered me guidance, invaluable tips and shared their perspectives. 
My sincere thanks go to Christopher Lazda, Dami\'an Gvirtz-Chen, 
and Daniel Loughran for inviting me to engaging conferences 
and generously sharing their insights.

Many thanks to Dingli Liang, Wei Zhou for their generous help in London.  
I thank Yu Min, Longke Tang, Ziquan Yang for the invaluable conversations during my PhD. 
I thank Martin Ortinez, Nina Wojtek, Bence Hevesi and Jesse Pajwani for 
organizing the London Junior Number Theory Seminar. 
I thank all LSGNT friends 
for creating such a vibrant and supportive academic community.

This work was supported by the Engineering and Physical Sciences Research Council [EP/S021590/1]. 
The EPSRC Centre for Doctoral Training in Geometry and Number Theory (The London School of Geometry and Number Theory), University College London.

\newpage
\noindent{\bf Declaration:} 
The results of Part \ref{Part1} and Section \ref{5} closely follow 
the previously appeared preprint by the author \cite{Yang} 
which is submitted for publication. 
Some of the results of Part \ref{Part2} will appear in a joint work with Livia Grammatica 
and Alexei Skorobogatov \cite{Skor}.

\tableofcontents


\chapter*{Notation}
\begin{itemize}
    \item $k$, an algebraically closed field of characteristic $p>0$.
    \item $W=W(k)$, the ring of Witt vectors of $k$.
    \item $K=W[\frac{1}{p}]$, the field of fractions of $W$. 
    \item $A[\ell^\infty]$, the $\ell$-primary torsion subgroup of any abelian group $A$. 
    \item $W_{\sigma}[F,V]$, the (non-commutative) Cartier-Dieudonné ring, defined by the relations $Fa=a^\sigma F$, $aV=Va^\sigma$, $FV=VF=p$. 
    \item $R=R^0\oplus R^1$, the Cartier-Dieudonné-Raynaud ring (Definition \ref{IRDring}). 
    \item $M[i](j)$: for $M\in D^b(R)$, shifted $i$ degrees in complex grading, shifted $j$ degrees in module grading (see Definition \ref{shifting}).
    \item $\HH^i(X,\Z_p(1)):=\varprojlim_{n}\HH^i_{\fppf}(X,\mu_{p^n})$ (\ref{HiZp(1)}). 
    \item $m_{\lambda}(M)$: the multiplicity of the slope $\lambda$ in the (iso)crystal $M$. 
    \item $h^{ij}:=\dim_k \HH^j(X,\Omega^i_{X/k})$: the Hodge numbers of $X$. 
    \item $b_i:=\dim_{\Q_{\ell}}\HH^i_{\et}(X,\Q_{\ell})$: the $i$-th Betti number of $X$. 
\end{itemize}

\chapter*{Introduction}



The (cohomological) Brauer group of a variety $X$ is defined by Grothendieck as $\Br(X):=\HH^2_{\et}(X,\G_m)$. 
It plays a fundamental role in governing both rational points and divisors (see \cite{CTS}). 
An enduring problem in arithmetic geometry is to understand the structure of $\Br(X)$ for a variety $X$ over a field $k$. 

This thesis focuses on a narrower aspect of this question: the study of
the $p$-primary part $\Br(X)[p^\infty]$, 
where $X$ is a proper smooth variety over an algebraically closed field $k$ of characteristic $p$. 
Indeed, in this setting, it is well known that the cohomological Brauer group 
$\Br(X)=\HH^2_{\et}(X,\G_m)$ is a torsion abelian group. For a prime $\ell\neq p$, 
using the Kummer sequence for the \'etale topology, Grothendieck showed that
$\Br(X)[\ell^\infty]$ is a direct sum of $(\Q_{\ell}/\Z_{\ell})^{b_2-\rho}$ 
and $\HH^3_{\et}(X,\Z_{\ell}(1))[\ell^\infty]$, where  $b_2$ be the second Betti number of $X$ 
computed using $\ell$-adic \'etale cohomology, and
$\rho=\dim_\Q(\Pic(X)\otimes\Q)$ is the Picard number of $X$ (see, e.g., \cite[Proposition 5.2.9]{CTS}). 
In contrast, the structure of the $p$-primary part $\Br(X)[p^\infty]$ is considerably more subtle. 
The goal of this thesis is to study this object in depth. 

Using the Kummer sequence for the fppf topology, one obtains a similar decomposition
for the $p$-primary torsion. Indeed, $\Br(X)[p^{\infty}]$ is also 
a direct sum of $(\Q_p/\Z_p)^{r-\rho}$ and $\HH^3(X,\Z_p(1))[p^\infty]$, 
where \begin{equation}\label{Zp1}
    \HH^i(X,\Z_p(1)):=\varprojlim_{n}\HH^i_{\fppf}(X,\mu_{p^n}),
\end{equation} 
following Illusie's definition \cite[Section II.5]{Il1}, and $r$ is the $\Z_p$-rank 
of the free part of $\HH^2(X,\Z_p(1))$ (see Section \ref{basiccalculation}).  
Thanks to foundational work by Illusie, Raynaud, Crew, and Ekedahl, the theory of the de Rham-Witt complex 
and its logarithmic subsheaves provides an effective description of the groups $\HH^i(X,\Z_p(1))$. 
The development and application of this theory form the main focus of Part \ref{Part1} of this thesis.

We begin with the following well-known structure theorem in Chapter \ref{Chap1}, based on results of Illusie and Raynaud (see Theorem \ref{Mainthm}): 

\begin{theorem*}
    Let $X$ be a smooth proper variety over an algebraically closed field $k$ of characteristic $p$. 
    Then the following statements hold. 

{\rm (1)} There is a split (non-canonically) exact sequence 
    \begin{equation*}
        0\rightarrow (\Q_p/\Z_p)^{r-\rho} \rightarrow \Br(X)[p^{\infty}] 
        \rightarrow \HH^3(X,\Z_p(1))[p^{\infty}] \rightarrow 0,
    \end{equation*} 
    where $r=m_1(\HH^2_{\cris}(X/W))$ is the multiplicity of slope $1$ of $\HH^2_{\cris}(X/W)$
.
    
{\rm (2)} For each $i\geqslant 1$ there is an exact sequence 
    \begin{equation*}
        0 \rightarrow U_i(k)\rightarrow \HH^i(X,\Z_p(1))[p^{\infty}] \rightarrow J_i \rightarrow 0,
    \end{equation*} 
    where $U_i$ is a finite dimensional connected commutative unipotent algebraic $k$-group, 
    and $J_i$ is a finite group.  
\end{theorem*}

Although the unipotent algebraic group $U_i$ is not uniquely determined by $X$, 
in the category of quasi-algebraic groups in the sense of Serre (\cite{Se})
it is well-defined up to isomorphism. 

The central objects of study in this thesis are thus the finite groups $J_i$ and the unipotent group $U_i$, 
particularly $J_3$ and $U_3$, which contrubute to the structure of $\Br(X)[p^\infty]$. 
For simplicity, we use $U$ and $J$ to denote $U_3$ and $J_3$, respectively. 

The dimension of 
the unipotent group $U$, according to Illusie's work, equals the so-called \emph{domino number} $T^{0,2}$ (see Definition \ref{Tij}). 
This number measures the non-degeneracy of the differential map 
\[\HH^2(X,WO_X)\rightarrow \HH^2(X,W\Omega^1_{X/k})\] 
in the slope spectral sequence. 
The domino numbers $T^{ij}$ are well-studied in the work of Crew and Ekedahl. 
These numbers, in some sense, measure the difference between Newton and Hodge polygons of $X$. 
In this context, we recall the notion of \emph{slope numbers} $m^{ij}$ (see Definition \ref{mij}), as introduced by Crew and Ekedahl: 

\begin{definition*}
    For a smooth proper variety $X/k$, \emph{the slope numbers} $m^{ij}$ are defined by:
\begin{equation}
    \begin{aligned}
    m^{ij}:=&\sum_{\lambda\in [i,i+1)}(i+1-\lambda)m_{\lambda}(\HH^{i+j}_{\cris}(X/W))\\
    &+\sum_{\lambda\in [i-1,i)}(\lambda-i+1) m_{\lambda}(\HH^{i+j}_{\cris}(X/W)).
    \end{aligned}
\end{equation}
\end{definition*}

We can then translate the work of Illusie, Crew and Ekedahl, into the following two nice statements (see Theorem \ref{CrewFor}, and Theorem \ref{Ekeineq}).

\begin{theorem*} Let $X$ be a proper smooth variety over an algebraically closed field $k$ of characteristic $p$. 

    {\rm (1)} When $X$ is a surface, we have Crew's formula: 
    \[\dim U=(h^{0,2}-h^{0,1})-(m^{0,2}-m^{0,1}).\]

    {\rm (2)} Assume $h^{0,2}+h^{1,1}+h^{2,0}=b_2$, then 
    \[\dim U=h^{0,2}-m^{0,2}.\]

\end{theorem*}

These results constitute the core content of Chapter \ref{Chap1}. 

In Chapter \ref{Chap2}, we present initial applications of the preceeding theory to various class of varieties.  
Rewriting Bloch-Kato's result, we show that if $X$ is ordinary, then $U=0$ (see Theorem \ref{Ordinary}). 
Using Milne's flat duality theorem for surfaces, we prove that (see Theorem \ref{surfaces}): 
\begin{theorem*}
    If $X$ is a surface, 
then $J$ is the Pontryagin dual of $\NS(X)[p^{\infty}]$. 
\end{theorem*} 
We also demonstrate that if $\HH^3_{\cris}(X/W)$ is torsion free, then $J=0$ (see Theorem \ref{AbeT}). 
Furthermore, we study $\Br(X)[p^\infty]$ for surfaces satisfying $q=-p_a$, a 
class that includes Enriques surfaces, Godeaux surfaces and (quasi)-hyperelliptic surfaces (see Theorem \ref{qpa}). 
Finally, using a K\"unneth formula due to Ekedahl, we show that if $X$ is ordinary, 
then the unipotent part of $\Br(X\times Y)$ is isomorphic to the unipotent part of $\Br(Y)$ (see Theorem \ref{prodord}). 

Our next goal, in Part \ref{Part2}, is to study the isogeny class of $U$ in depth. 
Indeed, while the dimension of $U$ is known from the work of Crew and Ekedahl, commutative unipotent groups in characteristic $p$ 
exhibit subtle and intricate behavior. For instance, the Witt group $\W_2$ (of dimension 2) is 
not annihilated by $p$, whereas $\G_a \times \G_a$ is. In fact, Serre showed that any connected 
commutative unipotent group $U$ is isogenous to a product of Witt vector groups $\W_{n_i}$, 
meaning there exists a surjective morphism $U \to \prod_i \W_{n_i}$ with finite kernel (see \cite[Theorem VII.2.10.1]{Se2}). 
The multiset of integers ${n_i}$--called the isogeny class of $U$--is uniquely determined by, 
and determines, the dimensions of the finite subgroups $U[p^n]$ killed by $p^n$ for all $n$.

In Chapter \ref{Ch3} we prove a formula of Skorobogatov that
allows one to explicitly compute $U$ in the case of abelian varieties when $p \neq 2$,
see Theorem \ref{Skorobogatov}. 
For an arbitrary abelian variety $A$ and any $n\geq 1$ there is an isomorphism
\[\Br(A)[p^n]\cong \Hom(A[p^n],A^\vee[p^n])^{\sym}/(\Hom(A,A^\vee)^{\sym}/p^n).\]
Here, $\sym$ denotes morphisms $f$ such that $f^\vee\cong f$, 
under the identifications given by the Weil pairings: 
\(w_A:A[p^n]^{\vee}\cong A^{\vee}[p^n], \quad w_{A^\vee}: A^{\vee}[p^n]^{\vee}\cong A^{\vee\vee}[p^n]\cong A[p^n].\)

In Chapter \ref{Chap4}, by solving equations for the relations of the Dieudonné modules, 
we compute the dimension of $U[p]$ for any principally polarized abelian variety $(A, \iota)$ 
from its Ekedahl-Oort type(see Theorem \ref{E-O}): 
\begin{theorem*}
Let $(A,\iota)$ be a principally polarized abelian variety of dimenision $g$  whose
Ekedahl-Oort type is represented by an elementary sequence $\varphi=\{\varphi(1),...,\varphi(g)\}$.
Setting $\varphi(0)=0$, let $P=\{1\leqslant m_1<...<m_h\leqslant g\}$ be the set of integers such that $\varphi(m)=\varphi(m-1)+1$. Then we have
    \[\dim U[p]=\frac{(g-h)(g-h-1)}{2}+\sum_{i=1}^h (m_i-i).\]
\end{theorem*}

Together with the dimension formula (\ref{Abe}), this enables us to determine the isogeny class of $U$ 
for any principally polarized abelian threefold (see Table \ref{abe3table}).
Combining this result with some more calculations for non-polarized group schemes using Kraft's classification theorem, 
we are able to classify the isogeny class of $U$ for all abelian threefolds, 
according to its Newton polygon and its $a$-number. We record the table as follows:

\begin{center}{\rm
    \begin{tabular}{ c|c|c|c|c }
    Newton polygon &    $a$      & $\dim(U[p])$& $\dim(U)$ & isogeny class of $U$\\ 
    supersingular, &    $3$      &  $3$       & $3$      & $\G_a\times \G_a\times \G_a$  \\ 
    supersingular, &    $2$      &  $3$       & $3$     & $\G_a\times \G_a\times \G_a$\\    
    supersingular, &    $1$      &  $2$       & $3$     & $\G_a\times \W_2$\\            
    $1/3$ type,    &    $2$      &  $2$       & $2$     & $\G_a\times \G_a$  \\ 
    $1/3$ type,    &    $1$      &  $2$       & $2$     & $\G_a\times \G_a$\\                 
    almost supersingular,& $2$   &  $1$       & $1$     & $\G_a$  \\ 
    almost supersingular & $1$   &  $1$       & $1$     & $\G_a$\\           
     almost ordinary   & $1$     &  $0$       & $0$     & $0$  \\ 
     ordinary      &    $0$      &  $0$       & $0$     & $0$\\           
    \end{tabular}}
\end{center}

Here, there are five different possibilities for the Newton polygon of $\HH^1_{\cris}(A/W)$, 
namely: the ordinary case (black), the almost ordinary case (blue), the almost supersingular case (orange), the $1/3$ case (green), and the supersingular case (red). 
\begin{center}
\begin{tikzpicture}\label{NewtonPolygon}
    \draw [step=1cm,gray,very thin] (0,0) grid (6,3);
    \draw (0,0) -- (3,0);
    \draw (3,0) -- (6,3);
    \draw [blue] (2,0) -- (4,1);
    \draw [orange] (1,0) -- (5,2);
    \draw [green](0,0) -- (3,1);
    \draw [green](3,1) -- (6,3);
    \draw [red](0,0) -- (6,3);
\end{tikzpicture}
\end{center}

In particular, we see from the table that for abelian threefolds, $U$ is not a product of $\G_a$'s if and only if $A$ is supersingular with $a=1$. 
Such abelian varieties were called \emph{supergeneral} in the literature. 
We propose a conjectural form of $U$ for supergeneral abelian varieties in Conjecture \ref{supergeneralconj}. 

The explicit structure of $U$, up to isomorphism (as a quasi-algebraic group) is highly intricate. 
For completeness, we include in Chapter \ref{Chap5} a detailed computation of $\Hom^{\circ}(A[p^2], A^{\vee}[p^2])^{\mathrm{sym}}$, 
the connected component of the group $\Hom(A[p^2], A^{\vee}[p^2])^{\mathrm{sym}}$
in the case of a principally polarized supergeneral abelian variety, 
illustrating this complexity. We also investigate the distribution of this Hom group across the moduli 
space of principally polarized abelian varieties (see Theorem \ref{distribution}).

In the third part, which contains only one chapter (Chapter \ref{Chap6}),
we present some related subjects and partial results of independent interest.
We establish a connection between fppf cohomology groups $\HH^i(X,\mu_{p^n})$ and the crystalline cohomology groups $\HH^i_{\cris}(X/W_n)$. 
We prove an injectivity criterion (see Theorem \ref{Inj}) that generalizes a result of Ogus.
We demonstrate the smoothness of Artin--Mazur formal groups (in particular, the formal Brauer group $\widehat{\Br(A)}$) for all abelian threefolds, and for a large class of abelian fourfolds, 
without invoking Ekedahl's deep diagonal $t$-structure theory (see Theorem \ref{flatAMFsmooth}).

\part{The General Form of $\Br(X)[p^\infty]$}\label{Part1}

\chapter{Fundamental Theory}\label{Chap1}

Throughout this chapter, $X$ is a proper smooth variety over an algebraically closed field $k$ of characteristic $p>0$. 

\section{Basic Calculation}\label{basiccalculation}

Recall that any group scheme represents an $\fppf$-sheaf (see \cite[Corollary II.1.7]{Mi}). 
In particular, we have the exact sequence of fppf sheaves
$$0\rightarrow \mu_{p^n}\rightarrow \mu_{p^{m+n}}\rightarrow \mu_{p^m}\rightarrow 0$$ for any pairs of positive integers $m,n$. 
This gives rise to two systems: the projective system $\{\mu_{p^n}\}_n$, and the injective system $\{\mu_{p^n}\}_n$. 

\begin{definition}\label{HiZp(1)}
    Define the following cohomology groups: \[\HH^i(X,\Z_p(1)):=\varprojlim_{n}\HH^i_{\fppf}(X,\mu_{p^n}), \quad \HH^i(X,\mu_{p^\infty}):=\varinjlim_{n}\HH^i_{\fppf}(X,\mu_{p^n}).\]
\end{definition}
Since free $\Z_p$-modules are projective, and divisible abelian groups are injective, 
we have the following (non-canonical) decompositions: 
\[\begin{aligned}
    \HH^i(X,\Z_p(1))&\cong \HH^i(X,\Z_p(1))_{\mathrm{free}}\oplus \HH^i(X,\Z_p(1))_{\mathrm{tors}},\\
    \HH^i(X,\mu_{p^\infty})&\cong \HH^i(X,\mu_{p^\infty})_{\mathrm{div}}\oplus \HH^i(X,\mu_{p^\infty})_{\mathrm{tors}},
\end{aligned}\]
where $\HH^i(X,\Z_p(1))_{\mathrm{tors}}$ is the torsion part of $\HH^i(X,\Z_p(1))$, and $\HH^i(X,\Z_p(1))_{\mathrm{free}}$ is the free quotient; 
$\HH^i(X,\mu_{p^\infty})_{\mathrm{div}}$ is the maximal divisible subgroup of $\HH^i(X,\mu_{p^\infty})$, and $ \HH^i(X,\mu_{p^\infty})_{\mathrm{tors}}$ is the non-divisible quotient. 
\begin{proposition}\label{proinj}The following statements hold: 
    
    {\rm (1)} Let $\T_p$ be the Tate module notation. Then we have \[\begin{aligned}
        &\HH^i(X,\Z_p(1))_{\mathrm{free}}\cong \varprojlim_n\HH^i(X,\mu_{p^\infty})_{\mathrm{div}}[p^n]\cong \T_p( \HH^i(X,\mu_{p^\infty})_{\mathrm{div}})\\
        &\HH^i(X,\mu_{p^\infty})_{\mathrm{div}}\cong \varinjlim_n (\HH^i(X,\Z_p(1))_{\mathrm{free}}/p^n)\cong \HH^i(X,\Z_p(1))_{\mathrm{free}}\otimes \Q_p/\Z_p.
    \end{aligned}\]

    {\rm (2)} There is a natural isomorphism: \[\HH^i(X,\Z_p(1))_{\mathrm{tors}}\cong \HH^{i-1}(X,\mu_{p^\infty})_{\mathrm{tors}}. \]
\end{proposition}
\begin{proof}
    Consider the following commutative diagram of exact sequences: 
\begin{center}
    \begin{tikzcd}
        0 \arrow{r} & \mu_{p^{m+1}} \arrow{r}\arrow[d,"p"] &\mu_{p^{m+n+1}} \arrow[d,"p"]\arrow[r,"p^{m+1}"] &\mu_{p^n}\arrow{r}\arrow{d} &0\\
        0 \arrow{r} & \mu_{p^m} \arrow{r}          &\mu_{p^{m+n}}          \arrow[r,"p^m"] &\mu_{p^n}\arrow{r}          &0.
    \end{tikzcd}
\end{center}
Taking the projective limit over $m$, we obtain an exact sequence: 
\[0\rightarrow \HH^i(X,\Z_p(1))/p^n\rightarrow \HH^i_{\fppf}(X,\mu_{p^n}) \rightarrow\HH^{i+1}(X,\Z_p(1))[p^n]\rightarrow 0.\]
Then taking the injective limit over $n$, we obtain:
\[0\rightarrow \HH^i(X,\Z_p(1))\otimes \Q_p/\Z_p \rightarrow \HH^i_{\fppf}(X,\mu_{p^\infty})\rightarrow \HH^{i+1}(X,\Z_p(1))[p^\infty]\rightarrow 0.\]
Now consider another commutative diagram: 
\begin{center}
    \begin{tikzcd}
        0 \arrow{r} & \mu_{p^m} \arrow{r}\arrow[d,"="] &\mu_{p^{m+n}} \arrow{d}\arrow[r,"p^{m}"] &\mu_{p^n}\arrow{r}\arrow{d} &0\\
        0 \arrow{r} & \mu_{p^m} \arrow{r}          &\mu_{p^{m+n+1}}          \arrow[r,"p^{m}"] &\mu_{p^{n+1}}\arrow{r}          &0
    \end{tikzcd}
\end{center}
Taking the injective limit over $n$, we obtain:
\[0\rightarrow \HH^{i-1}(X,\mu_{p^\infty})/p^m\rightarrow \HH^i_{\fppf}(X,\mu_{p^m}) \rightarrow\HH^{i}(X,\mu_{p^\infty})[p^m]\rightarrow 0.\]
Then taking the projective limit over $m$, 
\[0\rightarrow  \HH^{i-1}(X,\mu_{p^\infty})_{\mathrm{tors}}\rightarrow \HH^i(X,\Z_p(1))\rightarrow \T_p(\HH^{i}(X,\mu_{p^\infty}))\rightarrow 0.\]
Here, as will be discussed in Theorem \ref{rep}, the cohomology groups $\HH^i_{\fppf}(X,\mu_{p^n})$ are the group of $k$-points of finite type commutative group schemes, 
hence satisfy Mittag-Leffler conditions, so taking projective limits preserves exactness. 
\end{proof} 

We also have the Kummer sequence of fppf sheaves: 
\[0\rightarrow \mu_{p^n} \rightarrow \G_m \xrightarrow{p^n} \G_m\rightarrow 0.\]
Taking fppf cohomology, this yields the short exact sequence: 
\[0\rightarrow \Pic(X)/p^n\rightarrow \HH^2_{\fppf}(X,\mu_{p^n})\rightarrow \Br(X)[p^n] \rightarrow 0.\]
Passing to the direct limit over $n$, we obtain:
\[0\rightarrow \NS(X)\otimes (\Q_p/\Z_p) \rightarrow \HH^2_{\fppf}(X,\mu_{p^\infty}) \rightarrow \Br(X)[p^\infty]\rightarrow 0.\]

We now consider the following commutative diagram of exact sequences:
\begin{small}
\begin{center}
    \begin{tikzcd}
        0\arrow{r} & \NS(X)\otimes \Q_p/\Z_p \arrow{r}\arrow{d} &\HH^2_{\fppf}(X,\mu_{p^\infty}) \arrow{r}\arrow{d} &\Br(X)[p^\infty]\arrow{r}\arrow{d} &0\\
        0\arrow{r} & \HH^2(X,\Z_p(1))\otimes \Q_p/\Z_p \arrow{r} &\HH^2_{\fppf}(X,\mu_{p^\infty}) \arrow{r} &\HH^3(X,\Z_p(1))[p^\infty]\arrow{r}     &0
    \end{tikzcd}
\end{center}
\end{small}Applying the snake lemma to this diagram, we deduce the exact sequence:
\[0\rightarrow \mathrm{coker} \rightarrow \Br(X)[p^\infty]\rightarrow \HH^3(X,\Z_p(1))[p^\infty] \rightarrow 0,\]
where $\mathrm{coker}$ is the cokernel of the natural map $\NS(X)\otimes \Q_p/\Z_p\rightarrow\HH^2(X,\Z_p(1))\otimes \Q_p/\Z_p$.
As we will see later in \ref{Mainthm}, $\mathrm{coker}\cong (\Q_p/\Z_p)^{r-\rho}$, where $r$ is the multiplicity of slope 1 in $\HH^2_{\cris}(X/W)$ 
and $\rho$ is the Picard number of $X$. 

To further motivates the study of the de Rham-Witt complex, 
we will show that under the $d\log$ map, $\HH^i(X,\Z_p(1))$ is identfied with the logarithmic Hodge-Witt cohomology groups 
$\HH^{i-1}(X,W\Omega_{X/k,\log}^1)$, which in turn corresponds to the Frobenius invariant subgroup 
of $\HH^{i-1}(X,W\Omega^1_{X/k})$. 
The \emph{coherence} of the Hodge-Witt cohomology groups then yields the desired decomposition of $\Br(X)[p^\infty]$. 
To set the stage, we begin by reviewing Illusie's definition of the de Rham-Witt complex, the theory of logarithmic de Rham-Witt sheaves, 
the coherence theory of Hodge-Witt cohomology, and finally Crew's formula and Ekedahl's inequality.

\section{De Rham-Witt Complex}

Let $X/k$, as before, be a proper smooth vareity over an algebraically closed field $k$ of characteristic $p$. 
The de Rham-Witt complex $W_n\Omega_{X/k}^\bullet$ is
a complex of (Zariski or \'etale) sheaves of $W_nO_X$-modules on $X$. 
As $n$ varies, these form a pro-complex of sheaves, 
and their limit is denoted $W\Omega_{X/k}^\bullet$. 
They are defined as the initial object in the category of $V$-pro complexes, see \cite[Definition I.1.4]{Il1}, or Definition \ref{Vpro}. 
These sheaves are equipped with Frobenius and Verschiebung operators,
$F: W_{n+1}\Omega^i_{X/k}\rightarrow W_n\Omega^i_{X/k},\quad V: W_n\Omega^i_{X/k}\rightarrow W_{n+1}\Omega^i_{X/k}$ 
which satisfy the standard relations: 
\[FV=VF=p, \quad Fa=a^{\sigma}F, \quad Va=a^{\sigma^{-1}}V, \quad FdV=d.\]
A fundamental feature of the de Rham-Witt complex is the existence of a canonical quasi-isomorphism
\[Rf_{\cris,*}(O_{X/W_n}) \simeq W_n\Omega_{X/k}^\bullet,\]
where  $f_{\cris}:Cris(X/W_n) \rightarrow Zar_{X}$ is the projection from the crystalline site to the 
Zariski site, and $O_{X/W_n}$ is the structure crystal on $Cris(X/W_n)$ \cite[Theorem II.1.4]{Il1}. 
Passing to the limit, we obtain a quasi-isomorphism $Rf_{\cris,*}(O_{X/W}) \simeq W\Omega_{X/k}^\bullet$. 
On the right hand side, define the absolute Frobenius operator $\varphi:=p^iF$ on $W\Omega^i_{X/k}$, 
then this is compatible with the Frobenius on the crystalline cohomology side. 
Consequently, the hypercohomology $\Hc^*(W_n\Omega_{X/k}^\bullet)$ 
computes the crystalline cohomology $\HH^*_{\cris}(X/W_n)$, 
and we obtain the so-called slope spectral sequences:  
\begin{equation}
    \label{SSf} \begin{aligned}
        E_1^{pq}: \quad\HH^q(X,W_n\Omega^p_{X/k})\Longrightarrow \HH^{p+q}_{\cris}(X/W_n),\\
    \end{aligned}
\end{equation}
\begin{equation}
    \label{SS} \begin{aligned}
        E_1^{pq}:\quad\HH^q(X,W\Omega^p_{X/k})\Longrightarrow \HH^{p+q}_{\cris}(X/W).
    \end{aligned}
\end{equation}

Since there is a canonical filtration (\ref{CanFil}) on $W\Omega_{X/k}^i$ 
whose graded pieces can be expressed as extensions of coherent $O_X^{(p^n)}$ modules (Proposition \ref{Gradedpieces}), 
one can verify that the projective systems $W_\bullet\Omega^i_{X/k}$ and $\HH^j(X,W_\bullet\Omega^i_{X/k})$ 
satisfy the Mittag-Leffler condition. 
Consequently, for any pair $(i,j)$ there is a natural isomorphism 
$\varprojlim_n \HH^j(X,W_n\Omega^i_{X/k})\cong \HH^j(X,\varprojlim_{n}W_n\Omega^i_{X/k})$. 
We denote this cohomology group by $\HH^j(X,W\Omega^i_{X/k})$. 
The main result in \cite{Il1} can be stated as follows: 

\begin{theorem}[Illusie \cite{Il1}]\label{SSbasic}
    The slope spectral sequence $(\ref{SS})$ degenerate at $E_1$ page after inverting $p$. 
    Moreover, the isocrystal $\HH^j(X,W\Omega_{X/k}^i)\otimes_W K$ coincide with the slope $[i,i+1)$ part of $\HH^{i+j}_{\cris}(X/W)\otimes_WK$, with the Frobenius action twisted by $p^i$. 
\end{theorem}

However, without inverting $p$--that is, when one takes into account the $p$-torsion, the spectral sequence $(\ref{SS})$ 
fails to degenerate at the $E_1$ page. 
In fact, it is the presence of infinite torsion submodules that 
obstructs degeneration. 
These torsion submodules are organized into structures known as \emph{dominos}. 
While they satisfy certain finiteness conditions, their size can be measured by the so-called domino numbers $T^{ij}$ (see Definition \ref{Tij}). 

\medskip

\section{Interlude: Flat and de Rham-Witt cohomology} \label{WOmegalog}

The goal of this section is to define the logarithmic de Rham-Witt sheaves $W\Omega_{X/k,\log}^i$, 
and to explain that the cohomology group $\HH^i(X,\Z_p(1))$ 
is isomorphic to the (\'etale) cohomology group $\HH^{i-1}(X,W\Omega^1_{X/k,\log})$. 

Recall (see \cite[Section IV.3]{IR}) that for each $i$, the \'etale sheaf $W\Omega_{X/k,\log}^i$ (resp. $W_n\Omega_{X/k,log}^i$) 
is defined as the subsheaf of $W\Omega^i_{X/k}$ (resp. $W_n\Omega^i_{X/k}$)
\'etale locally generated by sections of the form 
\[d\underline{x_1}/\underline{x_1}\wedge...\wedge d\underline{x_i}/\underline{x_i}\]
where $x_1,...,x_i\in O_X^{\times}$ and $\underline{x_i}=(x_i,0,...)$ is the Teichm\"uller lift of $x_i$ to $WO_X$(resp. $W_nO_X$). 
Then the flat cohomology group $\HH^i(X,\Z_p(1))$ (resp. $\HH^i_{\fppf}(X,\mu_{p^n})$) 
is related to the cohomology of $W\Omega^1_{X/k}$(resp. $W_n\Omega^1_{X/k,\log}$) via the folloing construction. 
Recall that the Kummer sequence is exact for fppf sheaves:
\[0 \rightarrow \mu_{p^n} \rightarrow \G_m \rightarrow \G_m \rightarrow 0. \] 
Let $\epsilon: X_{\fl} \rightarrow X_{\et}$ be the natural projection from the fppf site to the \'etale site. 
By a theorem of Grothendieck, higher cohomology groups of the complex $R^i\epsilon_* \G_m$ vanish for $i\geqslant 1$, 
and since $X$ is reduced, we also have $\epsilon_*(\mu_{p^n})=0$.
Applying $\epsilon_*$ to the Kummer sequence yields the short exact sequence of \'etale sheaves
\[1 \rightarrow \G_m \rightarrow \G_m \rightarrow R^1\epsilon_*\mu_{p^n} \rightarrow 0,\]
along with $R^i\epsilon_*(\mu_{p^n})=0$ for $i\geqslant 2$. 
The Leray spectral sequence then gives an isomorphism $\HH^{i}_{\fppf}(X,\mu_{p^n})\cong \HH^{i-1}_{\et}(X,\G_m/\G_m^{p^n})$ for $i\geqslant 1$.

We now define the $d \log$ map: 
\begin{equation}
    d\log\colon \G_m \rightarrow W_n\Omega_{X/k}^1, 
\end{equation} 
which locally sends a section $x$ to $d\underline{x}/\underline{x}$, 
where $\underline{x}=(x,0,...,0)$ is the Teichm\"uller lift of $x$ to $W_n(O_X)$. 

This induces an isomorphism of sheaves (\cite{Il1}, Corollary I.3.27, Proposition \ref{GmWOmega1log})
\begin{equation}\label{dlog1}
    d\log: \G_m/\G_m^{p^n} \stackrel{\sim}\longrightarrow W_n\Omega_{X/k,\log}^1,
\end{equation} 
and hence an isomorphism on the cohomology groups
\begin{equation*}
    d\log: \HH^i_{\fppf}(X,\mu_{p^n})\xrightarrow[]{\cong} \HH^{i-1}_{\et}(X, \G_m/\G_m^{p^n}[-1]) \cong \HH^{i-1}(X,W_n\Omega^1_{X/k,\log}).
\end{equation*}
Passing to the projective limit over $n$, we obtain the following well-known result:

\begin{proposition}\label{log}
    We have an isomorphism: 
    \begin{equation}
        \HH^i(X,\Z_p(1))\cong \HH^{i-1}(X,W\Omega^{1}_{X/k,\log}).
    \end{equation}    
\end{proposition}

The fact that the \(d\log\) map (\ref{dlog1}) is an isomorphism is far from trivial. 
Its proof relies on an inductive argument and a careful analysis 
of the pieces of the canonical filtration of the de Rham-Witt complex. 
Full details are provided in Appendix A \ref{GmWOmega1log}.  


%



\section{Coherent $R$-modules}

In this section, we review the structural results of the slope spectral sequence introduced in Section $\ref{SS}$. 
We begin by recalling the definition of the graded \emph{Cartier-Dieudonné-Raynaud} ring $R$, 
along with the notion of a coherent graded $R$-module. 

\begin{definition}\label{IRDring}
    Define $R^0\cong W_{\sigma}[F,V]$, the (non-commutative) ring generated by semi-linear Frobenius $F$ 
    and Verschiebung $V$, subject to the relation $FV=VF=p$. 
    Let $R^1$ be the two-sided $R^0$-module generated by an element $d$, 
    satisfying the relations $d^2=0$, $FdV=d$. 
    The resulting $\Z$-graded ring $R=R^0\oplus R^1$ is called the {\rm Cartier--Dieudonn\'e--Raynaud ring}. 
\end{definition}

More concretely, the components of the ring $R$ can be described as: 
\begin{equation}
    \begin{aligned}
        R^0\cong \{\sum_{n>0}a_{-n}V^n+ a_0+\sum_{n>0} a_n F^n, \text{almost all $a_n$ are $0$}\},\\
        R^1\cong \{\sum_{n>0}a_{-n}dV^n+a_0d+\sum_{n>0}a_nF^nd, \text{almost all $a_n$ are $0$}\}. 
    \end{aligned}
\end{equation}

Given a complex $...\rightarrow M^i\rightarrow M^{i+1} \rightarrow ...$ 
where each $M^i$ is a $R^0$-module and the compatibility relation $FdV=d$ holds, 
one can endow the total module $\oplus M^i$ with a graded (left) $R$-module structure. 
Conversely, any such graded (left) $R$-module gives rise to a compatible complex.  

For each graded $R$-module $M$, define the shifted module $M(i)$ by $M(i)^n=M^{i+n}$, and differential $d_{M(i)}=(-1)^id_{M}$. 

\begin{remark}
    Our converntions of shifting follow those of Ekedahl \cite{Ek1,Ek2,Ek3}. 
    In contrast, Illusie and Raynaud \cite{Il1,IR} denote the degree shift by $M[i]$. 
\end{remark}

Now we introduce the fundamental object known as the domino $U_t$.

\begin{definition}\label{Ut}
Let \[U_t: k[[V]]\rightarrow k[[dV]]\] be the following $R$-module concentrated in degree $0$ and $1$: 
\begin{itemize}
    \item $k[[V]]:=k_{\sigma}[[V]]$ denote the $k$-vector space of formal power series $\sum_{n\geqslant 0} a_n V^n$, 
         subject to the relation $Va=a^{\sigma^{-1}}V$ for all $a\in k$. 
        On this space, the Frobenius $F$ acts as $0$. 
    \item  $k[[dV]]:=k_{\sigma}[[dV]]$ denote the $k$-vector space of formal power series $\sum_{n\geqslant 0}a_n dV^n$  
        where $V$ acts as $0$, and 
        $F$ acts by mapping $dV^n$ to $dV^{n-1}$ if $n\geq 1$, and by mapping $d$ to $0$, and $Fa=a^{\sigma}F$ for all $a\in k$. 
    \item For an integer $t$, define $U_t$ to be the $\Z$-graded $R$-module concentrated in degrees 0 and 1, 
with $U_t^0\cong k_{\sigma}[[V]]$ and $U_t^1\cong k_{\sigma}[[dV]]$, and 
differential given by $V^n \rightarrow dV^{n-t}$, where we set $dV^{n-t}=0$ whenever $n-t<0$.
\end{itemize} 

\end{definition}
\begin{definition}
    \label{Coherence}
A $\Z$-graded $R$-module $M$ is called {\em coherent}, if it admits a finite filtration 
by $\Z$-graded $R$-modules whose successive quotients 
are isomorphic to one of the following:
    \begin{itemize}
        \item $I_a$: a finite length $W$-module concentrated in a single degree, with $V$ nilpotent;
        \item $I_b$: a finitely generated free $W$-module concentrated in a single degree, 
        with $V$-topologically nilpotent;
        \item $II$: a shifted domino module $U_t(j)$.
    \end{itemize}
\end{definition}

The following result, due to Illusie, Raynaud and Ekedahl, is fundamental: 

\begin{theorem}[Illusie, Raynaud, Ekedahl]\label{coherent}
    The graded $R$-module associated to the complex
    \begin{equation*}\label{Cplx}
        \HH^j(X,WO_{X/k}) \xrightarrow{d} \HH^j(X,W\Omega^1_{X/k}) \xrightarrow{d} ... \xrightarrow{d} 
        \HH^j(X,W\Omega^n_{X/k})
    \end{equation*} is coherent. 
\end{theorem}

\begin{proof}
    See \cite{IR}.
\end{proof}

We now recall several invariants associated with the Hodge--Witt cohomology groups $\HH^j(X,W\Omega^i_{X/k})$, 
for $i,j\geqslant 0$, as introduced by Illusie and Raynaud \cite{IR} and further studied by Crew \cite{Cr} and Ekedahl \cite{Ek3}. 

\begin{itemize}
    \item 
The \emph{heart} of $\HH^j(X,W\Omega^i_{X/k})$, denoted
\begin{equation*}
    \text{c\oe ur}(\HH^j(X,W\Omega^i_{X/k})):=V^{-\infty}Z_2^{ij}/F^{\infty}B_2^{ij},
\end{equation*}
where
\begin{equation*}
    V^{-\infty}Z_2^{ij}:=\bigcap_{n} \mathrm{ker}(dV^n: \HH^j(X,W\Omega^i_{X/k})\rightarrow \HH^j(X,W\Omega^{i+1}_{X/k})),
\end{equation*} 
\begin{equation*}
    F^{\infty}B_2^{ij}:=\bigcup_{n} F^n \mathrm{im}(d: \HH^j(X,W\Omega^{i-1}_{X/k})\rightarrow \HH^j(X,W\Omega^{i}_{X/k})).
\end{equation*}

One can show that, \(V^{-\infty}Z_2^{ij}\) is the maximal sub $R^0$-module of 
\(\mathrm{ker}(d:\HH^j(X,W\Omega^i_{X/k})\rightarrow \HH^j(X,W\Omega^{i+1}_{X/k}))\), 
and $F^\infty B_2^{ij}$ is the minimal $R^0$-module containing \(\mathrm{im}(d: \HH^j(X,W\Omega^{i-1}_{X/k})\rightarrow \HH^j(X,W\Omega^{i}_{X/k}))\). 
It was shown in \cite[Theorem II.3.1]{IR} that the heart $\text{c\oe ur}(\HH^j(X,W\Omega^i_{X/k}))$ 
is a finitely generated $W$-module. 

\item The $(i,j)$-th domino: \[\HH^j(X,W\Omega^i_{X/k})/V^{-\infty}Z_2^{ij}\rightarrow F^\infty B_2^{i+1,j},\]
which is an iterated extension of modules of the form $U_t$. 
The number of such summands is called the \emph{dimension} of the domino, denoted $T^{ij}$ (see below). 

\item The integers  \begin{equation}
    \label{Tij} T^{ij}:=\dim_k \HH^j(X,W\Omega^i_{X/k})/(V^{-\infty}Z_2^{ij}+
V\HH^j(X,W\Omega^i_{X/k})).
\end{equation} $T^{ij}$ measure the dimension of the $(i,j)$-th domino. 
They count the number of time that modules of the form $U_t(-i)$ for varying $t$ appear in the successive 
quotients of any filtration satisfying the coherence condition, see \cite[Proposition I.2.18]{IR}. 
Notably, we will see that $T^{02}$ is the dimension of the unipotent part $U$ in the Brauer group. 
\end{itemize}

The following statements, due to Illusie, Raynaud \cite{IR}, and Ekedhal \cite{Ek1}, then summarize the first properties of these invariants. 

\begin{theorem}[Illusie]\label{HodgeWittn}
    If for some $n$, the groups \(\HH^j(X,W\Omega^i_{X/k})\) are finitely generated over \(W\) for all \(i+j=n\), 
    then we have a canonical decomposition: 
    \[\HH^n_{\cris}(X/W)\cong \bigoplus_{i+j=n}\HH^j(X,W\Omega^i_{X/k}),\] 
    and the Frobenius on \(\HH^n_{\cris}(X/W)\) corresponds to \(p^iF\) on \(\HH^j(X,W\Omega^i_{X/k})\). 
\end{theorem}

\begin{proof}
    See \cite[Theorem IV.4.5]{IR}.
\end{proof}

This is called Hodge-Witt in degree $n$. 

\begin{proposition}[Illusie, Ekedahl]\label{Survivesofcoeur}
    The heart $\text{{\rm c\oe ur}}(\HH^j(X,W\Omega^i_{X/k}))$ survives to the infinity page. That is, 
    there is a canonical sequence of inclusions: 
    \[B_\infty^{ij}\subseteq F^{\infty}B_2^{ij}\subseteq V^{-\infty} Z_2^{ij}\subseteq Z_\infty^{ij}, \]
    and the quotients $F^{\infty}B_2^{ij}/B_\infty^{ij}$ and $Z_\infty^{ij}/V^{-\infty} Z_2^{ij}$ are finite length as $W$-modules.
\end{proposition}

\begin{proof}
    See \cite[Theorem II.3.4]{IR}. 
\end{proof}

\begin{proposition}[Illusie, Ekedahl]\label{Ekedahlduality}\label{EasyVan}
    The domino numbers $T^{ij}$ satisfy Ekedahl's duality \[T^{ij}=T^{N-i-2,N-j+2},\] 
    where $N$ is the dimension of $X$. 
    As a result, we have $T^{ij}=0$ for $i\geqslant N-1$ or $j\leqslant 1$.
\end{proposition} 

\begin{proof}
    See \cite[Corollary II.2.17, Propositions II.2.19 and II.3.12]{Il1},\cite[IV. Corollary 3.5.1]{Ek1}.
\end{proof}

Recall that we have defined the logarithmic de Rham-Witt sheaves \(W\Omega^i_{X/k,\log}\) in Section \ref{WOmegalog}. 
We have the following main theorem concerning the cohomology groups of \(W\Omega^i_{X/k,\log}\).

\begin{theorem}[Illusie-Raynaud]\label{IllusieRaynaud}
    Let $X$ be a smooth proper variety over $k$. Then we have the following statements.

{\rm  (1)} For $i,j\geqslant 0$, there is a short exact sequence 
    \[0\rightarrow \HH^j(X,W\Omega_{X/k,\log}^i) \rightarrow \HH^j(X,W\Omega^{i}_{X/k}) \xrightarrow{1-F} \HH^j(X,W\Omega^{i}_{X/k}) \rightarrow 0.\] 

{\rm  (2)} For $i,j\geqslant 0$, there is a short exact sequence 
    \[0 \rightarrow \HH^j(X,W\Omega_{X/k,\log}^i)^0 \rightarrow \HH^j(X,W\Omega_{X/k,\log}^i) \rightarrow D^{ij} \rightarrow 0,\]
    where $\HH^j(X,W\Omega_{X/k,\log}^i)^0$ is the group of $k$-points of a connected commutative unipotent algebraic $k$-group of dimension $T^{{i-1},j}$,
    and $D^{ij}$ is the finitely generated $\Z_p$-module
$ (\text{\rm {c\oe ur}}(\HH^j(X,W\Omega^i_{X/k})))^{F=1}$. 
\end{theorem}

\begin{proof} This is \cite[Theorem IV.3.3]{IR}. (1) There is an exact sequence of pro-sheaves: 
\[0\rightarrow W_{\bullet}\Omega^i_{\log,X/k} \rightarrow W_\bullet\Omega^i_{X/k} \xrightarrow{F-1}  W_\bullet\Omega^i_{X/k} \rightarrow 0.\]
The cohomology groups of these sheaves at each level are of finite type, so they satisfy the Mittag-Leffler condition. 
By taking limit and taking long exact sequences, 
it suffices to show that $F-1$ is surjective on $\HH^j(X,W\Omega^i_{X/k})$. 
By Theorem \ref{coherent}, the cohomology group $\HH^j(X,W\Omega^i_{X/k})$ 
is an iterated extensions of $k_{\sigma}[[V]]$, $k_{\sigma}[[dV]]$ and finitely generated $W$-modules. 
We can see that $F-1$ is surjective on any such iterated extensions of $k_{\sigma}[[V]]$, 
$k_{\sigma}[[dV]]$ and finitely generated $W$-modules, by d\'evissage. 

(2) Consider the filtration $F^\infty B_2^{ij}\subseteq V^{-\infty}Z_2^{ij}\subseteq \HH^j(X,W\Omega^i_{X/k})$ 
on the cohomology group $\HH^j(X,W\Omega^i)$. 
$F^\infty B_2^{ij}$ is an iterated extention of $k[[dV]]$'s, $V^{-\infty}Z_2^{ij}/F^\infty B_2^{ij}$ 
is the heart of $\HH^j(X,W\Omega^i_{X/k})$ which is finitely generated, 
The map $F-1$ is surjective on all the graded pieces of this filtration; 
$$\HH^j(X,W\Omega^i_{X/k})/ V^{-\infty}Z_2^{ij}$$ is an iterated extension of $k[[V]]$'s.
The action of $F-1$ is surjective on all of these graded pieces, 
so the kernel of $F-1$ is the iterated extension of the kernel of $F-1$ on each pieces. 
The kernel of $F-1$ on $k[[dV]]$ is the group $\G_a(k)$, the kernel of $F-1$ on $\text{c\oe ur}(\HH^j(X,W\Omega^i_{X/k}))$ is the group $D^{ij}$ in the proposition, 
the kernel of $F-1$ on $k[[V]]$ is trivial. 
\end{proof}
    
Although the unipotent algebraic groups in Theorem \ref{IllusieRaynaud} are not uniquely determined by $X$, 
if we pass to the category of quasi-algebraic group in the sense of Serre (see \cite{Se}, or Section \ref{Secrep}),
we can naturally select canonical ones. 
Moreover, the isomorphism of Proposition \ref{log} between abstract abelian groups 
can be refined to an isomorphism between pro-quasi-algebraic groups. 

More precisely, consider the direct image map 
\[f_{\perf,*}: (X/k)_{\perf}\rightarrow k_{\perf}.\]
Here \((X/k)_{\perf}\) denotes the relative perfect site: 
its objects are pairs \((Y,T)\) where \(T\) is a perfect \(k\)-scheme and 
\(Y\) is \'etale over \(X_T\). The site \(k_{\perf}\) is the perfect site over \(k\), 
whose objects are perfect \(k\)-schemes, and equipped with \'etale topology. 
Then \(W_n\Omega_{X/k}^i\) can be regarded as a sheaf on \((X/k)_{\perf}\), and 
\(R^jf_{\perf,*}(W_n\Omega_{X/k}^i)\) is precisely the sheaf on 
\((k)_{\perf}\) represented by the \(W_n\)-module scheme 
\(\HH^j(X,W_n\Omega^i_{X/k})\), 
which assigns to every perfect \(k\)-algebra \(A\) the group \(\HH^j(X,W_n\Omega^i_{X/k})\otimes_{W}W(A)\). 
The exact sequence of pro-sheaves: 
\[0\rightarrow W_{\bullet}\Omega^i_{\log,X/k} \rightarrow W_\bullet\Omega^i_{X/k} \xrightarrow{F-1}  W_\bullet\Omega^i_{X/k} \rightarrow 0,\]
induces a long exact sequence of pro-sheaves on \(k_{\perf}\): 
\[...\rightarrow R^jf_{\perf,*}(W\Omega^i_{\log,X/k})\rightarrow R^jf_{\perf,*}(W\Omega_{X/k}^i) \xrightarrow{F-1}R^jf_{\perf,*}(W\Omega_{X/k}^i) \rightarrow....\]

Similar to Theorem \ref{IllusieRaynaud} (2), the map \(F-1\) on \(R^jf_{\perf,*}(W\Omega_{X/k}^i)\) is surjective. 
Hence the sheaf \(R^jf_*(W\Omega^i_{\log,X/k})\) is precisely the sheaf represented by the kernel of \(F-1\): 
a pro-quasi-algebraic group over \(k\). 

In this setting, Theorem \ref{log} lifts to an isomorphism between pro-quasi-algebraic groups. 
Indeed, consider the following commutative diagram of sites:  

\[\begin{tikzcd}
	{X_{\fppf}} & {X_{\Et}} & {X_{\perf}} \\
	{k_{\fppf}} & {k_{\Et}} & {k_{\perf},}
	\arrow["{\beta_X}", from=1-1, to=1-2]
	\arrow["{f_{\fppf}}"', from=1-1, to=2-1]
	\arrow["{\alpha_X}", from=1-2, to=1-3]
	\arrow["{f_{\Et}}"', from=1-2, to=2-2]
	\arrow["{f_{\perf}}"', from=1-3, to=2-3]
	\arrow["{\beta_{k}}", from=2-1, to=2-2]
	\arrow["{\alpha_k}", from=2-2, to=2-3]
\end{tikzcd}\]
we have \[R(\alpha_{k}\circ\beta_k)_*\circ Rf_{\fppf *}(\mu_{p^n})\cong Rf_{\perf,*}(W_n\Omega^1_{\log,X/k}[-1]).\]
Moreover, for every group scheme \(G\) over \(k\), \(R(\alpha\circ\beta)_*(G)\) will be the sheaf represented by the quasi-algebraic group associated to \(G\) (see, for example, \cite[p. 221]{Be}). 
As a result, the quasi-algebraic group associated to \(R^if_{\fppf,*}(\mu_{p^n})\) 
(this sheaf is represented by a group scheme, see Theorem \ref{rep}.)
agree with the quasi-algebraic group associaetd to \(R^{i-1}f_{\perf,*}(W_n\Omega^1_{\log,X/k})\). 

\begin{corollary}\label{H2Zp1}
    Let $X$ be a proper smooth variety over an algebraically closed field $k$ of characteristic $p$. Then: 
    \[\HH^2(X,\Z_p(1))\cong \Z_p^{r}\oplus \NS(X)[p^\infty],\] 
    where $r$ is the multiplicity of slope $1$ in $\HH^2_{\cris}(X/W)$. 
\end{corollary}

\begin{proof} By Theorem \ref{log}, $\HH^2(X,\Z_p(1))\cong \HH^1(X,W\Omega^1_{X/k,\log})$. 
By Theorem \ref{IllusieRaynaud} and Remark \ref{EasyVan}, $\HH^1(X,W\Omega^1_{X/k,\log})$ 
is the $F$-invariant subgroup $\HH^1(X,W\Omega^1_{X/k})$, which, by Theorem \ref{SSbasic}, has rank $r$. 
By Proposition \ref{proinj}, the torsion submodule of $\HH^2(X,\Z_p(1))$ 
is isomorphic to the non-divisible part of $\HH^1(X,\mu_{p^\infty})$, which, via the Kummer sequence, is $\NS(X)[p^\infty]$. 
\end{proof}

\begin{theorem}\label{Mainthm}
    Let $X$ be a smooth proper variety over an algebraically closed characteristic $p$ field $k$. 
    Then the following statements hold. 

{\rm (1)} There is a split exact sequence 
    \begin{equation*}
        0\rightarrow (\Q_p/\Z_p)^{r-\rho} \rightarrow \Br(X)[p^{\infty}] 
        \rightarrow \HH^3(X,\Z_p(1))[p^{\infty}] \rightarrow 0,
    \end{equation*} 
    where $r=m_1(\HH^2_{\cris}(X/W))$ is the multiplicity of slope $1$ of $\HH^2_{\cris}(X/W)$. 
    
{\rm (2)} For each $i\geqslant 1$ there is an exact sequence 
    \begin{equation*}
        0 \rightarrow U_i(k)\rightarrow \HH^i(X,\Z_p(1))[p^{\infty}] \rightarrow J_i \rightarrow 0,
    \end{equation*} 
    where $U_i$ is a connected commutative unipotent algebraic $k$-group of dimension $T^{02}$, 
    and $J_i$ is the finite group
    \[(\text{\rm {c\oe ur}}(\HH^j(X,W\Omega^i_{X/k})))^{F=1}[p^\infty].\]  
\end{theorem}

\begin{proof}: Part (1) is a direct consequence of Section \ref{basiccalculation} and Proposition \ref{H2Zp1}. 
Part (2) follows from
Proposition \ref{log} and Theorem \ref{IllusieRaynaud}. 
\end{proof}

\medskip

Our next goal is to extract precise information about the domino numbers $T^{ij}$ 
from the geometry of $X$. The works of Crew and Ekedahl then showed that these domino numbers, 
in some sense, 
measures the discrepancy between the Newton and Hodge polygons of $X$--
although not always in a precise manner, 
as we will clarify. To make this notion more precise, we now recall the category \( D^b_c(R) \) (see \cite{I2}, \cite{Ek3}).

The category of graded (left) $R$-modules is abelian. 
Let $D^b(R)$ denote the bounded derived category of complexes of left $R$-modules, 
let $D^-(R)$ denote the derived category of complexes of left $R$-modules bounded from above. 
A complex of $R$-modules $M$ is said to be coherent if $M\in D^b(R)$ 
and all of its cohomology groups $\HH^i(M)$ are all coherent $R$-modules. 
We write $D_c^b(R)$ for the full subcategory of $D^b(R)$ consisting of coherent complexes.

\begin{definition}\label{shifting}
Given a complex $M$ of $R$-modules, we define $M[i](j)$ to be the complex obtained by shifting $M$ by $i$ degrees
in the complex grading and by $j$ degrees in the $R$-module degree. 
The differentials are defined as $d=(-1)^id$ or $d=(-1)^jd$, depending on whether the differential is 
taken in complex degree or in module degree, respectively.
\end{definition} 

Again, our conventions follow those of Ekedahl \cite{Ek1,Ek2,Ek3}. 

For each $n\geqslant 1$, define the ring $R_n:=R/(V^nR+dV^nR)$, 
then $R_n$ can viewed as a right $R$-module and a left module over the ring of dual numbers $W_n[d]/(d^2)$ 
(which, for simplicity, we always write $W_n[d]$). 
For any $R$-module (or sheaf of $R$-modules) $M$, the tensor product 
$R_n\otimes M=M/(V^nM+dV^nM)$ will be a left $W_n$-complex. 
This construction extends to a derived funtor: 
\[R_n\otimes^{\Ld}_R: D^-(R)\rightarrow D^-(W_n[d]).\]
For any $R$-module $M$, we define \[\widehat{M}:=R\varprojlim_n R_n\otimes^{\Ld}_R M.\]
There is a natural map $M\rightarrow \widehat{M}$, and we say $M$ is complete if this map is an isomorphism. 
We have the following identifications: 
\[W_n\Omega^\bullet_{X/k}\cong R_n\otimes_R W\Omega^\bullet_{X/k},\quad W_n\Omega^\bullet_{X/k}\cong R_n\otimes^{\Ld}_R W\Omega^\bullet_{X/k}.\]

Illusie and Ekedahl gave the following characterization of coherent complexes 
in \cite[Proposition 2.4.7]{I2}, explaining the term \emph{coherent $R$-module}. 

\begin{theorem}
    Let $M\in D^b(R)$, then the following conditions are equivalent:

{\rm (1)} $M\in D_c^b(R)$;

{\rm (2)} $M$ is complete and $R_n\otimes^{\Ld}_R M\in D^b_{c}(W_n[d])$ for all $n$, here $D^b_c(W_n[d])$ consisting of objects 
$N$ such that $\HH^i(N)$ is finitely generated over $W_n$ for all $i$; 

{\rm (3)} $M$ is complete and $R_1\otimes^{\Ld}_R M\in D^b_c(k[d])$. 
\end{theorem}

Since the operations $F,V,d$ on the de Rham-Witt complex satisfies the relations in the Raynaud ring $R$, 
the de Rham-Witt complex $W\Omega^\bullet_{X/k}$ can be regarded as a sheaf of $R$-modules on $X$. 
The global sections funtor $\Gamma$ then yields a graded $R$-module. 
Applying the derived functor $R\Gamma$, we obtain an object $R\Gamma(X,W\Omega^\bullet_{X/k})\in D(R)$, 
whose $j$-th cohomology $R^j\Gamma(X,W\Omega^\bullet_{X/k})$ is just the graded $R$-module associated to the complex (see \cite{I2}): 
\begin{equation*}
    \HH^j(X,WO_{X/k}) \xrightarrow{d} \HH^j(X,W\Omega^1_{X/k}) \xrightarrow{d} ... \xrightarrow{d} 
    \HH^j(X,W\Omega^n_{X/k}).
\end{equation*}

\section{Crew's Formula and Ekedahl's Inequality}

Crew, in \cite{Cr}, explored the relationship between the integers $T^{ij}$ 
and the Hodge numbers $h^{ij}$ by computing functor \(\Tor^R\)
for the graded $R$-module $\HH^j(X,W\Omega^\bullet_{X/k})$. The key idea involves the isomorphism
\[R_1\otimes^{\Ld}_R R\Gamma(X,W\Omega^\bullet_{X/k})\cong R\Gamma(X, R_1\otimes^{\Ld}_{R}W\Omega^\bullet_{X/k})\cong R\Gamma(X,\Omega^\bullet_{X/k}),\]
where both sides are viewed as objects in $D(k[d])$, here $k[d]:=k[d]/(d^2)$ is the ring of dual numbers. This leads to a spectral sequence
\[E_2^{ab}=\Tor^R_a(R_1,R\Gamma^b(X,W\Omega^\bullet_{X/k}))\Rightarrow R^{b-a}\Gamma(X,\Omega^\bullet_{X/k}),\] 
all terms are interpreted as graded $R$-modules, i.e, complexes of $R^0$-modules satisfying $FdV=d$. 

Fixing an integer $i$, consider the $i$-th term of all of these graded $R$-modules, 
consider the length function, we obtain the identity:  
\begin{equation}\label{lengtheq}
    \begin{aligned}
    &\sum_{a,b}(-1)^{b-a}\length_W(\Tor^R_a(R_1,R\Gamma^b(X,W\Omega^\bullet_{X/k}))^i)\\
    &=\sum_{j}(-1)^j\length_W(R^{j}\Gamma(X,\Omega^\bullet_{X/k})^i),
    \end{aligned}
\end{equation}
where the right-hand side equals the Euler characteristic $\chi(\Omega^i_{X/k})$. 
Here, for an $R$-module $M$, we regard $M$ as the complex \[...\rightarrow M^i\rightarrow M^{i+1}\rightarrow ....\] 
Since the length function of the functor \(\Tor^R\) is additive, 
to compute the left-hand side it suffices to compute the functor $\Tor^R$ for $R$-modules of type $I_a$, $I_b$, $II$ respectively. 
Following notations in \cite{Cr}, we define $\ell(M)(i)$ to be 
\[\ell(M)(i):=\sum_{j}(-1)^j\length_W\Tor^R_j(R_1,M)^i. \]

\begin{proposition} The following statements hold:

{\rm (1)} If $M$ is of type $I_a$ in degree $0$, then for every $i$, we have $\ell(M)(i)=0$. 

{\rm (2)} If $M$ is of type $I_b$ in degree $0$, then we have \[\ell(M)(0)=\length_W(M/VM),\qquad \ell(M)(1)=-\length_W(M/FM),\]
\[ \ell(M)(i)=0 \quad \text{for}\; i\neq 0,1.\]

{\rm (3)} If $M$ is one of the $U_t$, then we have 
\[\ell(M)(0)=\ell(M)(2)=1,\qquad \ell(M)(1)=2.\]
\[ \ell(M)(i)=0 \quad \text{for}\; i\neq 0,1,2.\]

\end{proposition}

\begin{proof}
    See \cite{Cr}. 
\end{proof}

For $M$ of type $I_b$, the quantities $\length_W(M/VM)$ and $\length_W(M/FM)$ can be computed from the slopes of $M$: 
\[\begin{aligned}
    \length _W(M/VM)&=\sum_{\lambda\in [0,1)}(1-\lambda)m_{\lambda}(M),\\
     \length _W(M/FM)&=\sum_{\lambda\in [0,1)}\lambda m_{\lambda}(M),
\end{aligned}\]
where $m_{\lambda}(M)$ denotes the multiplicity of the slope $\lambda$ in the (iso)crystal $M$. 
This naturally motivates the definition of the \emph{slope numbers} $m^{ij}$, introduced in \cite{Cr} 
and also discussed in \cite{I2}. 

\begin{definition}\label{mij}
    For a smooth proper variety $X/k$, \emph{the slope numbers} $m^{ij}$ are defined by:
\begin{equation}
    \begin{aligned}
    m^{ij}:=&\sum_{\lambda\in [0,1)}(1-\lambda)m_{\lambda}(\HH^j(X,W\Omega^i_{X/k})\otimes K)\\
    &+\sum_{\lambda\in [0,1)}\lambda m_{\lambda}(\HH^{j+1}(X,W\Omega^{i-1}_{X/k})\otimes K).
    \end{aligned}
\end{equation}
\end{definition}

These numbers have a geometric interpretation (see Proposition IV.2.5 \cite{Ek3}): 
for each $n$: form the \emph{Hodge-Newton polygon} 
by assigning slope $i$ multiplicity $m^{i,n-i}$, then this Hodge--Newton polygon is the
maximal convex polygon lying below the Newton polygon of $\HH^n_{\cris}(X/W)$,
with the same end points, integral slopes and integral breaking points.

To compare the slope numbers with the Hodge numbers, 
Crew and Ekedahl introduced the Hodge--Witt numbers $h_W^{ij}$. 

\begin{definition}
    The Hodge--Witt numbers are defined as: \[h_W^{ij}:=m^{ij}+T^{ij}-2T^{i-1,j+1}+T^{i-2,j+2}.\]
\end{definition}

A direct computation from equation (\ref{lengtheq}) yields: 

\begin{theorem}[Crew's formula]\cite[Theorem 4]{Cr}\label{CrewFor} For any $i$, we have 
    \begin{equation}
        \sum_{j}(-1)^j h^{ij}_W=\sum_j(-1)^j h^{ij}.
    \end{equation}
\end{theorem}

More strikingly, Ekedahl in \cite{Ek3} developed a theory of diagonal $t$-structure 
on the derived category of $R$-modules. Using this framework, he proved a fundamental inequality, 
revealing deep connections between Hodge-Witt and Hodge numbers.

\begin{theorem}[Ekedahl's inequality] \cite[Theorem IV 3.3]{Ek3}\label{Ekeineq} For any $i$ and $j$, we have $h_W^{ij}\leqslant h^{ij}$. 
\end{theorem}

This result provides new insight into the behavior of $h_W^{ij}$. 
Notably, a direct computation shows
\[\sum_{i+j=n}h^{ij}_W=\sum_{i+j=n}m^{ij}={\rm rank} \, \HH^n_{\cris}(X/W).\]
Therefore, if $\sum_{i+j=n} h^{ij}={\rm rank}\, \HH^n_{\cris}(X/W)$, 
then necessarily $h^{ij}_W=h^{ij}$ for all $i+j=n$. 

\begin{proposition}[Ekedahl]
    Suppose $\HH^n_{\cris}(X/W)$ and $\HH^{n+1}_{\cris}(X/W)$ are both torsion-free
    and the Fr\"olicher spectral sequence 
    \[E_1^{ij}=H^j(X,\Omega^i_{X/k})\Rightarrow H^{i+j}_{\dR}(X/k)\]
    degenerates at the $E_1$ page. Then for all $p,q$ such that $p+q=n$,
    \[h^{pq}=m^{pq}+T^{pq}-2T^{p-1,q+1}+T^{p-2,q+2}.\] 
    In particular, we have $T^{0n}=h^{0n}-m^{0n}$. 
\end{proposition}

\begin{proof} Under the stated assumptions \begin{equation*}
    \sum_{p+q=n}h^{pq}=\dim_k(\HH^n_{\dR}(X/k))={\rm rank}\,\HH^n_{\cris}(X/W)=\sum_{p+q=n}h_W^{pq}.
\end{equation*}
Thus in this case $h^{pq}_W=h^{pq}$ for all $p+q=n$, 
which gives $m^{0i}+T^{0i}=h^{0i}$. \end{proof}

This applies, for instance, to all abelian varieties. 
That is, for any abelian variety, we have \[h^{ij}=m^{ij}+T^{ij}-2T^{i-1,j+1}+T^{i-2,j+2}.\]
This relation allows for the computation of all $T^{ij}$ directly from the 
slopes of an abelian variety (\cite[Corollary 6.3.12]{I2}). 

To illustrate, consider an abelian 3-fold $X$ whose slopes of $\HH^1_{\cris}(X/W)$ 
are $\{0,\frac{1}{2},\frac{1}{2},\frac{1}{2},\frac{1}{2},1\}$. 
Such an example can be obtained as the product of an ordinary elliptic curve and a supersingular abelian surface. 
The Newton polygon of $\HH^1_{\cris}(X/W)$ is pictured in orange below: 
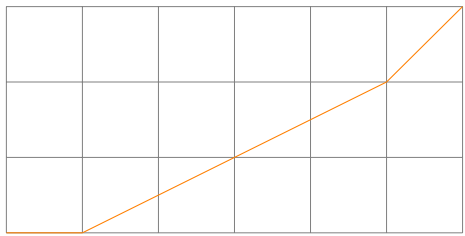
\begin{center}
    \begin{tikzpicture}\label{almostsupersingular}
        \draw [step=1cm,gray,very thin] (0,0) grid (6,3);
        \draw [orange] (0,0) -- (1,0);
        \draw [orange] (1,0) -- (5,2);
        \draw [orange](5,2) -- (6,3);
    \end{tikzpicture}
\end{center}
We draw the Newton polygon of $\HH^2_{\cris}(X/W)$ also in orange, 
the Hodge--Newton polygon in blue, 
and the Hodge (or Hodge--Witt) polygon in black, for comparison. 
The Hodge--Newton polygon is visibly the uppermost convex polygon below the Newton polygon, 
with integral slopes and break points. From this, we see that $T^{02}=1$. 

\begin{center}
    \begin{tikzpicture}\label{almostsupersingular2}
        \draw [step=0.5cm,gray,very thin] (0,0) grid (7.5,7.5);
        \draw [orange] (0,0) -- (2,1);
        \draw [orange] (2,1) -- (5.5,4.5);
        \draw [orange] (5.5,4.5) -- (7.5,7.5);
        \draw [blue] (0,0)--(1,0);
        \draw [blue] (1,0)--(6.5,5.5);
        \draw [blue] (6.5,5.5)--(7.5,7.5);
        \draw [black] (0,0)--(1.5,0);
        \draw [black] (1.5,0)--(6,4.5);
        \draw [black] (6,4.5)--(7.5,7.5);
    \end{tikzpicture}
\end{center}
\medskip

These mysterious Hodge-Witt numbers \( h_W^{ij} \) appear 
to be the true analogues of Hodge numbers for proper smooth varieties in characteristic \( p \), 
since the actual Hodge numbers \( h^{ij} \) often exhibit poor behavior. 
Ekedahl showed that the Hodge-Witt polygon of degree \( n \)--
the polygon obtained by assigning slope \( i \) multiplicity \( h_W^{i,n-i} \)--
lies below the so-called \emph{abstract Hodge polygon} of degree \( n \) (see \cite{Ek3}).
Concretely, consider the crystalline cohomology group \( \HH^n_{\mathrm{cris}}(X/W) \), and define
\[
M := \HH^n_{\mathrm{cris}}(X/W) / \HH^n_{\mathrm{cris}}(X/W)_{\mathrm{tors}}.
\]
Taking the elementary divisor decomposition of \( M/FM \),
\[
M/FM \cong \bigoplus_i (W/p^iW)^{h_{\mathrm{abs}}^{i,n-i}},
\]
then the Hodge-Witt polygon lies below the Hodge polygon defined by these \( h_{\mathrm{abs}}^{i,n-i} \).

Nevertheless, many questions remain unanswered. 
For instance, it is still unknown whether the Hodge-Witt numbers satisfy Hodge symmetry or not. 
A deeper question concerns the mixed characteristic setting: suppose \( \mathcal{X}/O_K \) 
is a proper smooth formal scheme over a \( p \)-adic ring \( O_K \). 
What are the relationships among the abstract Hodge polygon, the Hodge-Witt polygon of the special fiber, 
and the Hodge polygon of the generic fiber? We can conjecture that
\[
\text{abstract Hodge} \geqslant \text{generic Hodge} \geqslant \text{Hodge-Witt},
\]
but no part of this inequality has been proven yet.

\chapter{First Applications}\label{Chap2}

In the context of Theorem \ref{Mainthm} and the theory of Hodge-Witt numbers, 
this chapter examines various examples to illustrate the theory.  
Among these, ordinary varieties form the most tractable class in characteristic $p$.

\section{Ordinary Varieties}

A proper smooth variety is called {\em ordinary} if, for all $j$, the cohomology groups
of the subsheaf of exact forms vanish: $\HH^i(X,B\Omega^j_{X/k})=0$ for all $i$. 
Here, $B\Omega^j$ denotes the image of differential map $d:\Omega^{j-1}\rightarrow \Omega^j$. 

\begin{theorem}[Illusie]\label{Illusie}
Let $X$ be an ordinary smooth proper variety over $k$. Then
the slope spectral sequence split degenerates at the $E_1$ page. 
\end{theorem}
\begin{proof} By \cite[Theorem IV.4.13]{IR}, 
if $X$ is ordinary, then $\HH^q(X,BW\Omega^p_{X/k})=0$ for all $q$.
By \cite[Theorem 3.4.1]{I2}, $X$ is Hodge--Witt in the sense that $\HH^q(X,W\Omega_X^p)$ 
is a finitely generated $W$-module for all $q$. 
This implies that the slope spectral sequence for $X$ split degenerates at the $E_1$ page
by \ref{HodgeWittn}, \cite[Theorem II.3.7]{Il1}. 
\end{proof} 

As a consequence, if $X$ is ordinary, then all cohomology groups $\HH^q(X,W\Omega^p)$ are finitely generated $W$-modules, 
and the crystalline cohomology group $\HH^i_{\cris}(X/W)$ decomposes as the direct sum of $\HH^q(X,W\Omega_{X/k}^p)$, where $p+q=i$. 

\begin{corollary}\label{ordcor}
    Let $X$ be an ordinary smooth proper variety over $k$. Then for all $i,j$
we have $T^{ij}=0$
    and $\text{\rm c\oe ur}(\HH^j(X,W\Omega^i_{X/k})) \cong \HH^j(X,W\Omega^i_{X/k})$. 
\end{corollary}
\begin{proof}
Since the slope spectral sequence degenerates at the $E_1$ page, 
all differentials vanish. 
By definition, we have $F^{\infty}B_2^{ij}=0$ and  $V^{-\infty}Z_2^{ij}=\HH^j(X,W\Omega^i_{X/k})$. 
This implies \[T^{ij}=\dim_k \HH^j(X,W\Omega^i_{X/k})/(V^{-\infty}Z_2^{ij}+V\HH^j(X,W\Omega^i_{X/k}))=0,\]
and $\text{c\oe ur}(\HH^j(X,W\Omega^i_{X/k})) \cong \HH^j(X,W\Omega^i_{X/k})$ for all $i,j$.
\end{proof} 

\begin{theorem}[Bloch--Kato]\label{Bloch-Kato}
Let $X$ be an ordinary smooth proper variety over $k$.
Then the following statements hold.
    
{\rm (1)} $\HH^q(X,W_n\Omega_{X/k,\log}^r) \otimes_{\Z/p^n} W_n \cong \HH^q(X,W_n\Omega_{X/k}^r)$ for all $q,r$ and $n$. 

{\rm (2)} $\HH^q(X,W\Omega_{X/k,\log}^r) \otimes_{\Z_p} W \cong \HH^q(X,W\Omega_{X/k}^r)$ for all $q$. 

{\rm (3)} $F:\HH^q(X,W\Omega_{X/k}^r)\rightarrow \HH^q(X,W\Omega_{X/k}^r)$ is bijective for all $q$ and $r$.
\end{theorem}
\begin{proof}
    This is \cite[Proposition 7.3]{BK}.
\end{proof}  

\medskip

\begin{corollary}\label{slopeint}
    Let $X$ be an ordinary smooth proper variety over $k$. 
    Then all slopes of its crystalline cohomology are integers. 
    More precisely, the multiplicity of slope $i$ of $\HH^n_{\cris}(X/W)$ 
    equals the rank of $\HH^{n-i}(X,W\Omega_{X/k}^i)$. 
\end{corollary}
\begin{proof}
    This follows directly from Theorem \ref{Illusie} and Theorem \ref{Bloch-Kato} (2). 
    The isomorphism $\HH^q(X,W\Omega_{X/k,\log}^1) \otimes_{\Z_p} W \cong \HH^q(X,W\Omega_{X/k}^1)$
    implies that $\HH^q(X,W\Omega_{X/k}^1)$ is pure of slope $0$. Since
    $\HH^i_{\cris}(X/W)$ is isomorphic to the direct sum of 
    (Frobenius twists) of $\HH^q(X,W\Omega_{X/k}^p)$ for $p+q=i$, the corollary follows.
\end{proof} 

Note that the converse does not hold: having integral slopes in crystalline cohomology does not imply ordinarity. In fact, there exist varieties with integral slopes that are not even Hodge-Witt.

Now apply Theorem \ref{Mainthm}, we obtain

\begin{proposition}\label{Ordinary}
    Let $X$ be an ordinary smooth proper variety over $k$. Then we have
    $\Br(X)[p^\infty]\cong (\Q_p/\Z_p)^{t-\rho}\oplus J$,
    where $t$ is the rank of the $W$-module 
$\HH^1(X,W\Omega_{X/k}^1)$ and $J$ is a finite $p$-group such that 
    $J\otimes W \cong \HH^2(X,W\Omega_{X/k}^1)[p^{\infty}]$.
\end{proposition}
\begin{proof}
By Theorem \ref{Mainthm} (1), $\Br(X)[p^\infty]$ has the desired form, 
where $t$ is the multiplicity of slope $1$ in $\HH^2_{\cris}(X/W)$,
which equals the rank of $\HH^1(X,W\Omega^1)$ by Corollary \ref{slopeint}.
The finite exponent group is an extension of a finite group $J$ by $U(k)$. 
By Corollary \ref{ordcor}, the unipotent dimension $T^{02}=0$, hence $U(k)=0$,
and we have
$$J\cong (\text{c\oe ur}(\HH^2(X,W\Omega^1_{X/k})))^{F=1}[p^{\infty}] \cong (\HH^2(X,W\Omega^1_{X/k}))^{F=1}[p^{\infty}].$$ 
Now Theorem \ref{Bloch-Kato} (2) and Theorem \ref{IllusieRaynaud} (1)  tell us that  \[\HH^2(X,W\Omega^1_{X/k})\cong (\HH^2(X,W\Omega^1_{X/k}))^{F=1}\otimes_{\Z_p}W.\]
Taking the $p$-primary torsion subgroups, we obtain \[J\otimes W \cong \HH^2(X,W\Omega^1_{X/k})[p^{\infty}].\] 
\end{proof}

\begin{corollary}
Let $X$ be an ordinary smooth proper variety over $k$.
If $\HH^3_{\cris}(X/W)$ is torsion-free, 
    then $\Br(X)[p^\infty] \cong (\Q_p/\Z_p)^{t-\rho}$,
    where $t$ is the rank of the $W$-module $\HH^1(X,W\Omega^1_{X/k}).$
\end{corollary}
\begin{proof}
    In this case, $\HH^2(X,W\Omega^1_{X/k})$ is a direct summand of 
$\HH^3_{\cris}(X/W)$, which is a free $W$-module. 
So $J=0$ by Proposition \ref{Ordinary}. 
\end{proof} 
\begin{proposition}
   Let $X$ be an ordinary smooth proper variety over $k$. Then we have
$\NS(X)[p^{\infty}]\otimes W\cong \HH^1(X,W\Omega^1_{X/k})[p^{\infty}]$.
\end{proposition}
\begin{proof} By Kummer sequence, $\HH^1(X,\mu_{p^\infty})\cong \Pic(X)[p^\infty]$. 
By Theorem \ref{proinj}, Theorem \ref{log} and Theorem \ref{IllusieRaynaud},
we have $\HH^1(X,\mu_{p^\infty})\cong (\Q_p/\Z_p)^{n}\oplus T$, 
where $n$ is the rank of $\HH^0(X,W\Omega^1_{X/k})$ and $T=\HH^2(X,\Z_p(1))[p^{\infty}]$. 
This gives an isomorphism $\NS(X)[p^{\infty}]\cong \HH^2(X,\Z_p(1))[p^{\infty}]$. 
By Theorem \ref{Mainthm}(2), $\HH^2(X,\Z_p(1))[p^{\infty}]$ 
is an extension of $(\text{c\oe ur}(\HH^1(X,W\Omega^1_{X/k})))^{F=1}[p^{\infty}]$ by $U(k)$, 
where $U$ is of dimension $T^{01}$. But $T^{01}=0$ holds for all proper smooth varieties, 
and $\text{c\oe ur}(\HH^1(X,W\Omega^1_{X/k}))$ is isomorphic to $\HH^1(X,W\Omega^1_{X/k})$,  
so $T=\HH^2(X,\Z_p(1))[p^{\infty}]$ is isomorphic to $\HH^1(X,W\Omega^1_{X/k})^{F=1}[p^{\infty}]$. 
Again Theorem \ref{Bloch-Kato} (2) and Theorem \ref{IllusieRaynaud} (1) 
tell us that  \[\HH^1(X,W\Omega^1_{X/k})\cong (\HH^1(X,W\Omega^1_{X/k}))^{F=1}\otimes_{\Z_p}W.\]
Taking $p$-primary torsion subgroups, 
we obtain $T\otimes W \cong \HH^1(X,W\Omega^1_{X/k})[p^{\infty}]$. 
\end{proof} 

\medskip

\section{Interlude: Representability and Flat Duality Theorem}\label{Secrep}

To proceed further we now review the representability theorem and Milne's flat duality theorem for surfaces, which we will repeatly use. 
The first result concerns the higher direct images of pushforward for the big fppf site.

\begin{theorem}[Artin, Bragg, Olsson] \label{rep} 
    Let $f:X\rightarrow Spec(k)$ be a projective finite type scheme over $k$.
    Let $G$ be a finite flat abelian group scheme over $X$. Let $R^if_*G$ be the $i$-th 
    fppf-derived pushforward of the sheaf represented by $G$. 
    Then $R^if_*G$ is representable by a finite type group scheme over $k$. 
\end{theorem}
\begin{proof} This is proved by Bragg in \cite[Corollary 1.6]{BO}. \end{proof}

\medskip

Berthelot showed that the group of $k$-points of $R^if_*\mu_{p^n}$ is 
exactly the fppf cohomology group $\HH^i_{\fppf}(X,\mu_{p^n})$, see \cite[Corollary 2.10]{Be}. 
As a consequence of this, we may consider the quasi-algebraic group, (see
\cite{Se}), associated to $R^if_*\mu_{p^n}$. 
Recall that, a quasi-algebraic group is an equivalent class of algebraic groups, where two algebraic groups $G_1$ and $G_2$ are 
regarded as equivalent
if and only if there is a homomorphism $G_1\rightarrow  G_2$ such that on the value of $k$-points, it induces an bijection $G_1(k)\rightarrow G_2(k)$. 
As a result, the group of $k$-points only depends on this equivalence class. 
The quasi-algebraic group associated to $R^if_*\mu_{p^n}$ 
is denoted by $\uH^i(X,\mu_{p^n})$. In this way, $\HH^i_{\fppf}(X,\mu_{p^n})$ is identified with $\uH^i(X,\mu_{p^n})(k)$, the group of k-points of $\uH^i(X,\mu_{p^n})$. 
The advantage is that the category of quasi-algebraic groups over $k$ is abelian. 
Another way is to consider the associated perfect algebraic groups, 
and view them as sheaves on the perfect cite 
(see for example, \cite[Section 2.7]{Be}, or \cite[p.~192]{IR}). 


Since $R^if_*\mu_{p^n}$ is a sheaf annihilated by a power of $p$, by the classification theorem \cite[Proposition 7, p.~10]{Se},
the quasi-algebraic group $\uH ^i(X,\mu_{p^n})$ is an extension of an \'etale group
by a connected unipotent group: 
\begin{equation}
0\rightarrow \uU ^i(X,\mu_{p^n}) \rightarrow \uH ^i(X,\mu_{p^n}) \rightarrow \uD ^i(X,\mu_{p^n}) \rightarrow 0
\end{equation}
Since $k$ is algebraically closed, the \'etale part is just an abelian group of $p$-power order, 
and the connected part is an iterated extension of the additive $k$-group $\G_a$. 

The flat duality is the following statement.
\begin{theorem}[Artin, Milne, Berthelot \cite{Be}] \label{FlatDuality}
     Let $X$ be a smooth proper surface over $k$.
For $i,n\geqslant 0$, there are isomorphisms
\[\begin{aligned}
\uD^i(X,\mu_{p^n})\cong \mathrm{Hom}(\uD^{4-i}(X,\mu_{p^n}),\Q_p/\Z_p),\\
    \uU^i(X,\mu_{p^n}) \cong \mathrm{Ext}^1(\uU^{5-i}(X,\mu_{p^n}),\Q_p/\Z_p).
\end{aligned}
\]
\end{theorem}

Here $\mathrm{Ext}^1$ is taken in the category of quasi-algebraic group schemes over $k$. 
It has the following property: for each connected unipotent quasi-algebraic group $U$ over $k$, 
we can equip the abelian group $\mathrm{Ext}^1(U,\Q_p/\Z_p)$ with a \emph{canonical} connected unipotent quasi-algebraic group structure
(\cite[Section 8.4, p. 391]{Se}). For example, we have $\mathrm{Ext}^1(\G_a, \Q_p/\Z_p)\cong k$, and $\mathrm{Ext}^1(\W_n,\Q_p/\Z_p)\cong W_n$. 

We remark that the unipotent torsion of the logarithmic Hodge-Witt cohomology groups 
$\HH^q(X,W_n\Omega^p_{X/k,\log})$ for a proper smooth variety of any dimension also satisfy an analogous duality, 
due to Milne. 

\section{Surfaces}

\begin{theorem}\label{surfaces}
    Let $X$ be a proper smooth surface over $k$. Then the following statements holds: 

 {\rm   (1)} $J$ is the Pontryagin dual of $\NS(X)[p^\infty]$;

 {\rm   (2)} we have Crew's formula for surfaces: $\dim(U)=(h^{02}-h^{01})-(m^{02}-m^{01}).$
\end{theorem}

\begin{proof}
When $X$ is a surface, Proposition \ref{EasyVan} implies that all $T^{ij}$ are zero except
possibly $T^{02}$. Then
(\ref{CrewFor}) for $i=0$ gives Crew's formula for surfaces, proving
Theorem \ref{surfaces} (2). 

Let us prove Theorem \ref{surfaces} (1). By Theorem \ref{Mainthm} (1), (2), (3), 
it remains to prove that the finite group $J$ is the Pontryagin dual of $\NS(X)[p^{\infty}]$. 
By Theorem \ref{FlatDuality}, the \'etale part of $\HH^3_{\fppf}(X,\mu_{p^n})$ is the Pontryagin dual of the \'etale part of $\HH^1_{\fppf}(X,\mu_{p^n})$. 
Consider the isomorphism 
\[\HH^1_{\fppf}(X,\mu_{p^n})\cong \Pic(X)[p^n]\]
induced by the Kummer sequence. For  $n$ large enough so that $p^n$ annihilates $\NS(X)[p^\infty]$, 
the group $\Pic(X)[p^n]$ is isomorphic to 
$(\Z/p^n\Z)^g\oplus \NS(X)[p^{\infty}]$, 
where $g$ is the $p$-rank of the Picard variety of $X$. 
As a result, for $n>> 0$, the \'etale part of $\HH^3_{\fppf}(X,\mu_{p^n})$ is of the form 
\[(\Z/p^n\Z)^g\oplus \NS(X)[p^{\infty}]^{\vee}.\]
Next, the $k$-vector space $\uU^4(X,\mu_{p})$ is the dual (Theorem \ref{FlatDuality}) of $\uU^1(X,\mu_{p})$, 
so it must be trivial. This implies $\uU^4(X,\mu_{p^n})=0$, by induction on $n$.
Finally, $\uD^4(X,\mu_{p^n})$ is the Pontryagin dual of $\uD^0(X,\mu_{p^n})$, 
so it must also be trivial. 
Thus for surfaces we have $\HH^4_{\fppf}(X,\mu_{p^n})=0$, and $\HH^4(X,\Z_p(1))=0$. 
Now consider the exact sequence 
\[...\rightarrow \HH^3(X,\Z_p(1))\xrightarrow{p^n}\HH^3(X,\Z_p(1))\rightarrow \HH^3_{\fppf}(X,\mu_{p^n})\rightarrow \HH^4(X,\Z_p(1))\rightarrow ...\]
we have
$\HH^3_{\fppf}(X,\mu_{p^n})\cong\HH^3(X,\Z_p(1))/p^n$. 
For any $n$ large enough so that $p^n$ annihilates $\HH^3(X,\Z_p(1))[p^{\infty}]$, 
consider the snake lemma of the multiplication by $p^n$ map on the (split) short exact sequence 
\[0\rightarrow \HH^3(X,\Z_p(1))[p^\infty]\rightarrow \HH^3(X,\Z_p(1))\rightarrow \Z_p^g\rightarrow 0,\]
we get \[ \HH^3(X,\Z_p(1))/p^n\cong (\Z/p^n\Z)^{g}\oplus \HH^3(X,\Z_p(1))[p^{\infty}],\]  
the \'etale part of which is the direct sum of $(\Z/p^n\Z)^{g}$ and $J$. 
This group is also $(\Z/p^n\Z)^{g}\oplus \NS(X)[p^\infty]^\vee$, so $J$ is isomorphic to the Pontryagin dual of $\NS(X)[p^{\infty}]$. This proves this theorem. 
\end{proof}

\medskip

Let $X$ now be a surface with $q=-p_a$, 
where $q=\dim \Pic_{X/k}$ is the irregularity and $p_a=h^{02}-h^{01}$ is the arithmetic genus. 
In \cite[Corollary 5]{Su}, Suwa showed that the slope spectral sequence 
$\HH^j(X,W\Omega_{X/k}^i)\Rightarrow \HH^{i+j}_{\cris}(X/W)$ (split) degenerates at $E_1$ page.

\begin{corollary}\label{qpa}
    Let $X$ be a smooth proper surface over $k$ such that $q=-p_a$, then $U=0$. 
\end{corollary}

\begin{proof}
The degeneration of the slope spectral sequence implies $T^{02}=0$ hence $U=0$. 
\end{proof} 

We give a few examples of surfaces with $q=-p_a$ and compute the relavent invariants, see also \cite{Su}. 

(1) Rational or ruled surfaces.  

(2) Enriques surface over an algebraically closed field
$k$ of characteristic $p$ (\cite[Section 3-5]{BM3}). In this case the Picard number is $\rho(X)=10$.
There are three types of Enriques surfaces: classical, singular, and supersingular,
the last two types occuring only for $p=2$.
The Picard scheme ${\rm Pic}_{X/k}$ isomorphic to
the product of $\Z^{\oplus 10}$ and a finite flat group scheme
$\Z/2\Z$, $\mu_2$, $\alpha_2$, when $X$ is classical, singular, and supersingular, respectively. 


\begin{corollary}\label{Enriques}
    Let $X$ be an Enriques surface over $k$. Then $r=\rho$ and $U=0$. Moreover, 
the following statements hold.

{\rm    (1)} If $p\neq 2$, then $\Br(X)[p^{\infty}]=0$.

{\rm    (2)} If $p=2$ and $X$ is a classical Enriques surface, then $\Br(X)[p^{\infty}]=\Z/2$.

{\rm    (3)} If $p=2$ and $X$ is a non-classical Enriques surface, then $\Br(X)[p^{\infty}]=0$.
\end{corollary}
\begin{proof}
    For Enriques surfaces the number $T^{02}=0$ by Theorem \ref{qpa}. We have $r=\rho=10$. 
    Then $\NS(X)[p^{\infty}]$ is non-trivial only when $p=2$ and when $X$ is classical, in which case $\NS(X)[p^{\infty}]\cong \Z/2\Z$. 
    Now the result follows from Theorem \ref{surfaces}. 
\end{proof}

\medskip

(3) Hyperelliptic or quasi-hyperelliptic surfaces \cite[Section 3]{BM2} \cite[Section 2]{BM3}. 
In the Enriques classification of surfaces, they are defined to be a special kind of surfaces with $\kappa=0$, 
other than K3 surfaces, abelian surfaces, and Enriques surfaces. 
They are surfaces with $q=-p_a=1$, hence $\Alb(X)$ is an elliptic curve, moreover 
the fibres of the canonical map \[\pi:X\rightarrow \mathrm{Alb}(X)\]
are either all non-singular elliptic or all rational with one cusp. 
They are called hyperelliptic and quasi-hyperelliptic according to whether the fibres are all non-singular elliptic or all rational with one cusp. 
The authors in \cite{BM2} and \cite{BM3} proved that, 
hyperelliptic surfaces are always a quotient of a product of two elliptic curves by a finite group scheme; 
and quasi-hyperelliptic surfaces are always a quotient of $E\times C_0$, 
where $C_0$ is rational with one cusp, by a finite group scheme. 

\begin{corollary}\label{hyperelliptic}
    Let $p=2$ or $3$, let $X/k$ be a hyperelliptic or quasi-hyperelliptic surface. Then $U=0$ and $J=\NS(X)[p^\infty]^\vee$. 
\end{corollary}

An example of a hyperellptic surface when $p=2$ is an Igusa surface 
$I:=(E_1\times E_2)/(\Z/2\Z)$ where $E_1$ is an ordinary elliptic curve, $E_2$ is a supersingular elliptic curve, 
and $\Z/2\Z\cong E_1[p](k)$ acts on $E_1\times E_2$ via $\alpha(x,y)=(x+\alpha,-y)$. In this case, $U=J=0$, see Proposition \ref{Igusa}. 

(4) Godeaux surfaces. In this case, $q=-p_a=0$, and $r=\rho=9$ \cite{La}. Similar to the Enriques case, 
there are two types of Godeaux surfaces: classical and singular,
where singular ones occuring only for $p=5$.
The Picard scheme ${\rm Pic}_{X/k}$ isomorphic to
the product of $\Z^{\oplus 9}$ and a finite flat group scheme
$\Z/5\Z$ or $\mu_5$, when $X$ is classical and singular respectively. 

\begin{corollary}
    Let $X/k$ be a Godeaux surface. Then $r=\rho=9$ and $U=0$. Moreover, 
    the following statements hold. 

{\rm    (1)} If $p\neq 5$, then $\Br(X)[p^{\infty}]=0$.

{\rm    (2)} If $p=5$ and $X$ is classical, then $\Br(X)[p^{\infty}]=\Z/5$.

{\rm    (3)} If $p=5$ and $X$ is non-classical, then $\Br(X)[p^{\infty}]=0$.
\end{corollary}

\begin{proof}
     For Godeaux surfaces the number $T^{02}=0$ by Theorem \ref{qpa}. We have $r=\rho=9$. 
Then $\NS(X)[p^{\infty}]$ is non-trivial only when $p=5$ and when $X$ is classical, in which case $\NS(X)[p^{\infty}]\cong \Z/5\Z$. 
Now the result follows from Theorem \ref{surfaces}.
\end{proof} 

\medskip

Another easy consequence is that $h^{02}\leqslant 1$ implies $T^{02}\leqslant 1$, 
so $U(k)$ is either trivial or isomorphic to $k$. 
If, moreover, $\Pic_{X/k}$ is smooth, then the formal Brauer group $\widehat{\Br(X)}$ is a smooth formal group of dimension $1$, so it either has finite height or is $\widehat{\G_a}$. 

\begin{corollary}\label{NcSfs}
    Let $X$ be a proper smooth surface over $k$ such that the Picard scheme of $X$ is smooth, and 
    $h^{2,0}\leqslant 1$. Let $h$ be the height of $\widehat{\Br(X)}$ (we set $h=0$ if $h^{2,0}=0$).
    Then we have 
    \[\Br(X)[p^\infty]\cong\begin{cases}
        (\Q_p/\Z_p)^{b_2-2h-\rho} \oplus \NS(X)[p^{\infty}]^\vee, \ \text{if}\ h<\infty\\
        (\Q_p/\Z_p)^{b_2-\rho} \oplus H, \ \text{if}\ h=\infty
    \end{cases}\]
    where $H$ is an extension of the Pontryagin dual of $\NS(X)[p^{\infty}]$ by $k$. 
\end{corollary}

This corollary can be applied to abelian surfaces and K3 surfaces. 
For K3 surfaces, Corollary \ref{NcSfs} was proved in \cite{BY}.
Note that the Tate conjecture for K3 surfaces implies that for the supersingular K3 surfaces we have $\rho=b_2$, which equals $22$.
Corollary \ref{NcSfs} can also be applied to classical Enriques surfaces when $p\neq 2$.

\section{Abelian Varieties}


In this section we study $\Br(X)[p^{\infty}]$ when $X$ is an abelian variety. We start by showing that the finite group 
$J:=\text{c\oe ur}(\HH^2(X,W\Omega^1_{X/k}))^{F=1}[p^{\infty}]$ 
in Theorem \ref{Mainthm} (2) is trivial when $X$ is an abelian variety. 

\begin{proposition}\label{AbeT}
    Suppose $X/k$ is a proper smooth variety over an algebraically closed field $k$ of characteristic $p$, 
    such that $\HH^n_{\cris}(X/W)$ is free. Then for any $i+j=n$, 
    $\text{\rm{c\oe ur}}(\HH^j(X,W\Omega^i_{X/k}))^{F=1}$
    is a free $\Z_p$-module. 
    The logarithmic de Rham-Witt cohomology $\HH^j(X,W\Omega^i_{\log,X/k})$ 
    is then an extension of the free $\Z_p$-module $\text{\rm{c\oe ur}}(\HH^j(X,W\Omega^i_{X/k}))^{F=1}$ 
    by the group of $k$-points of a connected commutative unipotent group $U$ of dimension $T^{i-1,j}$.
\end{proposition}
\begin{proof}: By \cite[Theorem II.3.4]{IR} and Proposition \ref{Survivesofcoeur}, the heart $\text{c\oe ur}(\HH^j(X,W\Omega^i_{X/k}))$ 
survives to the $E_{\infty}$-page: 
\begin{equation*}
    B_{\infty}^{ij}\subseteq F^{\infty}B_2^{ij} \subseteq V^{-\infty}Z_2^{ij}\subseteq Z_{\infty}^{ij},
\end{equation*} and the quotient $ F^{\infty}B_2^{ij}/B_{\infty}^{ij}$ and $Z_{\infty}^{ij}/ V^{-\infty}Z_2^{ij}$ are finite length $W$-modules. 
We need to prove that $F$ cannot have non-zero
fixed points on the torsion part of $\text{c\oe ur}(\HH^j(X,W\Omega^i_{X/k}))$.  

To see this, consider the filtration on $\HH^n_{\cris}(X/W)$ induced by the slope spectral sequence: 
\[0\subseteq \Fil^0\subseteq \Fil^1 \subseteq ...\Fil^{j-1} \subseteq \Fil^j \subseteq ...\subseteq \Fil^n \cong \HH^n_{\cris}(X/W).\]
Each $\Fil^i$ is a free $W$-module, and we have 
$\Fil^j/\Fil^{j-1}\cong E_{\infty}^{ij}\cong Z_{\infty}^{ij}/B_{\infty}^{ij}$. 
Write $\mathcal{F}$ for the Frobenius acting on $\HH^n_{\cris}(X/W)$.
Then $\frac{\mathcal{F}}{p^i}$ acts on $\Fil^j/\Fil^{j-1}$, and this action is induced by 
the action of $F$ on $\HH^j(X,W\Omega^i_{X/k})$. 
By this equality, we can find $W$-submodules 
$\Fil^{j-1} \subseteq B'\subseteq Z'\subseteq \Fil^j$ such that 
$B'/\Fil^{j-1}=F^{\infty}B_2^{ij}/B_{\infty}^{ij}$ and $Z'/\Fil^{j-1}=V^{-\infty}Z_2^{ij}/B_{\infty}^{ij}$. 
It suffices to prove that $\frac{\mathcal{F}}{p^i}$ cannot have non-zero fixed points on the torsion part of $Z'/B'$. 

Since the quotient $B'/\Fil^{j-1}$ is of finite length, 
the crystals $(B',\mathcal{F})$ and $(\Fil^{j-1},\mathcal{F})$ have the same slopes. 
By Illusie's fundamental results in \cite{Il1}, 
the slopes of $(\Fil^{j-1},\mathcal{F})$ 
can be identified with the slopes of $(\HH_{\cris}^n(X/W),\mathcal{F})$ that are
contained in the interval $[j,\infty)$. 
Now suppose that $[x]\in Z'/B'$ is a torsion element fixed by $\frac{\mathcal{F}}{p^i}$, 
where $x\in Z'$. Then $p^kx\in B'$ for some $k$. 
Since the slopes of $F=\frac{\mathcal{F}}{p^i}$ on $B'$ are $\geqslant 1$, 
for any integer $m$ we can write $F^m(p^kx)=p^m(y)$ for some $y\in B'$. 
Choose $m\geqslant k$, then $F^m(p^kx)=p^my$ implies $F^m(x)=p^{m-k}y\in B'$
since $B'$ is torsion-free.
But $Fx-x\in B'$ implies $F^{m}x-x\in B'$, so $x$ must be contained in $B'$. 
Hence $[x]=0$ in $Z'/B'$. \end{proof}
\medskip

\begin{theorem}\label{Abe}
    Let $A$ be an abelian variety of dimension $g$ over $k$. Then we have \[\Br(A)\cong (\Q_p/\Z_p)^{r-\rho}\oplus U,\]
    where the dimension of the connected unipotent group $U$ is $\dim(U)=h^{02}-m^{02}$, where $h^{02}=\frac{g(g-1)}{2}$, and $J=0$. 
\end{theorem}

\begin{proof}
    Combine Theorem \ref{Mainthm} (1), (2), (3) with Proposition \ref{AbeT}. 
\end{proof}

\medskip

As an illustration,
for an abelian 3-fold $A$ over $k$ we calculate all possible values of the unipotent dimension $T^{0,2}$.
There are five different possibilities for the Newton polygon of $\HH^1_{\cris}(A/W)$, 
namely: the ordinary case (black), the almost ordinary case (blue), the almost supersingular case (orange), the $1/3$ case (green), and the supersingular case (red). 
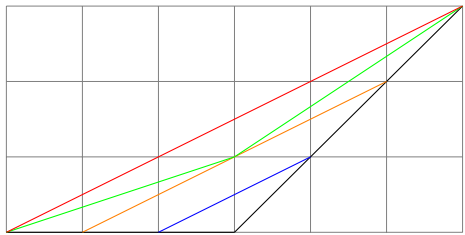
\begin{center}
\begin{tikzpicture}\label{NewtonPolygon}
    \draw [step=1cm,gray,very thin] (0,0) grid (6,3);
    \draw (0,0) -- (3,0);
    \draw (3,0) -- (6,3);
    \draw [blue] (2,0) -- (4,1);
    \draw [orange] (1,0) -- (5,2);
    \draw [green](0,0) -- (3,1);
    \draw [green](3,1) -- (6,3);
    \draw [red](0,0) -- (6,3);
\end{tikzpicture}
\end{center}
Then, by Theorem \ref{Mainthm}(3), we have $h^{02}=h_W^{02}=m^{02}+T^{02}$.
We have $h^{02}=3$. Explicitly calculating $m^{02}$ we find
\begin{equation*}
    m^{02}=\sum_{\lambda\in [0,1)}(1-\lambda)m_{\lambda}(\HH^2(A,WO_X)\otimes K)=\sum_{\lambda\in [0,1)}(1-\lambda)m_{\lambda}(\HH^2_{\cris}(A/W)).
\end{equation*}
We have the following table: 
\begin{center}
    \begin{tabular}{ c|c|c|c|c } 
     cases                & slopes of $\HH^2(A,WO_{X})\otimes K$                           & $m^{02}$ & $T^{02}$ & $T^{03}$\\ 
     ordinary             & $\{0,0,0\}$                                                    &  $3$     &   $0$    &  $0$  \\ 
     almost ordinary      & $\{0,\frac{1}{2},\frac{1}{2},\frac{1}{2},\frac{1}{2}\}$        &  $3$     &   $0$    &  $0$  \\           
     almost supersingular &          $\{\frac{1}{2},\frac{1}{2},\frac{1}{2},\frac{1}{2}\}$ &  $2$     &   $1$    &  $1$  \\          
     $1/3$ type           & $\{\frac{2}{3},\frac{2}{3},\frac{2}{3}\}$                      &  $1$     &   $2$    &  $1$  \\                        
     supersingular        & $\emptyset$                                                    &  $0$     &   $3$    &  $1$  
    \end{tabular}
\end{center}
Thus the unipotent dimension of $\Br(A)[p^{\infty}]$ is $0,0,1,2,3$, respectively. Recall that we have Ekedahl's symmetry
\[T^{02}=T^{13}, \quad T^{03}=T^{12}.\]

\begin{corollary} Let \(A\) be an abelian \(3\)-fold over \(k\). The following statements hold:

{\rm  (1)} If $A$ is ordinary, 
    then $\Br(A)[p^{\infty}]\cong (\Q_p/\Z_p)^{9-\rho}$. 
    
{\rm  (2)} If $A$ is almost ordinary, 
    then $\Br(A)[p^{\infty}]\cong (\Q_p/\Z_p)^{5-\rho}$. 

{\rm  (3)} If $A$ is almost supersingular, 
    then $\Br(A)[p^{\infty}]\cong (\Q_p/\Z_p)^{7-\rho}\oplus k$.

{\rm  (4)} If $A$ is of $1/3$-type,
    then $\Br(A)[p^{\infty}]\cong (\Q_p/\Z_p)^{9-\rho}\oplus U(k)$,
    where $U$ is an extension of $\G_a$ by $\G_a$.

{\rm  (5)} If $A$ is supersingular, 
    then $\Br(A)[p^{\infty}]\cong (\Q_p/\Z_p)^{15-\rho}\oplus U(k)$,
    where $U/k$ is a connected commutative unipotent algebraic group of dimension $3$.
\end{corollary}
\begin{proof} This follows from Theorem \ref{Mainthm} (1), (2), (3) and Proposition \ref{AbeT}. \end{proof}

\begin{remark} {\rm
    Serre classified all connected commutative unipotent algebraic group of dimension $2$, see \cite[Lemme VII.2.9.3]{Se2}. 
    Apart from the Witt group of length $2$: $W_2$ and $\G_a\times \G_a$, 
    for each element \[\sum_{i\geq 0} c_i T^i\in k_{\sigma}[T],\] there is an extension of $k$ by $k$, given by the formula 
    \begin{equation}
        (a_1,b_1)+(a_2,b_2)=(a_1+a_2, b_1+b_2+\sum_{i} c_i F(a_1,a_2)^{p^i})
    \end{equation}
    where $F(x_1,x_2)=\frac{1}{p}(x_1^p+x_2^p-(x_1+x_2)^p)$. 
Moreover, Serre proved that $\mathrm{Ext}^1(\G_a,\G_a)\cong k_{\sigma}[T]$. 
We will see, in Table \ref{Table2}, for any abelian $3$-fold of $1/3$-type, $U$ is always $\G_a\times \G_a$. 
In dimension $3$, all connected commutative unipotent algebraic group is isogenous to either $\G_a\times \G_a\times \G_a$, or $\W_2\times \G_a$, or $\W_3$. 
We will see, in Table \ref{Table2}, for any supersingular abelian $3$-fold, if its $a$-number is $2$ or $3$, then $U\cong \G_a\times \G_a\times \G_a$, 
if its $a$-number is $1$, then $U$ is isogenous to $\W_2\times \G_a$. 
}
\end{remark}

Recall that an abelian variety $A/k$ of dimension $g$ is called {\em superspecial}
if $A\cong\prod_{i=1}^g  E_i$, where $E_i$ for $i=1,...,g$ are
supersingular elliptic curves. When \(g\geqslant 2\), it is known 
that all superspecial abelian varieties of given dimension are isomorphic \cite[Theorem 2]{Oo2a}.

\begin{proposition}\label{Superspecial}
Let $A$ be the superspecial abelian variety of dimension $g$.
    Then we have $\HH^2_{\fppf}(X,\mu_p)\cong (\Z/p)^{g^2}\oplus k^{\frac{g(g-1)}{2}}$
and $\Br(A)\cong k^{ \frac{g(g-1)}{2}}$. 
\end{proposition}

\begin{lemma}\label{Homlemma}
    Let $E$ be a supersingular elliptic curve, and let $E[p]$ be its $p$-torsion group scheme. 
    Then 
    \begin{equation*}
        \Hom_k(E[p],E[p])\cong \F_{p^2} \oplus k,
    \end{equation*}
    with algebra structure given by $(x,a)\cdot(y,b)=(xy,xb+ay^p)$. 
\end{lemma}
\begin{proof} The contravariant Dieudonné module functor is an anti-equivalence 
between the category of finite commutative group schemes of $p$-power order over $k$ 
and the category of finite length modules over the Dieudonné ring (see, for example, \cite[Theorem III.1]{Fon}). In particular, we have 
\begin{equation*}
    \Hom_k(E[p],E[p])=\Hom_{W_{\sigma}[F,V]}(\mathbb{D}(E[p]),\mathbb{D}(E[p])). 
\end{equation*}
It is known (\cite[Chapter III]{De}) that $\mathbb{D}(E[p])$ is a 
$2$-dimensional vectors space over $k$ with basis $\{e_1,e_2\}$ such that 
\begin{equation*}
    Fe_1=e_2, \quad Fe_2=0, \quad Ve_1=e_2, \quad Ve_2=0.
\end{equation*}
Suppose $f:\mathbb{D}(E[p])\rightarrow \mathbb{D}(E[p])$ is
such that $f(e_1)=a_1e_1+a_2e_2$, $f(e_2)=b_1e_1+b_2e_2$. Then
$f$ satisfies $f(F(e_i))=F(f(e_i))$ and $f(V(e_i))=V(f(e_i))$. 
This is equivalent to $a_1=b_2^{\sigma}=b_2^{\sigma^{-1}}$, and $b_1=0$. 
So $a_1$ can be any element of $\F_{p^2}$, and $b_2=a_1^{\sigma}=a_1^{\sigma^{-1}}$, $b_1=0$, and $a_2$ can be any element in $k$. 
This shows that 
\[ \Hom_{W_{\sigma}[F,V]}(\mathbb{D}(E[p]),\mathbb{D}(E[p]))\cong \F_{p^2}\oplus k, \]
proving the lemma. \end{proof}

\noindent{\bf Proof of Proposition \ref{Superspecial}:} We use induction on $g$, and we use the product formula in \cite{Sk}. 
For $g=1$, we have $\HH^2_{\fppf}(X,\mu_p)\cong \Pic(X)/p\cong \Z/p$. 

By \cite[Theorem 1.3]{Sk}, for pointed projective varieties $X$ and $Y$, we have
\[ \HH^2_{\fppf}(X\times Y,\mu_{p})\cong \HH^2_{\fppf}(X,\mu_{p})\oplus \HH^2_{\fppf}(Y,\mu_{p}) \oplus \Hom_k(S_X,S_Y^{\vee}),\]
where $S_X$ is the finite commutative group $k$-scheme whose Cartier dual $S_X^{\vee}$ is $\Pic_{X/k}[p]:=ker(\Pic_{X/k}\xrightarrow{p} \Pic_{X/k})$, and similarly for $S_Y$.
Taking $X$ and $Y$ to be $E$ and $E^{g-1}$ respectively, 
we have \begin{equation}
    \HH^2_{\fppf}(E\times E^{g-1},\mu_p)\cong \HH^2_{\fppf}(E,\mu_{p})\oplus \HH^2_{\fppf}(E^{g-1},\mu_{p}) \oplus \Hom_k(E[p]^{\vee},E^{g-1}[p]).
\end{equation} 
The group scheme $E[p]$ is self dual by Weil pairing. By Lemma \ref{Homlemma}, we have 
\begin{equation}
    \Hom_k(E[p]^{\vee},E^{g-1}[p])\cong (\Z/p)^{ 2g-2}\oplus k^{g-1}.
\end{equation}
Induction on the dimension $g$ now gives the desired formula for $\HH^2_{\fppf}(E^{g},\mu_p)$.

Using the Kummer sequence and the Tate conjecture for abelian varieties, 
we see that $\Br(X)[p^\infty]$ only contains the unipotent part. 
By counting the dimension of the unipotent part, we conclude that $\Br(X)[p^{\infty}]\cong k^{\frac{g(g-1)}{2}}$. $\square$

\medskip

In the rest of this section, we compute the domino numbers $T^{ij}$ for all abelian $4$-folds, 
for documentation. There are $9$ different kinds of Newton polygons: 
the supersingular case (red), the supersingular $1/3$ type (green), the $1/4$ type (violet), the almost supersingular case (yellow), 
the ordinary $1/3$ case (brown), the half-ordinary-half-supersingular (orange), the almost ordinary (blue) and the ordinary case (black). 

\begin{center}
    \begin{tikzpicture}\label{NewtonPolygon}
        \draw [step=1cm,gray,very thin] (0,0) grid (8,4);
        \draw (0,0) -- (4,0);
        \draw (4,0) -- (8,4);
        \draw [blue] (3,0) -- (5,1);
        \draw [orange] (2,0) -- (6,2);
        \draw [brown]  (1,0) -- (4,1);
        \draw [brown]  (4,1) -- (7,3);
        \draw [yellow] (1,0) -- (7,3);
        \draw [purple] (0,0) -- (4,1);
        \draw [purple] (4,1) --(8,4);
        \draw [green] (0,0) -- (3,1);
        \draw [green] (3,1) -- (5,2);
        \draw [green] (5,2) -- (8,4);
        \draw [red](0,0) -- (8,4);
    \end{tikzpicture}
\end{center}

For the supersingular abelian 4-folds: 

\begin{center}
    \begin{tabular}{ c|c|c|c } 
     $T^{ij}$             &   $i=0$  & $i=1$  & $i=2$   \\ 
     $j=4$                &    $1$   &   $4$  &  $6$    \\ 
     $j=3$                &    $4$   &   $18$ &  $4$    \\           
     $j=2$                &    $6$   &   $4$  &  $1$    \\          
    \end{tabular}
\end{center}

For the supersingular $1/3$ slope (green): 

\begin{center}
    \begin{tabular}{ c|c|c|c } 
     $T^{ij}$             &   $i=0$  & $i=1$  & $i=2$   \\ 
     $j=4$                &    $1$   &   $4$  &  $4$    \\ 
     $j=3$                &    $4$   &   $12$  &  $4$    \\           
     $j=2$                &    $4$   &   $4$  &  $1$    \\          
    \end{tabular}
\end{center}

For the almost supersingular slope (yellow): 

\begin{center}
    \begin{tabular}{ c|c|c|c } 
     $T^{ij}$             &   $i=0$  & $i=1$  & $i=2$   \\ 
     $j=4$                &    $1$   &   $4$  &  $3$    \\ 
     $j=3$                &    $4$   &   $8$  &  $4$    \\           
     $j=2$                &    $3$   &   $4$  &  $1$    \\          
    \end{tabular}
\end{center}

For the $1/4$ (purple): 

\begin{center}
    \begin{tabular}{ c|c|c|c } 
     $T^{ij}$             &   $i=0$  & $i=1$  & $i=2$   \\ 
     $j=4$                &    $1$   &   $3$  &  $3$    \\ 
     $j=3$                &    $3$   &   $9$  &  $3$    \\           
     $j=2$                &    $3$   &   $3$  &  $1$    \\          
    \end{tabular}
\end{center}

For the ordinary $1/3$ slope (brown): 

\begin{center}
    \begin{tabular}{ c|c|c|c } 
     $T^{ij}$             &   $i=0$  & $i=1$  & $i=2$   \\ 
     $j=4$                &    $1$   &   $3$  &  $2$    \\ 
     $j=3$                &    $3$   &   $6$  &  $3$    \\           
     $j=2$                &    $2$   &   $3$  &  $1$    \\          
    \end{tabular}
\end{center}

For the half-ordinary-half-supersingular (orange): 

\begin{center}
    \begin{tabular}{ c|c|c|c } 
     $T^{ij}$             &   $i=0$  & $i=1$  & $i=2$   \\ 
     $j=4$                &    $1$   &   $2$  &  $1$    \\ 
     $j=3$                &    $2$   &   $4$  &  $2$    \\           
     $j=2$                &    $1$   &   $2$  &  $1.$    \\          
    \end{tabular}
\end{center}

\section{Product with an Ordinary Variety.}

In general, it is difficult to formulate a K\"unneth type formula for the $p$-primary torsion subgroup $\Br(X\times Y)[p^\infty]$, 
due to the complexity of the Künneth formula for Hodge-Witt cohomology groups $\HH^q(X,W\Omega^p_{X/k})$, as developed in \cite{Ek2}. 
The resulting expression is often highly nontrivial and involves intricate interactions between the cohomology groups of $X$ and $Y$. 

To illustrate this phenomenon, we consider a concrete example given by Ekedahl: the self-product of an Igusa surface. 
This example highlights the subtle behavior of the unipotent part of the Brauer group under products. 
Using Ekedahl's Künneth formula, we show that if $X$ is ordinary, then the unipotent part of $\Br(X \times Y)$ is naturally isomorphic to that of $\Br(Y)$. 

\begin{proposition}\label{Igusa} Suppose $p=2$. Let $I=(E_1\times E_2)/(\Z/2)$ be an Igusa surface over $k$ in \ref{hyperelliptic}. 
    Then the following statements hold: 

    {\rm    (1)} when $X=I$, then $U=J=0$. 
    
    {\rm    (2)} when $X=I\times I$, then $U\cong \G_a$ and $J=0$. 
\end{proposition}

\begin{proof} This follows directly from Ekedahl’s computation \cite[III.8]{Ek2}, together with Proposition \ref{Mainthm}. \end{proof}

Let $X$ and $Y$ be two proper smooth variety over $k$. According to Ekedahl \cite{Ek2}, 
there is an quasi-isomorphism of $R$-complexes in $D_c^b(R)$: 
\[R\Gamma(X,W\Omega^\bullet_{X/k}) \text{\ding{84}}^{\Ld}_R R\Gamma(Y,W\Omega^\bullet_{Y/k})\cong R\Gamma(X\times Y,W\Omega^\bullet_{X\times Y/k}) \]
where \ding{84} is the universal product of $R$-modules defined in \cite{Ek2}. 
In general, this K\"unneth formula is highly nontrivial, for instance, 
the domino $U_1$ will be the $0$-th cohomology of $E_{1/2}\text{\ding{84}}^{\Ld}_{R}E_{1/2}$. 
However, there is a special case where the formula simplifies considerably: when one of the varieties is ordinary. 

\begin{theorem}[Ekedahl's K\"unneth formula]\label{X*Y}
    Let $X$ be an ordinary propor smooth variety over $k$, and let $Y$ be an arbitrary proper smooth variety over $k$, 
    then we have the following exact sequence for each pair of $i,j$: 
    \begin{equation*}
        \begin{aligned}
            0\rightarrow \bigoplus_{\substack{
            i_1+i_2=i\\
            j_1+j_2=j}}\HH^{j_1}(X,W\Omega^{i_1}_{X/k})\otimes \HH^{j_2}(Y,W\Omega^{i_2}_{X/k})\rightarrow \HH^j(X\times Y,W\Omega^i_{X\times Y/k})\\
            \rightarrow \bigoplus_{\substack{
            i_1+i_2=i\\
            j_1+j_2=j+1}}\Tor_1^W(\HH^{j_1}(X,W\Omega^{i_1}_{X/k}),\HH^{j_2}(Y,W\Omega^{i_2}_{X/k}))\rightarrow 0.
        \end{aligned}
    \end{equation*}
\end{theorem}

Here, for two $R^0$ modules $M$, $N$, assume $M$ is finitely generated and the Frobenius $F$ acts bijectively on $M$, 
then both $M\otimes N$ and $\Tor_{1}^W(M,N)$ can be equipped with natural $R^0$ module structures as follows. 
Since $F$ is bijective on $M$, we may choose a free resolution of $M$ of length $2$, $0\rightarrow W^{n_1} \rightarrow W^{n_2} \rightarrow M$, 
where the Frobenius action on each free module is also bijective. 
Tensoring this complex with $N$, we obtain the complex $0\rightarrow N^{\oplus n_1}\rightarrow N^{\oplus n_2}$,
on which the actions of $F$ and $V$ are induced from those on $N$. 
Explicitly, the operators acts as $F(x\otimes y)=F(x)\otimes F(y),\quad V(x\otimes y)=F^{-1}(x)\otimes V(y).$
Taking cohomology, we obtain $\Tor^W_i(M,N)$ with an induced $R^0$-module structure. 

\begin{theorem}\label{prodord}
    Let $X$ be an ordinary propor smooth variety over $k$, 
    and let $Y$ be an arbitrary proper smooth variety over $k$. 
    The following statements hold: 
    
{\rm (1)} The unipotent parts of $\Br(X\times Y)[p^\infty]$ and $\Br(Y)[p^\infty]$ are isomorphic. 

{\rm (2)} Assume both $\Pic_{X/k}$ and $\Pic_{Y/k}$ are smooth, and the following two groups are free: 
\[ \HH^1(X,W\Omega^1_{X/k}),\quad \HH^1(Y,W\Omega^1_{X/k}),\]
then the non-divisible part of $\Br(X\times Y)[p^\infty]$, $\HH^3(X\times Y,\Z_p(1))[p^\infty]$, 
is isomorphic to the direct sum of $\HH^3(X,\Z_p(1))[p^\infty]$ and $\HH^3(Y,\Z_p(1))[p^\infty]$.  
\end{theorem}

\begin{proof} For short, let us use $\HH^{p,q}(X)$ to denote $\HH^q(X,W\Omega^p_{X/k})$ in this proof
(and similarly for $Y$). 
All these cohomology groups are successive extensions of 
$R^0$-modules of type I or $k[[V]]$, or $k[[dV]]$. 
The map $(F-1)$ is surjective on all of these modules, 
so the $F$-invariant subspace of the extension is the extension of the $F$-invariant subspaces. 
We will use this observation implicitly throughout this proof. 

For (1), we have $\HH^{0,0}(X)\cong \HH^{0,0}(Y)\cong W$, 
so the $F=1$ parts of $\HH^{1,2}(X) \otimes \HH^{0,0}(Y)$ and $\HH^{0,0}(X)\otimes \HH^{1,2}(Y)$ 
are $\HH^3(X,\Z_p(1))[p^\infty]$ and $\HH^3(Y,\Z_p(1))[p^\infty]$, respectively. 
It suffices to show that the $F=1$ parts of the other twelve terms in the K\"unneth formula do not contain unipotent parts. 

By Remark \ref{EasyVan}, $\HH^{1,0}(Y)$, $\HH^{0,1}(Y)$ and $\HH^{1,1}(Y)$ are all finitely generated, 
and the cohomologies of $X$ are always finitely generated,
so $\HH^{0,2}(X)\otimes \HH^{1,0}(Y)$, $\HH^{1,1}(X)\otimes \HH^{0,1}(Y)$, $\HH^{0,1}(X)\otimes \HH^{1,1}(Y)$
are all finitely generated $W$-modules as well, so their $F$-invariant subgroups have trivial unipotent parts. 
For $\HH^{1,0}(X)\otimes \HH^{0,2}(Y)$, 
$\HH^{1,0}(X)$ can be written as a direct sum of several copies of $W$ and torsion $W$-modules of the form $W/p^n$ (on which $F$ maps $1$ to $1$).   
By a easy calculation we have $(W/p^n)\otimes \HH^{0,2}(Y)\cong \HH^{0,2}(Y)/p^n$, 
its $F$-fixed part cannot have unipotent part, because it has a filtration by Theorem \ref{coherent} whose graded pieces are either modules of type $I$ or $k[[V]]$. 
The  $F$-fixed part of $W\otimes \HH^{0,2}(Y)$ cannot have unipotent part either.
For the eight $\Tor_1^W$ groups, notice that $\HH^{0,0}(Y)$, $\HH^{0,1}(Y)$, $\HH^{0,1}(X)$ and $\HH^{0,0}(X)$ are free, 
their $\Tor_1^W$ groups with any $W$-modules are always trivial. 
The groups $\HH^{0,3}(X)$, $\HH^{1,0}(Y)$, $\HH^{0,2}(X)$ and $\HH^{1,1}(Y)$ are all finitely generated, 
so $\Tor_1^W(\HH^{0,3}(X),\HH^{1,0}(Y))$ and $\Tor_1^W(\HH^{0,2}(X),\HH^{1,1}(Y)$ are both finitely generated as well, 
so their $F$-fixed parts can't have unipotent group. 
For $\Tor_1^W(\HH^{1,1}(X),\HH^{0,2}(Y))$ and the group $\Tor_1^W(\HH^{1,0}(X),\HH^{0,3}(Y))$, 
we can write $\HH^{1,1}(X)$ and $\HH^{1,0}(X)$ as a direct sum of $W$'s and modules of the form $W/p^n$, 
then by a easy calculation $\Tor_1^W(W/p^n, \HH^{0,3}(Y))\cong \HH^{0,3}(Y)[p^n]$ and similarly 
$\Tor_1^W(W/p^n, \HH^{0,2}(Y))\cong \HH^{0,2}(Y)[p^n]$, so their $F=1$ part can't have unipotent parts as well. This proves (1). 

For (2), the smoothness of the $\Pic_{X/k}$ and $\Pic_{Y/k}$ implies that 
$\HH^{0,2}(X)$ and $\HH^{0,2}(Y)$ doesn't contain $V$-torsion(\ref{Eke81}). 
The two groups $\HH^{1,0}(X)$ and $\HH^{1,0}(Y)$ are free because 
$W\Omega^i_{X/k}$ doesn't have $p$-torsion. 
By our assumption, and the fact that $\HH^{0,1}(X)$ and $\HH^{0,1}(Y)$ are both free 
(they are $\TC(\Pic_{X/k})$ and $\TC(\Pic_{Y/k})$ respectively), 
the four tensor products are all free $W$ module. 
The eight $\Tor_1^W$ groups are all zero by our freeness assumption. \end{proof}

\begin{remark}
    From Ekedahl's duality theorem \cite{Ek1}, there is no nilpotent torsion part in $\HH^1(X,W\Omega^i_{X/k})$ for any $i$. 
    As a consequence the torsion part of $\HH^1(X,W\Omega^1_{X/k})$ is always isomorphic to 
    $\NS(X)[p^\infty]\otimes_{\Z_p} W$. 
    I don't know how to show this without invoking Ekedahl's duality. 
\end{remark}

\part{Brauer Groups of Abelian Varieties}\label{Part2}

The theory of the slope spectral sequence provides a description of the $p$-primary part of the Brauer group 
$\Br(X)[p^\infty]$ by expressing it as a direct sum of $(\Q_p/\Z_p)^{r-\rho}$ and $\HH^3(X, \Z_p(1))$. 
The latter is an extension of a finite group $J$ by the group of $k$-points 
of a finite-dimensional connected commutative unipotent group $U$. 
However, the structure of $U$ remains somewhat mysterious. 
While its dimension is known from the work of Crew and Ekedahl, commutative unipotent groups in characteristic $p$ 
exhibit subtle and intricate behavior. For instance, the Witt group $\W_2$ (of dimension 2) differs significantly from the product $\G_a \times \G_a$: $\W_2$ is 
not annihilated by $p$, whereas $\G_a \times \G_a$ is. In fact, Serre showed that any connected 
commutative unipotent group $U$ is isogenous to a product of Witt vector groups $\W_{n_i}$, 
meaning there exists a surjective morphism $U \to \prod_i \W_{n_i}$ with finite kernel. 
The multiset of integers ${n_i}$--called the isogeny class of $U$--is uniquely determined by, 
and determines, the dimensions of the finite subgroups $U[p^n]$ killed by $p^n$ for all $n$.

Our goal is to determine the isogeny class of $U$, at least for certain easy varieties. 
To this end, we employ a formula introduced by Skorobogatov, 
which allows explicit computation of $U$ in the case of abelian varieties when $p \neq 2$. 
This formula expresses the $p$-torsion subgroup $\Br(A)[p^n]$ (for any $n$) in terms of a certain group of homomorphisms 
from $A[p^n]$ to $A^\vee[p^n]$. By solving equations for the relations of the Dieudonné modules, 
we compute the dimension of $U[p]$ for any principally polarized abelian variety $(A, \iota)$ 
from its Ekedahl-Oort type (see Theorem \ref{E-O}). 
Together with the dimension formula (\ref{Abe}), this enables us to determine the isogeny class of $U$ 
for any principally polarized abelian threefold (see Table \ref{abe3table}). 
Combining this result with some more calculations for non-polarized group schemes using Kraft's classification theorem, 
we are able to classify the isogeny class of $U$ for all abelian threefolds, 
in terms of their Newton polygons and $a$-numbers (see Table \ref{Table2}).

The explicit structure of $U$ is highly intricate. 
For completeness, we include a detailed computation of $\Hom(A[p^2], A^{\vee}[p^2])^{\mathrm{sym}}$ 
in the case of a principally polarized supergeneral abelian variety, 
illustrating the complexity of the problem. 
We also investigate the distribution of this Hom group across the moduli 
space of principally polarized abelian varieties.

\chapter{A Formula by Skorobogatov} \label{Ch3}

Let $k$, as before, be an algebraically closed field of characteristic $p>0$, 
and let $A$ be an abelian variety of dimension $g$ over $k$. 
In Theorem \ref{Abe}, using Illusie's theory of slope spectral sequence, 
we showed that $\Br(A)[p^\infty]$ is a direct sum of the group of $k$-points 
of a connected commutative unipotent group $U$ of dimension $T^{0,2}$ 
and a divisible part $(\Q_p/\Z_p)^{r-\rho}$, where $r$ is the multiplicity of slope $1$ in $\HH_{\cris}^2(A/W)$, 
and $\rho$ is the Picard number. 
In \cite{Skor}, Skorobogatov provided the following explicit formula to 
determine the structure of $U$, when $p\neq 2$. 

\begin{theorem}\label{Skorobogatov}{\bf (Skorobogatov)}
    Assume $p\neq 2$, and $k$ is algebraically closed of characteristic $p$. For any abelian variety $A/k$, 
    we have the following isomorphism: 
    \[\Br(A)[p^n]\cong \Hom_k(A[p^n],A^{\vee}[p^n])^{\sym}/(\Hom(A,A^{\vee})^{\sym}/p^n),\]
    where \emph{$\sym$} denotes morphisms $f$ such that $f^\vee\cong f$, 
    under the identifications given by the Weil pairings: 
    \[w_A:A[p^n]^{\vee}\cong A^{\vee}[p^n],\quad w_{A^\vee}: A^{\vee}[p^n]^{\vee}\cong A^{\vee\vee}[p^n]\cong A[p^n]\]
    those morphisms $f$ such that $f^{\vee}\cong f$. 
\end{theorem}

For completeness, we demonstrate the proof as following. The proof is completely due to Skorobogatov \cite{Skor}. 

Write $\HH^i(A):=\HH^i_{\fppf}(A,\mu_{p^n})$. Note that $\HH^i(k)=0$ for $i\geqslant 1$.
We have the following morphisms of abelian varieties $A\times A\to A$: the two projections
$\pi_1$, $\pi_2$, and the sum map $m$.
We also have morphisms of abelian varieties $A\to A\times A$: the two embeddings
$i_1(x)=(x,0)$ and $i_2(x)=(0,x)$, and the antidiagonal $a(x)=(x,-x)$.
Let $s\colon A\times A\to A\times A$ be the involution $s(x,y)=(y,x)$.

For a contravariant functor $F$ define
$F(A\times A)_\prim:=\Ker(i_1^*)\cap\Ker(i_2^*)$.
Then we have a natural direct sum decomposition
$$F(A\times A)=F(A)\oplus F(A)\oplus F(A\times A)_\prim.$$
Let $F(A\times A)^\sym$ be the $s^*$-invariant subgroup of $F(A\times A)$.

A contravariant functor $F$ from the category of abelian varieties over $k$
to the category of $p$-primary torsion abelian groups has degree $d$
if $F(\spec(k))=0$ and $[m]^*=m^d$.
An additive functor has degree 1. It is well-known that a functor is additive if and only if
it preserves finite direct sums, see, e.g., [Stacks, Lemma 0DLP].
The spectral sequence attached to the projection $A\times B\to A$ shows that
sending $A$ to $\HH^1(A)$ commutes with finite direct sums, so this is an additive functor.
Sending $A$ to $\HH^2(A)$ is a quadratic functor. This is well-known when $p\neq{\rm char}(k)$,
because in this case $\mu_{p^n}$ is a smooth group scheme, so
fppf cohomology coincides with \'etale cohomology which is the exterior algebra of $\HH^1(A)$
(up to a twist). 

\begin{lemma} \label{2}   
Assume that $p\neq 2$. For any quadratic functor with values in $p$-primary torsion abelian
groups the map $m^*-\pi_1^*-\pi_2^*\colon F(A)\to F(A\times A)$
is injective with image contained in $F(A\times A)_\prim^\sym$.
\end{lemma}

\begin{proof}
Composing with $i_1^*$ and $i_2^*$ 
we check that $m^*(x)-\pi_1^*(x)-\pi_2^*(x)\in F(A\times A)^\sym_\prim$.
Composing with $a^*$, we see that $[-1]:A\to A$ acts on $\Ker(m^*-\pi_1^*-\pi_2^*)$ as
$-1$. Since $F$ is quadratic, $[-1]^*=\mathrm{Id}_{F(A)}$, thus $2$ annihilates 
$\Ker(m^*-\pi_1^*-\pi_2^*)$, implying that $\Ker(m^*-\pi_1^*-\pi_2^*)=0$. 
\end{proof}

\begin{proposition}\label{3}
The map $m^*-\pi_1^*-\pi_2^*\colon \HH^2(A)\to \HH^2(A\times A)$ is surjective onto 
$(1+s^*)\HH^2(A\times A)_\prim$. (Here $k$ can be an arbitrary field.)
\end{proposition}

\begin{proof}
    Let $\T$ be the torsor $[p^n]\colon A\to A$ with structure group $A[p^n]$.
We denote by $[\T]$ the class of $\T$ in $\HH^1_\fppf(A,A[p^n])$.
The cup-product in fppf cohomology of $A$ produces a class
$$[\T]\cup[\T]\in\HH^2_\fppf(A,A[p^n]\otimes A[p^n]).$$
We have classes $\pi_i^*[\T]\in\HH^1_\fppf(A\times_kA,A[p^n])$ for $i=1,2$.
The cup-product in fppf cohomology of $A\times_kA$ produces a class
$$\pi_1^*[\T] \cup \pi_2^*[\T]\in \HH^2_\fppf(A\times_kA,A[p^n]\otimes A[p^n]).$$

On the one hand, by \cite[Proposition  1.3]{Sk}, the map sending
$\psi\in\Hom_k(A[p^n]\otimes A[p^n],\mu_{p^n})$ to the push-forward
$\psi_*(\pi_1^*[\T] \cup \pi_2^*[\T])\in\HH^2(A\times_kA)$
defines an isomorphism
$$\Hom_k(A[p^n]\otimes A[p^n],\mu_{p^n})\tilde\longrightarrow\HH^2(A\times_kA)_\prim.$$
This implies that every element of $(1+s^*)\HH^2(A\times A)_\prim$ can be written as 
$\psi_*(\pi_1^*[\T] \cup \pi_2^*[\T])+\psi_*(\pi_2^*[\T] \cup \pi_1^*[\T])$ for some
$\psi$.

On the other hand, functoriality of cup-product gives
$m^*(\psi_*([\T] \cup [\T]))=\psi_*(m^*[\T]\cup m^*[\T])$, and similarly with $m^*$ replaced by 
$\pi_1^*$ or $\pi_2^*$. By the additivity of $\HH^1$, we have
$m^*[\T]=\pi_1^*[\T]+\pi_2^*[\T]$. It follows that
$$(m^*-\pi_1^*-\pi_2^*)(\psi_*([\T]\cup[\T]))=\psi_*(\pi_1^*[\T]\cup\pi_2^*[\T])+
\psi_*(\pi_2^*[\T]\cup\pi_1^*[\T]).$$
This proves the proposition.
\end{proof}
\begin{theorem}\label{1}
Assume that $p\neq 2$. Then the map $m^*-\pi_1^*-\pi_2^*$ induces an
isomorphism
\[\begin{aligned}
    &\HH^2(A)\tilde\longrightarrow\HH^2(A\times_kA)_\prim^\sym\\
    \cong&\Hom_k(A[p^n]\otimes A[p^n],\mu_{p^n})^\sym\cong\Hom_k(A[p^n],A^\vee[p^n])^\sym.
\end{aligned}\]
\end{theorem}
\begin{proof}
    Combine Lemma \ref{2} and Proposition \ref{3}.
\end{proof}

\medskip

The Kummer sequence
$$1\to\mu_{p^n}\to\G_m\to\G_m\to 1$$
is exact in the fppf topology. It gives a natural isomorphism
$\HH^1(A)\cong A^\vee[p^n](k)$
and an exact sequence
\begin{equation}
0\to \Pic(A)/p^n\to\HH^2(A)\to\Br(A)[p^n]\to 0.\label{e1}
\end{equation}
Now we show that: assume that $p\neq 2$. Then the map $m^*-\pi_1^*-\pi_2^*$ induces an
isomorphism 
\[\begin{aligned}
    &\Br(A)[p^n]\tilde\longrightarrow\Br(A\times A)[p^n]_\prim^\sym\\
    \cong& \Hom_k(A[p^n],A^\vee[p^n])^\sym/(\Hom(A,A^\vee)^\sym/p^n).
\end{aligned}\]

\noindent{\bf Proof of Theorem \ref{Skorobogatov}:} 
    The primitive part of (\ref{e1}) is the exact sequence
$$0\to\Hom(A,A^\vee)/p^n\to\HH^2(A\times_kA)_\prim\to\Br(A\times_kA)_\prim[p^n]\to 0.$$
Since $p$ is odd, we have $\HH^1(\Z/2,\Hom(A,A^\vee)/p^n)=0$, so passing to
the $\mathrm{s}^*$-invariant subgroups preserves exactness. It is well-known
that $m^*-\pi_1^*-\pi_2^*$ induces an isomorphism $\NS(A)\tilde\longrightarrow\Hom(A,A^\vee)^\sym$.
Using this and Theorem \ref{1} we see that the first two vertical maps
in the commutative diagram
\begin{small}
\begin{center}
    \begin{tikzcd}
        0\arrow{r}&\Hom(A,A^\vee)^\sym/p^n\arrow{r}&\HH^2(A\times_kA)_\prim^\sym\arrow{r}&\Br(A\times_kA)_\prim[p^n]\arrow{r}& 0\\
        0\arrow{r}& \Pic(A)/p^n\arrow{r}\arrow[u,"\cong"]&\HH^2(A)\arrow{r}\arrow[u,"\cong"]&\Br(A)[p^n]\arrow{r}\arrow{u}& 0
    \end{tikzcd}
\end{center}
\end{small}
are isomorphisms. Thus so is the third vertical map. \hfill $\Box$

\begin{remark}
    Similar to Theorem \ref{IllusieRaynaud}, 
    the isomorphism of Theorem \ref{Skorobogatov} 
    can be upgraded to an isomorphism between quasi-algebraic groups. 
    Indeed, following the proof of Theorem \ref{Skorobogatov}, one can construct a map between fppf sheaves: 
    \[R^2f_{\fppf,*}\mu_{p^n}/(\Pic(A)/p^n)\rightarrow \mathbf{Hom}(A[p^n],A^{\vee}[p^n])^{\sym}/(\Hom(A,A^{\vee})^{\sym}/p^n).\]
    On \(k\)-points, this map is an isomorphism, 
    and therefore it induces a isomorphism between their corresponding quasi-algebraic groups. 
    A detailed proof will be given in a forthcoming work. 
\end{remark}

\chapter{Dimension of $U[p]$}\label{Chap4}

From now on, we assume $p\neq 2$, and as before, let $k$ be an algebraically closed of characteristc $p$. 
In this chapter, we use Skorobogatov's formula to compute the dimension of $U[p]$, the $p$-torsion subgroup of $U$,  
for an arbitrary principally polarized abelian variety $(A,\iota)$ from the Ekedahl-Oort type of $A$. 
Note that $U[p]$ is always a finite product of $\G_a$. 
We also compute the dimension of $U[p]$ for an arbitrary (not necessarily polarized) abelian threefold, 
using Kraft's classification theorem. 
This allows us to determine the isogeny class of $U$ for all abelian threefolds. 

In the rest of this part, we carry out the following strategy. 
For two finite flat commutative group schemes \(G_1,G_2\) over \(k\), 
assume \(M_1, M_2\) are their covariant Dieudonn\'e modules, respectively. 
Let \(\mathbf{Hom}(G_1,G_2)\) be the Hom group scheme from \(G_1\) to \(G_2\), 
then this functor sends every \(k\)-scheme \(S\) to the abelian group \(\Hom_{S}(G_1\times S,G_2\times S)\). 
Define \(\mathbf{Hom}_{\mathbb{D}}(M_1,M_2)\)  to be the subgroup scheme of the 
\(W\)-module scheme \(\mathbf{Hom}_W(M_1,M_2)\) cut out by the Dieudonn\'e compatibility relations. 
Explicitly, for every \(k\)-algebra \(R\), 
the group scheme \(\mathbf{Hom}_{\mathbb{D}}(M_1,M_2)\) send \(R\) to the group of homomorphisms 
\[f\in \Hom_{W(R)}(M_1\otimes W(R),M_2\otimes W(R))\] 
such that the following two squares commute: 
\[\begin{tikzcd}
	{M_1\otimes W(R)} & {M_2\otimes W(R)} & {\sigma^*M_1\otimes W(R)} & {\sigma^*M_2\otimes W(R)} \\
	{\sigma^*M_1\otimes W(R)} & {\sigma^*M_2\otimes W(R)} & {M_1\otimes W(R)} & {M_2\otimes W(R).}
	\arrow["f", from=1-1, to=1-2]
	\arrow["{F_1\otimes id}"', from=1-1, to=2-1]
	\arrow["{F_2\otimes id}"', from=1-2, to=2-2]
	\arrow["f", from=1-3, to=1-4]
	\arrow["{V_1\otimes id}"', from=1-3, to=2-3]
	\arrow["{V_2\otimes id}"', from=1-4, to=2-4]
	\arrow["f", from=2-1, to=2-2]
	\arrow["f", from=2-3, to=2-4]
\end{tikzcd}\]

According to crystalline Dieudonn\'e theory, for each \(k\)-scheme \(S\), 
there is a contravariant functor \(\mathbb{D}\) from the category of 
finite flat commutative groups schemes over \(S\)
to the category of Dieudonné crystals over \(S\). 
We thus have a natural transformation
\[\mathbf{Hom}(G_1,G_2)\rightarrow \mathbf{Hom}_{\mathbb{D}}(M_1,M_2),\]
by mapping every \(x\in \mathbf{Hom}(G_1,G_2)(S)\) to the dual map of 
\(\mathbb{D}(x):\mathbb{D}(G_{2,S})\rightarrow \mathbb{D}(G_{1,S})\) 
evaluated at the canonical thickening \(W_m(S)\) with its natural power divided structure. 

On the level of \(k\)-points, by classical Dieudonné theory, it induces an isomorphism. 
Therefore upon passing to quasi-algebraic groups, the map 
\[\mathbf{Hom}(G_1,G_2)\rightarrow \mathbf{Hom}_{\mathbb{D}}(M_1,M_2)\]
induces an isomorphism. 
Consequently one may compute \(\mathbf{Hom}(G_1,G_2)\) by working with \(\mathbf{Hom}_{\mathbb{D}}(M_1,M_2)\). 
When \(p\neq 2\), the same argument applies after imposing the symmetric condition. 

\section{Principally Polarized Case: Ekedahl-Oort Type}

Recall that in \cite{Oo1}, Ekedahl and Oort stratified the moduli space of 
principally polarized abelian varities of dimension $g$ over $k$ 
according to the isomorphism class of the pair $(A[p],\iota|_{A[p]})$. 
There are precisely $2^g$ such isomorphism classes, 
each represented by an Ekedahl-Oort types $\varphi$, 
a sequence of integers $(\varphi(1),...,\varphi(g))$ satisfying 
\[\varphi(i)-\varphi(i-1)=0 \;\text{or} \;1,\quad \text{with}\;\varphi(0)=0.\]
The set of Ekedahl-Oort types corresponds bijectively to the set of subsets $P\subseteq \{1,...,g\}$, 
where the associated type $\varphi_P$ is defined recursively by: \[\varphi_P(i)=\varphi_P(i-1)+1 \quad \text{if and only if}\; i\in P.\]
For each given Ekedahl-Oort type $\varphi$, 
the $a$-number of $A[p]$ is $g-\varphi(g)$, and the $p$-rank of $A$ is the largest integer $f$ such that $\varphi(f)=f$. 

\begin{example} We highlight a few key examples: 
    \begin{itemize}
        \item If $A$ is ordinary, its Ekedahl-Oort type is $\{1,2,...,g\}$. 
        \item If $A$ is superspecial (supersingular with maximal $a$-number: $g$), the Ekedahl-Oort type is $\{0,0,...,0\}$. 
        \item If $A$ is supergeneral (supersingular with minimal $a$-number: $1$), the Ekedahl-Oort type is $\{0,1,2,...,g-1\}$. 
    \end{itemize}
\end{example}

We will prove: 
\begin{theorem} \label{E-O}
    Suppose $(A,\iota)$ is a principally polarized abelian variety of dimension $g$ 
    with Ekedahl-Oort type $\varphi=\varphi_P$, associated to a subset 
    $P=\{1\leqslant m_1<...<m_h\leqslant g\}$ of size $h$. 
    Then the dimension of $U[p]$ is given by: 
    \[\frac{(g-h)(g-h-1)}{2}+\sum_{i=1}^{h}(m_i-i).\]
\end{theorem}

This implies, for instance, let $(A,\iota)$ be a principally polarized abelian variety: 
\begin{itemize}
    \item If $(A,\iota)$ is ordinary, then $U[p]$ is trivial. This aligns with the fact that $U$ itself is trivial.
    \item If $(A,\iota)$ is superspecial, then $\dim U[p]=\frac{g(g-1)}{2}$, consistent with Theorem \ref{Superspecial}. 
    \item If $(A,\iota)$ is supergeneral, then $U[p]$ has dimension $g-1$. Here recall that $\dim U=\frac{g(g-1)}{2}$. 
\end{itemize}

Granting Theorem \ref{E-O}, it is easy to obtain: 

\begin{theorem}
    The possible isogeny classes of $U$ for principally polarized abelian threefolds are listed in the table below:
\begin{center}{\rm\label{abe3table}
    \begin{tabular}{ c|c|c|c|c } \label{Table2}
     E-O type             & Newton polygon                                                  & $dim(U[p])$& $dim(U)$ & isogeny class of $U$\\ 
     $\{0,0,0\}$          &   supersingular, $a=3$                                          &  $3$       & $3$     & $\G_a\times \G_a\times \G_a$  \\ 
     $\{0,0,1\}$          &   supersingular, $a=2$                                          &  $3$       & $3$     & $\G_a\times \G_a\times \G_a$\\           
     $\{0,1,1\}$          &   $1/3$ type                                                    &  $2$       & $2$     & $\G_a\times \G_a$  \\ 
     $\{0,1,2\}$          &   supersingular, $a=1$                                          &  $2$       & $3$     & $\G_a\times \W_2$\\     
                          &   $1/3$ type                                                    &  $2$       & $2$     & $\G_a\times \G_a$\\                 
     $\{1,1,1\}$          &   almost supersingular                                          &  $1$       & $1$     & $\G_a$  \\ 
     $\{1,1,2\}$          &   almost supersingular                                          &  $1$       & $1$     & $\G_a$\\           
     $\{1,2,2\}$          &   almost ordinary                                               &  $0$       & $0$     & $0$  \\ 
     $\{1,2,3\}$          &   ordinary                                                      &  $0$       & $0$     & $0$\\           
    \end{tabular}}
\end{center}
\end{theorem}

\noindent\emph{Proof:} 
We need to show the only match-ups between E-O type and Newton polygon are the ones in the table. 
For $\varphi=\{0,0,0\}$, the $a$-number equals $3$ so it must be superspecial; 
for $\varphi=\{1,2,3\}$ and $\{1,2,2\}$, their $p$-rank equals $3$ and $2$ respectively 
so they must be ordinary and almost ordinary respectively. 
For $\varphi=\{1,1,2\}$ and $\{1,1,1\}$, their $p$-rank equals $1$ so they both are almost supersingular. 
For  $\{0,0,1\}$ and  $\{0,1,1\}$,
 by a result of Oort \cite[Theorem (8.3)]{Oo1} the former one is supersingular and the latter one is $1/3$ type. 
 Now the result follows from Theorem \ref{E-O}.
 $\square$

From the table above, we observe that for principally polarized abelian threefolds, the group 
$U$ is not annihilated by $p$ if and only if $A$ is supergeneral. 
In all other cases, $U$ is either trivial or a product of $\G_a$. 
In the supergeneral case, $U$ is isogenous to $\G_a\times \W_2$. 

The next two sections are devoted to prove Theorem \ref{E-O}.

\medskip

\section{The Hom Group}\label{SectionTheHomGroup}

Let $A$ be an arbitrary (not necessarily polarized) abelian variety. 
Let $M$ denote the covariant Dieudonn\'e module of $A[p^\infty]$, equipped with a $W$-basis $\mathcal{B}$. 
Assume the Frobenius $F$ and the Verschibung $V$ have matrices $[\mathcal{F}_M]$ and $[\mathcal{V}_M]$, respectively, with respect to $\mathcal{B}$. 
Then the Dieudonn\'e module of $A^\vee[p^\infty]$, which is isomorphic to $M^\vee$ via the first Weil pairing,
can be described as a free $W$ module with dual basis $\mathcal{B}^*$. 
Under this basis, the matrices of $\mathcal{F}_{M^\vee}$ and $\mathcal{V}_{M^\vee}$ are given by 
$[\mathcal{V}_M]^{t\sigma}$ and $[\mathcal{F}_M]^{t\sigma^{-1}}$, respectively, 
where $t$ denotes the transpose, and $\sigma$ and $\sigma^{-1}$ acts on matrix entries via the Frobenius and the inverse. 
Here, we are using the explicit formula: for every $g\in M^{\vee}:M\rightarrow k$ and every $x\in M$, we have
\[\mathcal{F}_{M^{\vee}}(g)(x)=g(\mathcal{V}_M(x))^{\sigma}, \quad \mathcal{V}_{M^{\vee}}(g)(x)=g(\mathcal{F}_M(x))^{\sigma^{-1}}.\]

Now consider a morphism $f:A[p^n]\rightarrow A^\vee[p^n]$, it corresponds to a map $f:M/p^n\rightarrow M^{\vee}/p^n$. 
With respect to the basis $\mathcal{B}$ and $\mathcal{B}^*$, 
this is given by a matrix $[f]\in \mathrm{Mat}_{2g\times 2g}(\W_{n})$ satisfying 

\[[f][\mathcal{F}_M]=[V_M]^{t\sigma}[f]^\sigma,\]
\[[f][\mathcal{V}_M]=[F_M]^{t\sigma^{-1}}[f]^{\sigma^{-1}}.\]

Here, we are using column-vector notation for coefficients: 
\[\begin{aligned}
    f(e_1,...,e_{2g})=(e_1,...,e_{2g})[f],\\
    \mathcal{F}(e_1,...,e_{2g})=(e_1,...,e_{2g})[\mathcal{F}_M],\\
    \mathcal{V}(e_1,...,e_{2g})=(e_1,...,e_{2g})[\mathcal{V}_M].
\end{aligned}\]

Next, we determine the matrix of the dual morphism $f^{\vee}: M/p^n\rightarrow M^{\vee}/p^n$. 
Recall that there are canonical identifications: 
\[A[p^n]^{\vee}\cong A^\vee[p^n],\quad A^\vee[p^n]^\vee\cong A^{\vee\vee}[p^n]\cong A[p^n],\]
but these isomorphisms arise from two different Weil pairings, differing by a sign. 
Consequently, under the dual basis $\mathcal{B}$ and $\mathcal{B}^*$, 
the matrix representing $f^\vee$ is given by $-[f]^t$. 
Therefore, the group of symmetric homomorphisms $\Hom(A[p^n],A^{\vee}[p^n])^{\sym}$ corresponds canonically to
the group of matrices $[f]\in \mathrm{Mat}_{2g\times 2g}(\W_{n})$ satisfying:  
\[[f][\mathcal{F}_M]=[V_M]^{t\sigma}[f]^\sigma,\]
\[[f][\mathcal{V}_M]=[F_M]^{t\sigma^{-1}}[f]^{\sigma^{-1}},\]
\[[f]=-[f]^t.\]

\begin{lemma}
    Under the condition $[f]=-[f]^t$, the two Dieudonné compatibility equations: 
\[[f][\mathcal{F}_M]=[V_M]^{t\sigma}[f]^\sigma,\quad [f][\mathcal{V}_M]=[F_M]^{t\sigma^{-1}}[f]^{\sigma^{-1}}.\]
    are equivalent to each other. 
\end{lemma}

\noindent\emph{Proof:} Take the $\sigma$-twisted transpose of the second equation yields the first. $\square$

In the where $A$ is principally polarized, the covariant and contravariant Dieudonné modules are canonically isomorphic.
For simplicity, we work in the covariant framework here. 
Note that in this exposition, we do not use the principal polarization to identify the dual basis 
$\mathcal{B}^*$ with the original basis $\mathcal{B}$. 
However, in \cite{Skor}such an identification was made.
While the notation differs, the resulting computations are equivalent.

\section{Calculation: Combinatorics}

Let $\varphi=\varphi_P$ be the elementary sequence associated with the subset $P=\{m_1<...<m_h\}\subseteq\{1,2,...,g\}$ of size $h$. 
Let $(A,\iota)$ be a principally polarized abelian variety with Ekedahl-Oort type $\varphi$. 
Note that the $a$-number of $A$ is equal to $g-h$. 
Our goal is to describe the group $\Hom(A[p],A^{\vee}[p])^{\sym}$, 
the symmetric part of the homomoprhism group; 
that is, the group of homomorphisms $f:A[p]\rightarrow A^\vee[p]$ such that $f^{\vee}=f$. 
From our analysis in the previous section, the group is isomorphic to the space of $2g\times 2g$ 
matrices in $k$ satisfying: 
\[[f][\mathcal{F}_M]=[V_M]^{t\sigma}[f]^\sigma,\]
\[[f]=-[f]^t.\]

We now describe the matrices $[\mathcal{F}_M]$ and $[\mathcal{V}_M]$. 
Given $\varphi$, we extend it to a final sequence $(\psi(1),...,\psi(2g))$ by symmetrically setting 
\[\psi(i)=\varphi(i), \;for\; i=1,...,g; \qquad \psi(2g-i)=\varphi(i)+g-i, \; for\; i=1,...,g.\]

\begin{proposition}\emph{\textbf{(Oort)}}
    The covariant Dieudonné module $M$ of the group scheme $G_{\varphi}$ 
    is a $2g$-dimensional $k$-vector space with basis $\{Z_i|i=1,2,...,2g\}$,
    and the operators $\mathcal{F}$ and $\mathcal{V}$ acts as follows: 
    \[\begin{aligned}
        &&\mathcal{F}(X_i)=Z_i&&&\mathcal{F}(Y_i)=0,\\
        &&\mathcal{V}(Z_i)=0,&&&\mathcal{V}(Z_{2g+1-i})=\pm Y_i,
    \end{aligned}\]
    where the sets $\{m_i\}$, $\{n_i\}$, $\{X_i\}$, $Y_i$ are defined as follows: 
    \begin{itemize}
        \item The indices $0<m_1<...<m_g\leqslant 2g$ are those $i\in\{1,...,2g\}$ such that $\psi(i)=\psi(i-1)+1$. 
        For these indices, we set $X_i=Z_{m_i}$.
        \item The indices $0<n_g<...<n_1\leqslant 2g$ are those $j\in\{1,...,2g\}$ such that $\psi(j)=\psi(j-1)$, and we set $Y_j=Z_{n_j}$.
    \end{itemize} 
    In particular, for $i=1,...,g$, we set $V(Z_{2g+1-i})=+Y_i$ if and only if $Z_{2g+1-i}\in \{Y_1,...,Y_g\}$.  
\end{proposition}

\begin{proof}
    This is Section 9 of \cite{Oo1}.
\end{proof}

Using the basis $\{X_1,...,X_g,Y_1,...,Y_g\}$, we now describe the matrix $[\mathcal{F}_M]$ and $[\mathcal{V}_M]$.  
Let $Q:=\{1,...,g\}\backslash P=\{n_{h+1}>...>n_g\}$, 
define two $g\times g$ matrices $U,V$ by setting:
\[\begin{aligned}
    &U_{i, m_i}=1,\; &for\; &i=1,...,h, &\quad &U_{r,l}=0 \quad \text{otherwise}; \\
    &V_{j,n_j}=1\; &for\; &j=h+1,...,g,& \quad &V_{r,l}=0 \quad \text{otherwise},
\end{aligned}\]

In other words, let $E_{i,j}$ be the standard elementary $(g\times g)$ matrices, 
with $(i,j)$-entry 1 and all other entries 0, then we define 
\[U=\sum_{i=1}^h E_{i,m_i}, \quad V=\sum_{j=h+1}^g E_{j,n_j}.\]
Notice that $U+V$ is always a permutation matrix. 
Then we can express the action of $F$ and $V$ on $M$ and $M^{\vee}$ as: 

\[[F_M]=\left(\begin{matrix}
    U&0\\
    V&0
\end{matrix}\right),\; [V_M]= \left(\begin{matrix}
    0&0\\
    -V^t&U^t
\end{matrix}\right).\]

This description can be interpreted as follows: 
\begin{itemize}
    \item The matrix $U$ occupies the first $h$ rows, with identity columns indexed by elements of $P$, 
    and zeros elsewhere. 
    \item The matrix $V$ occupies the last $g-h$ rows, with identity columns indexed by elements of $P$, 
    placed along the second diagonal. 
\end{itemize}

\begin{example}
    If $g=4,P=\{2,3\}$, $\varphi=\{0,1,2,2\}$, then 
\[U=\left(\begin{matrix}
    0&1&0&0\\
    0&0&1&0\\
    0&0&0&0\\
    0&0&0&0
\end{matrix}\right),\; V=\left(\begin{matrix}
    0&0&0&0\\
    0&0&0&0\\
    0&0&0&1\\
    1&0&0&0
\end{matrix}\right).\]
\end{example}

Suppose $[f]\in \mathrm{Mat}_{2g}(k)$ is written in block form as: 
\[[f]=\left(\begin{matrix}
    A&B\\
    C&D
\end{matrix}\right),\]
then the conditions for Dieudonn\'e relation become: 
\begin{equation}\label{F-equation}
    \left(\begin{matrix}
        A&B\\
        C&D
    \end{matrix}\right)\left(\begin{matrix}
        U&0\\
        V&0
    \end{matrix}\right)=\left(\begin{matrix}
        0&0\\
        -V^t&U^t
    \end{matrix}\right)^{t\sigma}\left(\begin{matrix}
        A&B\\
        C&D
    \end{matrix}\right)^\sigma
\end{equation}
\[A=-A^t,D=-D^t,C=-B^t.\]

So we are reduced to solve three matrices equations: 
\[AU+BV=V C^{\sigma}, \quad CU+DV=U C^{\sigma}, \quad A=-A^t, \quad UD^\sigma=VD^\sigma=0.\]
Since $U+V$ is a permutation matrix, $D$ has to be $0$.

\begin{proposition}\label{solutions}
    The solutions to the above equations are necessarily and precisely of the following form:
    the matrix $B$ must be a $g\times g$ matrix satisfying the following conditions: 

{\rm(B1)} For $1\leqslant i\leqslant h$ and $1\leqslant j\leqslant h$, we have $B_{i,j}=B_{m_i,m_j}^p$.

{\rm(B2)} For $h+1\leqslant i\leqslant g$ and $h+1\leqslant j\leqslant g$, we have $B_{i,j}=B_{n_j,n_i}^p$.

{\rm(B3)} For $1\leqslant i\leqslant h$ and $h+1\leqslant j \leqslant g$, we have $B_{i,j}=0$. 

{\rm(B4)} For $i\in Q$, $j\in P$, we have $B_{i,j}=0$. 

For each such matrix $B$, we must have $D=0$ and $C=-B^t$, and $A$ satiefying: 

{\rm(A1)} For  $1\leqslant i\leqslant h$ and $1\leqslant j\leqslant h$, $A_{i,j}=0$.

{\rm(A2)} For $h+1\leqslant i\leqslant g$, $1\leqslant j \leqslant h$, we have $A_{i,j}=-B_{m_j,n_i}^\sigma$, 

{\rm(A3)} $A^t=-A$.
\end{proposition}

\noindent\emph{Proof}: 
We solve the system explicitly. Equation~(\ref{F-equation}) is equivalent to the following four equations: 
\[AU+BV=V C^{\sigma}, \quad CU+DV=U C^{\sigma}, \quad A=-A^t, \quad UD^\sigma=VD^\sigma=0.\]
Since $U+V$ is a permutation matrix,  it follows immediately that $D=0$. 
Now, we analyze each of the remaining equations.

First we focus on the equation $B^tU=UB^{t\sigma}$. 
$B^tU$ is the matrix where the $m_i$-th column is the $i$-th row of $B$ for $i=1,...,h$, 
and all other columns are $0$. 
$UB^{t\sigma}$ is the matrix where $i$-th row is the $m_i$-th column of $B$, 
raised to the $p$-th power, and all other rows are $0$. 
So the equation $B^tU=UB^{t\sigma}$ is equivalent to: 

(B1) for $1\leqslant i\leqslant h$ and $1\leqslant j\leqslant h$, we have $B_{i,j}=B_{m_i,m_j}^p$. 

(B3) for $1\leqslant i\leqslant h$ and $h+1\leqslant j \leqslant g$, we have $B_{i,j}=0$. 

(B4) for each $i\in Q,j \in P$, $B_{i,j}=0$. 

Now we focus on the equation $AU+BV=-VB^{t\sigma}$. 
$AU$ is the matrix with $m_i$-th column equal to the $i$-th column of $A$ for $i=1,...,h$, and all other columns are $0$. 
$BV$ is the matrix with $n_j$-th column equal to the $j$-th column of $B$ for $j=h+1,...,g$, and all other columns are $0$. 
Also, $-VB^{t\sigma}$ is the matrix where the $j$-th row equals to the negative $p$-th power of the $n_j$-th column of $B$ for $j=h+1,...,g$, 
and all other rows are $0$. 
So the equation $AU+BV=-VB^{t\sigma}$ maybe described as: 

(B2) for $h+1\leqslant i\leqslant g$ and $h+1\leqslant j\leqslant g$, we have $B_{i,j}=B_{n_j,n_i}^p$.

(B5)=(B2) the upper right $h\times (g-h)$ block of $B$ is $0$. 

(A1) the upper left $h\times h$ part of $A$ is $0$, 

(A2) the lower left $(g-h)\times h$ part of $A$ is uniquely determined by $B$. 

We have one extra condition (A3) coming from $A^t=-A$. These finishes the proof. $\square$

Let $\mathrm{H}_{\varphi}$ denote the set of matrices $B$ satisfying the conditions in Theorem \ref{solutions}. 

\begin{theorem}\label{dimH}
    The set $\mathrm{H}_{\varphi}$ is naturally 
    a direct sum of a $\F_p$-vector space and a $k$-vector space of dimension 
    \[\sum_{i=1}^k(m_i-i).\] 

    For each $B$, the corresponding space of solutions of $A$ 
    is naturally a coset of a $k$-vector space of dimension $\frac{(g-h)(g-h-1)}{2}$. 
\end{theorem}

To prove this, we define a directed graph (permitting loops) $\mathcal{G}_{\varphi}=(V_{\varphi},E_{\varphi})$ 
that encodes all the combinatorial conditions from Proposition~\ref{solutions}.

\begin{definition}\label{GGphi}
    For each Ekedahl-Oort type $\varphi$, we define the directed graph (permitting loops) $\mathcal{G}_{\varphi}=(V_{\varphi},E_{\varphi})$ as follows: 

    {\rm(0)} For each $1\leqslant i\leqslant g$, $1\leqslant j\leqslant g$, define a vertex $v_{i,j}$; 

    {\rm(1)} For each $1\leqslant i\leqslant h$, $1\leqslant j\leqslant h$, draw a directed edge from $v_{i,j}$ to $v_{m_i,m_j}$; 

    {\rm(2)} For each $h+1\leqslant s\leqslant g$, $h+1\leqslant t\leqslant g$, draw a directed edge from $v_{i,j}$ to $v_{n_j,n_i}$. 

    We mark certain vertices as zero based on the vanishing conditions:  

    {\rm(3)} For $1\leqslant i\leqslant h$ and $h+1\leqslant j \leqslant g$, mark $v_{i,j}$ as $0$. 
    
    {\rm(4)} For $i\in Q$, $j\in P$, mark $v_{i,j}$ as $0$. 
\end{definition}

Note that conditions (3) and (4) may overlap. 
The relations in Proposition~\ref{solutions} correspond to this graph in the following way: 
a relation $B_{i_1,j_1}=B_{i_2,j_2}^p$ corresponds to a direct edge from $v_{i_1,j_1}$ to $v_{i_2,j_2}$; 
a condition $B_{i,j}=0$ corresponds to marking vertex $v_{i,j}$ as $0$. 
Thus, the directed graph $\mathcal{G}_{\varphi}=(\mathcal{V}_{\varphi}, E_{\varphi})$ captures all the combinatorial structure of the solution space described in Proposition~\ref{solutions}.

\begin{example}
    Suppose $g=4,P=\{2,3\},\varphi=\{0,1,2,2\}$, then the directed graph $\mathcal{G}_{\varphi}=(\mathcal{V}_{\varphi}, E_{\varphi})$ looks like: 

\begin{center}
    \begin{tikzpicture}
        \draw [step=1cm,gray,very thin] (0,0) grid (4,4);
        \filldraw[fill=gray, fill opacity=0.3, draw=black] (2,0) rectangle (4,2);
        \filldraw[fill=orange, fill opacity=0.5, draw=black] (0,0) rectangle (1,1);
        \filldraw[fill=orange, fill opacity=0.5, draw=black] (0,3) rectangle (1,4);
        \filldraw[fill=orange, fill opacity=0.5, draw=black] (3,0) rectangle (4,1);
        \filldraw[fill=orange, fill opacity=0.5, draw=black] (3,3) rectangle (4,4);
        \filldraw[fill=gray, fill opacity=0.3, draw=black] (0,2) rectangle (2,4);
        \filldraw[fill=yellow, fill opacity=0.5, draw=black] (1,1) rectangle (3,3);
        \node at (3.5,3.5){0};
        \node at (2.5,3.5){0};
        \node at (2.5,2.5){0};
        \node at (3.5,2.5){0};
        \node at (1.5,3.5){0};
        \node at (1.5,0.5){0};
        \node at (2.5,0.5){0};
        \node at (0.5,3.5){$v_{1,1}$};
        \node at (0.5,0.5){$v_{4,1}$};
        \node at (3.5,0.5){$v_{4,4}$};
        \node at (2.5,1.5){$v_{3,3}$};
        \node at (3.5,1.5){$v_{3,4}$};
        \node at (0.5,2.5){$v_{2,1}$};
        \node at (1.5,2.5){$v_{2,2}$};
        \node at (0.5,1.5){$v_{3,1}$};
        \node at (1.5,1.5){$v_{3,2}$};
        \node at (2.5,0.5){0};
        \draw[->] (2.8,1.5) to[bend left] (3.5,0.7);
        \draw[->] (3.2,0.5) to[bend left] (0.5,3.3);
        \draw[->] (3.8,1.5) to[bend right] (3.8,3.5);
        \draw[->] (2.5,0.3) to[bend left] (0.5,0.3);
        \draw[->] (0.8,3.5) to[bend left] (1.5,2.7);
        \draw[->] (1.8,3.5) to[bend left] (2.5,2.7);
        \draw[->] (0.8,2.5) to[bend left] (1.5,1.7);
        \draw[->] (1.8,2.5) to[bend left] (2.5,1.7);
    \end{tikzpicture}
\end{center}
In this example, it is straightforward to verify that the number of free orbits is $2$. 
\end{example}
 
There are several simple but important observations about the graph $\mathcal{D}_{\varphi}=(V_{\varphi},E_{\varphi})$.

\begin{proposition}
    Let $\mathcal{D}_{\varphi}=(V_{\varphi},E_{\varphi})$ be the graph defined as above. The following statements hold: 

    {\rm(1)} For each directed edge $e\in E_{\varphi}$ from $v_{i_1,j_1}$ to $v_{i_2,j_2}$, we have $|i_1-j_1|\leqslant |i_2-j_2|$. 
    In other words, edges move progressively away from the diagonal.

    {\rm(2)} The in-degree and out-degree of every vertex are at most $1$.

    {\rm(3)} The connected components of the graph $\mathcal{D}_{\varphi}$ are of two types: either cycles: 
    \[v_{i_1,j_1}\rightarrow v_{i_2,j_2}\rightarrow... \rightarrow v_{i_n,j_n}\rightarrow v_{i_1,j_1},\] 
    or chains: 
    \[v_{i_1,j_1}\rightarrow v_{i_2,j_2}\rightarrow... \rightarrow v_{i_n,j_n}.\]

    {\rm(4)} For each $i$, the diagonal vertex $v_{i,i}$ has both in-degree and out-degree equal to 1. 
    Moreover, any edge originating from a diagonal vertex ends at another diagonal vertex.

    {\rm(5)}  The diagonal vertices form several disjoint cycle components.

    {\rm(6)} The out-degree $deg^+(v_{i,j})=0$ if and only if 
    $1\leqslant i \leqslant h < j\leqslant g$ or $1\leqslant j \leqslant h < i\leqslant g$. 

    {\rm(7)} There are exactly $h(g-h)$ chain components above (and symmetrically, below) the diagonal. 

    {\rm(8)} In each cycle component, no vertex carries a mark of $0$.  

    {\rm(9)} Each chain component above the diagonal terminates at a vertex marked with $0$.

    {\rm(10)} Each vertex below the diagonal that carries a mark of $0$ is the starting point of a chain component.

    {\rm(11)} The number of chain components without a vertex marked 0 equals $\sum_{i=1}^h(m_i-i)$. 
\end{proposition}

\noindent\emph{Proof}: (1)(2)(4)(6) follow directly from the definitions of the graph and the construction rules. 
(3) follows from (2): since all vertices have in-degree and out-degree at most $1$, 
each component must either be a cycle or a chain. (5) follows from (4): since diagonal vertices only connect to other diagonal vertices, 
and their in-degree and out-degree are both 1, they must form cycle components. 
(7) is a consequence of (1)(6): we know exactly which vertices terminate chains: 
these are in the upper-right or lower-left blocks outside the \( h \times h \) corner, giving \( h(g - h) \) such chains in each region. 
For (8), any vertex marked $0$ is either of the form \( v_{s,t} \), or \( v_{n_j,m_i} \) with \( 1 \leq i \leq h < j \leq g \). 
The former has out-degree zero, the latter has in-degree zero. Thus, such vertices cannot appear in cycles. 
(9) is consequence of (6): a chain above the diagonal ends at a vertex \( v_{i,j} \) with \( i \leq h < j \), and such a vertex is marked $0$. 
For (10),  from the definition, a vertex below the diagonal with a mark of $0$ 
is of the form \( v_{n_j, m_i} \). This vertex cannot be the target of any edge (in-degree $0$), so it must be the start of a chain.
For (11), the number of chain components that include a vertex marked $0$ corresponds bijectively to the set
\[
\{ (n_j, m_i) \mid n_j > m_i, \ n_j \in Q, \ m_i \in P \},
\]
whose size is
\[
\sum_{i=1}^h (g - h + i - m_i).
\]
Thus, the number of chain components without any mark of 0 is
\[
h(g - h) - \sum_{i=1}^h (g - h + i - m_i) = \sum_{i=1}^h (m_i - i). \qed
\]

\noindent\textbf{Proof of Theorem \ref{dimH} (1).}
Each connected component of the directed graph
$\mathcal{D}_{\varphi}$  is either a cycle component or a chain component. 
In a cycle component
\[v_{i_1,j_1}\rightarrow v_{i_2,j_2}\rightarrow... \rightarrow v_{i_n,j_n}\rightarrow v_{i_1,j_1},\] 
the corresponding variables satisfy $B_{i_1,j_1}=B_{i_2,j_2}^p$, $B_{i_2,j_2}=B_{i_3,j_3}^p$, ... ,$B_{i_n,j_n}=B_{i_1,j_1}^p$. 
Therefore $B_{i_1,j_1}=B_{i_1,j_1}^{p^n}$, 
which implies that the entries associated with the cycle form an $\F_{p^n}$-vector space. 
In particular, the solution space for such a component is isomorphic to the finite field $\F_{p^n}$.

In a chain component whose, If any vertex in the component is marked with $0$, 
then by definition, all associated entries in the matrix $B$ must be $0$. 
If no vertex in the chain is marked with $0$, 
then the entries corresponding to the chain are freely determined, 
forming a copy of the additive groups $\G_a(k)$. 
Hence, the number of chain components with no $0$-mark determines the dimension of the 
additive group component in the solution space. 
By Proposition 2.1(11), this number is $\sum_{i=1}^h (m_i - i).$ \qed

\noindent\textbf{Proof of Theorem \ref{dimH} (2).}
Consider the structure of the matrix $A$, which depends on the matrix $B$ 
according to Conditions \ref{solutions}(A1) and (A2):
the upper left $h\times h$ part of $A$ is trivial, and upper right $h\times (g-h)$ part, 
down left $(g-h)\times h$ part are both determined by $B$ uniquely. 
The down right $(g-h)\times (g-h)$ matrix has to be anti-symmetric, 
so the solution is naturally a $k$-vector space of dimension $\frac{(g-h)(g-h-1)}{2}$. \qed 

\noindent\textbf{Proof of Theorem \ref{E-O}.}
This result follows directly from combining Skorobogatov’s formula with Theorem \ref{dimH}. \qed

\section{Non-Principally Polarized Abelian Threefolds}

We also aim to address the case of abelian varieties that are not principally polarized.
At present, we do not have a complete calculation for this more general setting, and further work is required to fully understand it.
Nevertheless, we do have some partial results that allow us to compute the isogeny class of $U$ for an arbitrary abelian threefold. 
We begin by recalling Kraft's classification theorem for finite commutative group schemes killed by $p$. 

As usual, let $k$ be an algebraically closed field of characteristic $p$, and let $p\neq 2$ in this section. 
Let $\BT_1$ denote the category of $1$-truncated Barsotti-Tate groups--that is, 
finite commutative group schemes over $k$ killed by $p$ satisfying the relations 
$\Ker(F)=\IM(V)$ and $\Ker(V)=\IM(F)$. 
A $1$-truncated Barsotti-Tate group is said to be \emph{indecomposable} if it cannot be written as the direct product of two non-trivial subgroup schemes.

Recall that to every \emph{circular word} $w$, composed of letters $F$ and $V$, 
Kraft associated a $\BT_1$ group $G_w$ to $w$ by defining its covariant Dieudonné module $M_w$. 
For example, the circular word $[VVFF]$ corresponds to a Dieudonn\'e module 
with basis $\{e_1,e_2,e_3,e_4\}$ where
the operators $F$ and $V$ act according to the diagram: 
\[e_1\xleftarrow{V}e_2 \xleftarrow{V} e_3 \xrightarrow{F}e_4 \xrightarrow{F} e_1. \]

We say a word $w$ is indecomposable if it is not of the form $w=w'^d$ for any word $w'$ and $d> 1$. 
Kraft then proved the following theorem: 

\begin{theorem}{\bf (Kraft)}
    For every indecomposable circular word $w$, the $\BT_1$ group scheme $G_w$ is indecomposable. 
    Every indecomposable $\BT_1$ group scheme is of the form $G_w$ for an indecomposable circular word $w$. 
\end{theorem}

\begin{proof}
    See for example, \cite[Theorem 5.3]{Oo3}. 
\end{proof}

\begin{proposition}
    Let $A$ be a supergeneral abelian variety over $k$. Then $A[p]$ is isomorphic to the $\BT_1$ group scheme associated with the word:
    \[[VV...VFF...F].\]
\end{proposition}

\begin{proof}
    Let $M$ denote the covariant Dieudonné module of $A[p^\infty]$. 
    According to \cite[Lemma 5.7]{LNV}, $M$ is generated by a single element $z$ as a module over $W_\sigma[F,V]$, 
    and there exists a relation
    \[f=a_0F^g+a_1F^{g-1}+...+a_g+b_1V+...+b_{g-1}V^{g-1}+b_gV^g\in W_\sigma[F,V]\]
    such that $fz=0$ and $M\cong W_{\sigma}[F,V]/(f)$ with $a_0,b_g\in W$ units 
    and all other coefficients are divisible by $p$. 
    Reducing modulo $p$, the resulting module corresponds to the word $[VV...VFF...F]$. 
\end{proof}

\begin{corollary}\label{supergeneraldimUp}
    Let $A$ be any supergeneral (not necessarily principally polarized) abelian variety of dimension $g$. 
    Then $\dim U[p]=g-1$. 
\end{corollary}

\begin{proof}
    For any supergeneral abelian variety $A$ of dimension $g$, 
    $A[p]$ is always the same group scheme. 
    Since the space $\Hom(A[p],A^\vee[p])^\sym$ does not depend on the polarization, 
    the result follows from Theorem \ref{E-O}. 
\end{proof}

We now turn to supersingular abelian varieties with $a$-number $a=g-1$, i.e, 
with almost maximal \(a\)-number.  
We compute in this case, $\dim(U[p])$. 

For integers $m\geqslant 0$ and $n\geqslant 1$, let $w_{m,n}$ be the circular word consisting of $m$-copies of $FV$, followed by 
$F$, then $n$ copies of $FV$, followed by $V$:
\[[(FV)(FV)...(FV)(F)(FV)(FV)...(FV)(V)].\]
Let $G_{m,n}$ be the corresponding group scheme of rank $2(m+n+1)$. 
The Cartier dual of $G_{m,n}$ is easily seen to be isomorphic to $G_{n-1,m+1}$. 

Recall from \cite{Oo4} that Oort studied when the $p$-kernel $X[p]$ of a p-divisible group $X$ 
determines $X$ itself. 
For each Newton polygon \(\beta\) with slopes in the interval \([0,1]\), 
he introduced a \emph{minimal} p-divisible group \(\HH_{\beta}\), 
characterized by the property that any p-divisible group $X$ such that $X[p]\cong \HH_{\beta}[p]$ is 
necessarily isomorphic to $\HH_{\beta}$. 
The corresponding \(\BT_1\) group schemes \(\HH_{\beta}[p]\), 
are thus also called \emph{minimal} because of this uniquely lifting property. 
In \cite[5.3]{Oo5}, the authors described, for each pair coprime integers \(m,n\), 
the covariant Dieudonné module of the minimal p-divisible group $\HH_{m,n}$ 
(isoclinic of slope \(\frac{n}{m+n}\), of multiplicity $m+n$), 
as a free $W$-module with basis $\{e_0,e_1,...,e_{m+n-1}\}$ 
and actions $F(e_i)=e_{i+n}, V(e_i)=e_{i+m}$, where $e_{i+m+n}=pe_i$. 
For instance, $\HH_{1,1}[p]$ corresponds to the circular word $[FV]$. 
It is also true that a direct product of minimal $\BT_1$ groups schemes is again minimal. 

\begin{lemma}\label{minimalnn+1}
    The minimal $\BT_1$ group schemes $\HH_{n,n+1}[p]$ and $\HH_{n+1,n}[p]$ 
    are isomorphic, respectively, to the group schemes associated with the circular words:
    \[[(FV)(FV)...(FV)(F)], \quad [(FV)(FV)...(FV)(V)]. \]
\end{lemma}

\begin{proof}
    The covariant Dieudonné module of the minimal p-divisible group $\HH_{n,n+1}$, 
    is a free $W$-module with basis $\{e_0,e_1,...,e_{2n}\}$, 
    and \[F(e_i)=e_{i+n}, V(e_i)=e_{i+n+1} \text{where} e_{i+2n+1}=pe_i.\]
    Tracing the $F$ and $V$ actions yields
    \[e_n\xrightarrow{F} e_{2n}\xleftarrow{V}e_{n-1}\xrightarrow{F}...\xleftarrow{V}e_0\xrightarrow{F}e_n,\]
    which corresponds exactly to the circular word \([(FV)...(FV)(F)]\). 

    Similarly, the covariant Dieudonné module of the minimal p-divisible group $\HH_{n+1,n}$, 
    is a free $W$-module with basis $\{e_0,e_1,...,e_{2n}\}$, 
    and \[F(e_i)=e_{i+n+1}, V(e_i)=e_{i+n} \text{where} e_{i+2n+1}=pe_i.\] 
    In this case the action diagram is
    \[e_{0}\xrightarrow{F} e_{n+1}\xleftarrow{V}e_{1}\xrightarrow{F}...\xleftarrow{V}e_n\xleftarrow{V}e_0,\]
    which corresponds to the word \([(FV)...(FV)(V)]\). 
\end{proof}

\begin{proposition}
    Let $A$ be a supersingular abelian variety of dimension $g$ with $a$-number $a(A)=g-1$. 
    Then there exist integers $m\geqslant 0,n\geqslant1$ with $m+n\leqslant g-1$, such that $A[p]$ is isomorphic to: 
    \[[FV]^{g-1-m-n}\oplus G_{m,n}.\]
\end{proposition}

\begin{proof}
    Decompose $A[p]$ into indecomposable $BT_1$ group schemes and write down their circular words. 
    Since $a(A)=g-1$, there are exactly $g-1$ copies of consecutive occurences of 
    $FV$ among all the words of $A[p]$. 
    Thus each word can only be formed from copies of $(FV)$ with at most one additional $(F)$ and at most one additional $(V)$. 
    The only possibilities of such indecomposable words are: 
    \[[FV], [(FV)...(FV)(F)], [(FV)...(FV)(V)], G_{m,n}.\]
    Consequently, $A[p]$ must either be of the form $[FV]^{g-1-m-n}\oplus G_{m,n}$, 
    or be of the form \[[FV]^k\oplus[(FV)(FV)...(FV)(F)]\oplus[(FV)(FV)...(FV)(V)].\] 
    The second form, by Lemma \ref{minimalnn+1} corresponds to the minimal $p$-divisible group for a Newton polygon that is not supersingular.
    In this case, the lifting to a $p$-divisible group is unique up to isomorphism, and the resulting Newton polygon is not supersingular-contradicting our hypothesis on $A$. 
\end{proof} 

\begin{corollary}\label{Abe3g=2}
    Let $A$ be a supersingular abelian threefold over $k$ with $a$-number $2$. 
    Then the isomorphism class of $A[p]$ is one of: 
    \[[FFVVFV],\quad [VVFFVF],\quad [FFVV][FV].\] 
\end{corollary}

We now proceed to compute $\dim(\Hom(A[p],A^\vee[p])^{\sym})$ for a supersingular abelian variety with $a=g-1$. 

\begin{proposition}\label{HomGmnsym}
    The group \[\Hom(G_{m,n},G_{m,n}^\vee)^\sym\] has dimension $\frac{(m+n+1)(m+n)}{2}-\min\{m,n-1\}$.
\end{proposition}

\begin{proof}
    The covariant Dieudonné module of $G_{m,n}$ has the form: 
    \[e_1\xrightarrow{F}e_2\xleftarrow{V}e_3\xrightarrow{F}...\xleftarrow{V}e_{2m+1}\xrightarrow{F}e_{2m+2}\xrightarrow{F}e_{2m+3}\xleftarrow{V}...\xleftarrow{V}e_{2m+2n+2}\xleftarrow{V}e_1\]  
    The dual basis \(\{e_i^*\}\) then satisfies:
    \[e_1^*\xleftarrow{V}e_2^*\xrightarrow{F}e_3^*\xleftarrow{V}...\xrightarrow{F}e_{2m+1}^*\xleftarrow{V}e_{2m+2}^*\xleftarrow{V}e_{2m+3}^*\xrightarrow{F}...\xrightarrow{F}e_{2m+2n+2}^*\xrightarrow{F}e_1^*.\] 

    Let $x=\sum_{i=1}^{2m+2n+2}a_ie_i^*$, then a straightforward computation gives:
    \begin{equation}\label{GmnFV}
        \begin{aligned}
        F(x)&=a_{2m+2n+2}^\sigma e_1^*+\sum_{i=1}^m a_{2i}^\sigma e_{2i+1}^*+\sum_{i=m+2}^{m+n+1} a_{2i-1}^\sigma e_{2i}^*, \\
        V(x)&=\sum_{i=1}^{m+1} a_{2i}^{\sigma^{-1}}e_{2i-1}^*+\sum_{i=m+1}^{m+n} a_{2i+1}^{\sigma^{-1}}e_{2i}^*.
        \end{aligned}
    \end{equation}
    
    Let $f: G_{m,n}\rightarrow G_{m,n}^\vee$. 
    Following notations in Section \ref{SectionTheHomGroup}, 
    for $1\leqslant j\leqslant 2m+2n+2$, write 
    \[f(e_j)=\sum_{i=1}^{2m+2n+2} a_{i,j}e_i^*, \]

    We need to find the relations these coefficients $a_{i,j}$ need to satisfy. 
    These coefficients are given by a symmetric map $f: G_{m,n}\rightarrow G_{m,n}^\vee$ 
    if and only if they satisfy the following conditions: 
    \begin{itemize}
        \item For \(j=1,...,m\), we have \(F(f(e_{2j-1}))=f(e_{2j})=V(f(e_{2j+1})).\) 
        \item \(F(F(f(e_{2m+1})))=F(f(e_{2m+2}))=f(e_{2m+3})=V(f(e_{2m+4})).\)
        \item For \(j=m+2,...,m+n-1\), we have \(F(f(e_{2j}))=f(e_{2j+1})=V(f(e_{2j+2}))\).
        \item \(F(f(e_{2m+2n}))=f(e_{2m+2n+1})=V(f(e_{2m+2n+2}))=V(V(f(e_1)))\). 
        \item The matrix $[f]$ formed by $\{a_{i,j}\}$ is anti-symmetric.
    \end{itemize}

    Just as the proof of \ref{E-O}, we draw a direct graph that encodes all the combinatorial conditions. 
    We draw a $(2m+2n+2)\times(2m+2n+2)$ 
    grid with $a_{i,j}$ in the $i$-th row and $j$-th column. 
    An arrow from $a_{i_1,j_1}$ to $a_{i_2,j_2}$ indicates that
    that $a_{i_2,j_2}=a_{i_1,j_1}^\sigma$ is required in the Dieudonn\'e relations. 
    (Unlike in the previous section, here we take the target to be the Frobenius of the source for convenience.)
    An entry $a_{i,j}$ is marked $0$ if and only if $a_{i,j}=0$ 
    is required in the Dieudonné relations. 

    Define: \[P:= \{2j|k=1,2,...,m+1\}\cup \{2j+1|j=m+1,...,m+n\},\]
    \[Q:=\{2j|k=1,2,...,m\}\cup \{2j+1|j=m+1,...,m+n\}\cup \{2m+2n+2\}.\]

    Then: 
    \begin{itemize}
        \item $\{e_i|i\in P\}$ spans $\IM{(F)}$ of the covariant Dieudonné module of $G_{m,n}$. 
        \item $\{e_i^*|i \notin P\}$ spans $\IM{F}$ of the covariant Dieudonné module of $G_{m,n}^\vee$. 
        \item $\{e_i|i\in Q\}$ spans $\IM{(V)}$ of the covariant Dieudonné module of $G_{m,n}$. 
        \item $\{e_i^*|i \notin Q\}$ spans $\IM{(V)}$ of the covariant Dieudonné module of $G_{m,n}^\vee$. 
    \end{itemize}
    We have \[\#P=\#Q=m+n+1, \#(P\cap Q)=m+n.\]
   
The Dieudonné module structure and the symmetry conditions then imply:  
\begin{itemize}
    \item All diagonal entries vanish, since we require the matrix to be anti-symmetric. 
    \item Each arrow indicates that the target equals the Frobenius of the source.
    \item For any $j\in P$, since $f(e_j)=F(f(e_{j-1}))$, formula \ref{GmnFV} forces: for any $i\in P$, $a_{i,j}=0$; 
    for any $i \notin P$, we draw an arrow from $a_{i-1,j-1}$ to $a_{i,j}$. 
    \item For any $j\in Q$, since $f(e_j)=V(f(e_{j+1}))$, formula \ref{GmnFV} forces: for any $i\in Q$, $a_{i,j}=0$;
    for each $i\notin Q$, we draw an arrow from $a_{ij}$ to $a_{i+1,j+1}$. 
    \item The matrix has to be anti-symmetric. 
\end{itemize}

    Let us first settle the case $n-1\leqslant m$. 
    We take $m=4$ and $n=3$ to illustrate. 
    We show that: in this case, under the above conditions, there are exactly
    \[\frac{(m+1)m}{2}+(m+1)(n-1)+(m+1)+\frac{(n-1)(n-2)}{2}=\frac{(m+n+1)(m+n)}{2}-(n-1)\]
    degrees of freedom for the coefficients $\{a_{i,j}\}$.

    \begin{center}
        \begin{tikzpicture}
            \draw [step=1cm,gray,very thin] (0,0) grid (16,16);
            \filldraw[fill=gray, fill opacity=0.3, draw=black] (0,7) rectangle (9,16);
            \filldraw[fill=gray, fill opacity=0.3, draw=black] (10,1) rectangle (15,6);
            \filldraw[fill=gray, fill opacity=0.3, draw=black] (0,1) rectangle (9,6);
            \filldraw[fill=gray, fill opacity=0.3, draw=black] (10,7) rectangle (15,16);
            \node at (0.5,15.5){0};
            \node at (1.5,14.5){0};
            \node at (1.5,12.5){0};
            \node at (1.5,10.5){0};
            \node at (1.5,8.5){0};
            \node at (1.5,6.5){0};
            \node at (1.5,5.5){0};
            \node at (1.5,3.5){0};
            \node at (1.5,1.5){0};
            \node at (1.5,0.5){0};
            \node at (3.5,14.5){0};
            \node at (3.5,12.5){0};
            \node at (3.5,10.5){0};
            \node at (3.5,8.5){0};
            \node at (3.5,6.5){0};
            \node at (3.5,5.5){0};
            \node at (3.5,3.5){0};
            \node at (3.5,1.5){0};
            \node at (3.5,0.5){0};
            \node at (5.5,14.5){0};
            \node at (5.5,12.5){0};
            \node at (5.5,10.5){0};
            \node at (5.5,8.5){0};
            \node at (5.5,6.5){0};
            \node at (5.5,5.5){0};
            \node at (5.5,3.5){0};
            \node at (5.5,1.5){0};
            \node at (5.5,0.5){0};
            \node at (7.5,14.5){0};
            \node at (7.5,12.5){0};
            \node at (7.5,10.5){0};
            \node at (7.5,8.5){0};
            \node at (7.5,6.5){0};
            \node at (7.5,5.5){0};
            \node at (7.5,3.5){0};
            \node at (7.5,1.5){0};
            \node at (7.5,0.5){0};
            \node at (9.5,14.5){0};
            \node at (9.5,12.5){0};
            \node at (9.5,10.5){0};
            \node at (9.5,8.5){0};
            \node at (9.5,6.5){0};
            \node at (9.5,5.5){0};
            \node at (9.5,3.5){0};
            \node at (9.5,1.5){0};
            \node at (10.5,14.5){0};
            \node at (10.5,12.5){0};
            \node at (10.5,10.5){0};
            \node at (10.5,8.5){0};
            \node at (10.5,6.5){0};
            \node at (10.5,5.5){0};
            \node at (10.5,3.5){0};
            \node at (10.5,1.5){0};
            \node at (10.5,0.5){0};
            \node at (12.5,14.5){0};
            \node at (12.5,12.5){0};
            \node at (12.5,10.5){0};
            \node at (12.5,8.5){0};
            \node at (12.5,6.5){0};
            \node at (12.5,5.5){0};
            \node at (12.5,3.5){0};
            \node at (12.5,1.5){0};
            \node at (12.5,0.5){0};
            \node at (14.5,14.5){0};
            \node at (14.5,12.5){0};
            \node at (14.5,10.5){0};
            \node at (14.5,8.5){0};
            \node at (14.5,6.5){0};
            \node at (14.5,5.5){0};
            \node at (14.5,3.5){0};
            \node at (14.5,1.5){0};
            \node at (14.5,0.5){0};
            \node at (15.5,14.5){0};
            \node at (15.5,12.5){0};
            \node at (15.5,10.5){0};
            \node at (15.5,8.5){0};
            \node at (15.5,5.5){0};
            \node at (15.5,3.5){0};
            \node at (15.5,1.5){0};
            \node at (15.5,0.5){0};

            \node at (2.5,13.5){0};
            \node at (4.5,11.5){0};
            \node at (6.5,9.5){0};
            \node at (8.5,7.5){0};
            \node at (11.5,4.5){0};
            \node at (13.5,2.5){0};

            \filldraw[fill=orange, fill opacity=0.3, draw=black] (0,13) rectangle (1,14);
            \filldraw[fill=orange, fill opacity=0.3, draw=black] (0,11) rectangle (1,12);
            \filldraw[fill=orange, fill opacity=0.3, draw=black] (0,9) rectangle (1,10);
            \filldraw[fill=orange, fill opacity=0.3, draw=black] (0,7) rectangle (1,8);

            \filldraw[fill=orange, fill opacity=0.3, draw=black] (2,11) rectangle (3,12);
            \filldraw[fill=orange, fill opacity=0.3, draw=black] (2,9) rectangle (3,10);
            \filldraw[fill=orange, fill opacity=0.3, draw=black] (2,7) rectangle (3,8);

            \filldraw[fill=orange, fill opacity=0.3, draw=black] (4,9) rectangle (5,10);
            \filldraw[fill=orange, fill opacity=0.3, draw=black] (4,7) rectangle (5,8);

            \filldraw[fill=orange, fill opacity=0.3, draw=black] (6,7) rectangle (7,8);

            \filldraw[fill=orange, fill opacity=0.3, draw=black] (2,15) rectangle (3,16);
            \filldraw[fill=orange, fill opacity=0.3, draw=black] (4,15) rectangle (5,16);
            \filldraw[fill=orange, fill opacity=0.3, draw=black] (6,15) rectangle (7,16);
            \filldraw[fill=orange, fill opacity=0.3, draw=black] (8,15) rectangle (9,16);

            \filldraw[fill=orange, fill opacity=0.3, draw=black] (4,13) rectangle (5,14);
            \filldraw[fill=orange, fill opacity=0.3, draw=black] (6,13) rectangle (7,14);
            \filldraw[fill=orange, fill opacity=0.3, draw=black] (8,13) rectangle (9,14);

            \filldraw[fill=orange, fill opacity=0.3, draw=black] (6,11) rectangle (7,12);
            \filldraw[fill=orange, fill opacity=0.3, draw=black] (8,11) rectangle (9,12);

            \filldraw[fill=orange, fill opacity=0.3, draw=black] (8,9) rectangle (9,10);

            \filldraw[fill=orange, fill opacity=0.3, draw=black] (0,2) rectangle (1,3);
            \filldraw[fill=orange, fill opacity=0.3, draw=black] (2,2) rectangle (3,3);
            \filldraw[fill=orange, fill opacity=0.3, draw=black] (4,2) rectangle (5,3);
            \filldraw[fill=orange, fill opacity=0.3, draw=black] (6,2) rectangle (7,3);
            \filldraw[fill=orange, fill opacity=0.3, draw=black] (8,2) rectangle (9,3);

            \filldraw[fill=orange, fill opacity=0.3, draw=black] (0,4) rectangle (1,5);
            \filldraw[fill=orange, fill opacity=0.3, draw=black] (2,4) rectangle (3,5);
            \filldraw[fill=orange, fill opacity=0.3, draw=black] (4,4) rectangle (5,5);
            \filldraw[fill=orange, fill opacity=0.3, draw=black] (6,4) rectangle (7,5);
            \filldraw[fill=orange, fill opacity=0.3, draw=black] (8,4) rectangle (9,5);

            \filldraw[fill=orange, fill opacity=0.3, draw=black] (11,7) rectangle (12,8);
            \filldraw[fill=orange, fill opacity=0.3, draw=black] (11,9) rectangle (12,10);
            \filldraw[fill=orange, fill opacity=0.3, draw=black] (11,11) rectangle (12,12);
            \filldraw[fill=orange, fill opacity=0.3, draw=black] (11,13) rectangle (12,14);
            \filldraw[fill=orange, fill opacity=0.3, draw=black] (11,15) rectangle (12,16);

            \filldraw[fill=orange, fill opacity=0.3, draw=black] (13,7) rectangle (14,8);
            \filldraw[fill=orange, fill opacity=0.3, draw=black] (13,9) rectangle (14,10);
            \filldraw[fill=orange, fill opacity=0.3, draw=black] (13,11) rectangle (14,12);
            \filldraw[fill=orange, fill opacity=0.3, draw=black] (13,13) rectangle (14,14);
            \filldraw[fill=orange, fill opacity=0.3, draw=black] (13,15) rectangle (14,16);

            \filldraw[fill=orange, fill opacity=0.3, draw=black] (11,2) rectangle (12,3);
            
            \filldraw[fill=orange, fill opacity=0.3, draw=black] (13,4) rectangle (14,5);

            \draw[->] (-0.3,15.5) to[bend left] (0.5,14.7);
            \draw[->] (-0.3,13.5) to[bend left] (0.5,12.7);
            \draw[->] (-0.3,11.5) to[bend left] (0.5,10.7);
            \draw[->] (-0.3,9.5) to[bend left] (0.5,8.7);
            \draw[->] (-0.3,7.5) to[bend left] (0.5,6.7);
            \draw[->] (-0.3,6.5) to[bend left] (0.5,5.7);
            \draw[->] (-0.3,4.5) to[bend left] (0.5,3.7);
            \draw[->] (-0.3,2.5) to[bend left] (0.5,1.7);

            \draw[->] (0.7,16.5) to[bend left] (1.5,15.7);
            \draw[->] (0.7,14.5) to[bend left] (1.5,13.7);
            \draw[->] (0.7,12.5) to[bend left] (1.5,11.7);
            \draw[->] (0.7,10.5) to[bend left] (1.5,9.7);
            \draw[->] (0.7,8.5) to[bend left] (1.5,7.7);
            \draw[->] (0.7,5.5) to[bend left] (1.5,4.7);
            \draw[->] (0.7,3.5) to[bend left] (1.5,2.7);
            \draw[->] (0.7,1.5) to[bend left] (1.5,0.7);
            \draw[->] (0.7,0.5) to[bend left] (1.5,-0.3);

            \draw[->] (1.7,15.5) to[bend left] (2.5,14.7);
            \draw[->] (1.7,13.5) to[bend left] (2.5,12.7);
            \draw[->] (1.7,11.5) to[bend left] (2.5,10.7);
            \draw[->] (1.7,9.5) to[bend left] (2.5,8.7);
            \draw[->] (1.7,7.5) to[bend left] (2.5,6.7);
            \draw[->] (1.7,6.5) to[bend left] (2.5,5.7);
            \draw[->] (1.7,4.5) to[bend left] (2.5,3.7);
            \draw[->] (1.7,2.5) to[bend left] (2.5,1.7);

            \draw[->] (2.7,16.5) to[bend left] (3.5,15.7);
            \draw[->] (2.7,14.5) to[bend left] (3.5,13.7);
            \draw[->] (2.7,12.5) to[bend left] (3.5,11.7);
            \draw[->] (2.7,10.5) to[bend left] (3.5,9.7);
            \draw[->] (2.7,8.5) to[bend left] (3.5,7.7);
            \draw[->] (2.7,5.5) to[bend left] (3.5,4.7);
            \draw[->] (2.7,3.5) to[bend left] (3.5,2.7);
            \draw[->] (2.7,1.5) to[bend left] (3.5,0.7);
            \draw[->] (2.7,0.5) to[bend left] (3.5,-0.3);

            \draw[->] (3.7,15.5) to[bend left] (4.5,14.7);
            \draw[->] (3.7,13.5) to[bend left] (4.5,12.7);
            \draw[->] (3.7,11.5) to[bend left] (4.5,10.7);
            \draw[->] (3.7,9.5) to[bend left] (4.5,8.7);
            \draw[->] (3.7,7.5) to[bend left] (4.5,6.7);
            \draw[->] (3.7,6.5) to[bend left] (4.5,5.7);
            \draw[->] (3.7,4.5) to[bend left] (4.5,3.7);
            \draw[->] (3.7,2.5) to[bend left] (4.5,1.7);

            \draw[->] (4.7,16.5) to[bend left] (5.5,15.7);
            \draw[->] (4.7,14.5) to[bend left] (5.5,13.7);
            \draw[->] (4.7,12.5) to[bend left] (5.5,11.7);
            \draw[->] (4.7,10.5) to[bend left] (5.5,9.7);
            \draw[->] (4.7,8.5) to[bend left] (5.5,7.7);
            \draw[->] (4.7,5.5) to[bend left] (5.5,4.7);
            \draw[->] (4.7,3.5) to[bend left] (5.5,2.7);
            \draw[->] (4.7,1.5) to[bend left] (5.5,0.7);
            \draw[->] (4.7,0.5) to[bend left] (5.5,-0.3);

            \draw[->] (5.7,15.5) to[bend left] (6.5,14.7);
            \draw[->] (5.7,13.5) to[bend left] (6.5,12.7);
            \draw[->] (5.7,11.5) to[bend left] (6.5,10.7);
            \draw[->] (5.7,9.5) to[bend left] (6.5,8.7);
            \draw[->] (5.7,7.5) to[bend left] (6.5,6.7);
            \draw[->] (5.7,6.5) to[bend left] (6.5,5.7);
            \draw[->] (5.7,4.5) to[bend left] (6.5,3.7);
            \draw[->] (5.7,2.5) to[bend left] (6.5,1.7);

            \draw[->] (6.7,16.5) to[bend left] (7.5,15.7);
            \draw[->] (6.7,14.5) to[bend left] (7.5,13.7);
            \draw[->] (6.7,12.5) to[bend left] (7.5,11.7);
            \draw[->] (6.7,10.5) to[bend left] (7.5,9.7);
            \draw[->] (6.7,8.5) to[bend left] (7.5,7.7);
            \draw[->] (6.7,5.5) to[bend left] (7.5,4.7);
            \draw[->] (6.7,3.5) to[bend left] (7.5,2.7);
            \draw[->] (6.7,1.5) to[bend left] (7.5,0.7);
            \draw[->] (6.7,0.5) to[bend left] (7.5,-0.3);

            \draw[->] (7.7,15.5) to[bend left] (8.5,14.7);
            \draw[->] (7.7,13.5) to[bend left] (8.5,12.7);
            \draw[->] (7.7,11.5) to[bend left] (8.5,10.7);
            \draw[->] (7.7,9.5) to[bend left] (8.5,8.7);
            \draw[->] (7.7,7.5) to[bend left] (8.5,6.7);
            \draw[->] (7.7,6.5) to[bend left] (8.5,5.7);
            \draw[->] (7.7,4.5) to[bend left] (8.5,3.7);
            \draw[->] (7.7,2.5) to[bend left] (8.5,1.7);

            \draw[->] (8.7,16.5) to[bend left] (9.5,15.7);
            \draw[->] (8.7,14.5) to[bend left] (9.5,13.7);
            \draw[->] (8.7,12.5) to[bend left] (9.5,11.7);
            \draw[->] (8.7,10.5) to[bend left] (9.5,9.7);
            \draw[->] (8.7,8.5) to[bend left] (9.5,7.7);
            \draw[->] (8.7,5.5) to[bend left] (9.5,4.7);
            \draw[->] (8.7,3.5) to[bend left] (9.5,2.7);
            \draw[->] (8.7,1.5) to[bend left] (9.5,0.7);
            \draw[->] (8.7,0.5) to[bend left] (9.5,-0.3);

            \draw[->] (9.7,16.5) to[bend left] (10.5,15.7);
            \draw[->] (9.7,14.5) to[bend left] (10.5,13.7);
            \draw[->] (9.7,12.5) to[bend left] (10.5,11.7);
            \draw[->] (9.7,10.5) to[bend left] (10.5,9.7);
            \draw[->] (9.7,8.5) to[bend left] (10.5,7.7);
            \draw[->] (9.7,5.5) to[bend left] (10.5,4.7);
            \draw[->] (9.7,3.5) to[bend left] (10.5,2.7);
            \draw[->] (9.7,1.5) to[bend left] (10.5,0.7);
            \draw[->] (9.7,0.5) to[bend left] (10.5,-0.3);

            \draw[->] (10.7,15.5) to[bend left] (11.5,14.7);
            \draw[->] (10.7,13.5) to[bend left] (11.5,12.7);
            \draw[->] (10.7,11.5) to[bend left] (11.5,10.7);
            \draw[->] (10.7,9.5) to[bend left] (11.5,8.7);
            \draw[->] (10.7,7.5) to[bend left] (11.5,6.7);
            \draw[->] (10.7,6.5) to[bend left] (11.5,5.7);
            \draw[->] (10.7,4.5) to[bend left] (11.5,3.7);
            \draw[->] (10.7,2.5) to[bend left] (11.5,1.7);

            \draw[->] (11.7,16.5) to[bend left] (12.5,15.7);
            \draw[->] (11.7,14.5) to[bend left] (12.5,13.7);
            \draw[->] (11.7,12.5) to[bend left] (12.5,11.7);
            \draw[->] (11.7,10.5) to[bend left] (12.5,9.7);
            \draw[->] (11.7,8.5) to[bend left] (12.5,7.7);
            \draw[->] (11.7,5.5) to[bend left] (12.5,4.7);
            \draw[->] (11.7,3.5) to[bend left] (12.5,2.7);
            \draw[->] (11.7,1.5) to[bend left] (12.5,0.7);
            \draw[->] (11.7,0.5) to[bend left] (12.5,-0.3);

            \draw[->] (12.7,15.5) to[bend left] (13.5,14.7);
            \draw[->] (12.7,13.5) to[bend left] (13.5,12.7);
            \draw[->] (12.7,11.5) to[bend left] (13.5,10.7);
            \draw[->] (12.7,9.5) to[bend left] (13.5,8.7);
            \draw[->] (12.7,7.5) to[bend left] (13.5,6.7);
            \draw[->] (12.7,6.5) to[bend left] (13.5,5.7);
            \draw[->] (12.7,4.5) to[bend left] (13.5,3.7);
            \draw[->] (12.7,2.5) to[bend left] (13.5,1.7);

            \draw[->] (13.7,16.5) to[bend left] (14.5,15.7);
            \draw[->] (13.7,14.5) to[bend left] (14.5,13.7);
            \draw[->] (13.7,12.5) to[bend left] (14.5,11.7);
            \draw[->] (13.7,10.5) to[bend left] (14.5,9.7);
            \draw[->] (13.7,8.5) to[bend left] (14.5,7.7);
            \draw[->] (13.7,5.5) to[bend left] (14.5,4.7);
            \draw[->] (13.7,3.5) to[bend left] (14.5,2.7);
            \draw[->] (13.7,1.5) to[bend left] (14.5,0.7);
            \draw[->] (13.7,0.5) to[bend left] (14.5,-0.3);

            \draw[->] (14.7,15.5) to[bend left] (15.5,14.7);
            \draw[->] (14.7,13.5) to[bend left] (15.5,12.7);
            \draw[->] (14.7,11.5) to[bend left] (15.5,10.7);
            \draw[->] (14.7,9.5) to[bend left] (15.5,8.7);
            \draw[->] (14.7,7.5) to[bend left] (15.5,6.7);
            \draw[->] (14.7,6.5) to[bend left] (15.5,5.7);
            \draw[->] (14.7,4.5) to[bend left] (15.5,3.7);
            \draw[->] (14.7,2.5) to[bend left] (15.5,1.7);

            \draw[->] (15.7,15.5) to[bend left] (16.5,14.7);
            \draw[->] (15.7,13.5) to[bend left] (16.5,12.7);
            \draw[->] (15.7,11.5) to[bend left] (16.5,10.7);
            \draw[->] (15.7,9.5) to[bend left] (16.5,8.7);
            \draw[->] (15.7,7.5) to[bend left] (16.5,6.7);
            \draw[->] (15.7,6.5) to[bend left] (16.5,5.7);
            \draw[->] (15.7,4.5) to[bend left] (16.5,3.7);
            \draw[->] (15.7,2.5) to[bend left] (16.5,1.7);

            \node at (0.5,14.5){$a_{2,1}$};
        \end{tikzpicture}
    \end{center}

    To understand the relations and the $0$'s in this grid, 
    we consider the following four different regions (the gray areas in the picture): 
    \begin{itemize}
        \item $S_1$, entries $a_{ij}$ such that $1\leqslant i,j\leqslant 2m+1$;
        \item  $S_2$, entries $a_{ij}$ such that $2m+3\leqslant i\leqslant 2m+2n+1, \quad 1\leqslant j\leqslant 2m+1$;
        \item  $S_3$, entries $a_{ij}$ such that $1\leqslant i\leqslant 2m+1, \quad 2m+3\leqslant j\leqslant 2m+2n+1$;
        \item  $S_4$, entries $a_{ij}$ such that $ 2m+3\leqslant i,j\leqslant 2m+2n+1$.
    \end{itemize}

    Within these regions, some coefficients have neither a marked $0$, nor any incoming or outgoing arrows. 
    We mark these coefficients with orange boxes. 
    They represent variables with no constraints. 
    There are exactly \[\frac{m(m+1)}{2}+(m+1)(n-1)+\frac{(n-1)(n-2)}{2}\] pairs (under symmetric along the diagonal) of such boxes. 

    When $m>n-1$, there are exactly $m+1$ symmetric pairs of free orbits, 
    namely those passing though $a_{2i,1}$, for $i=1,2,...,m+1$. 
    All other orbits, either start or end at a $0$. 
    Each pair of free orbits corresponds to one copy of $\G_a$.
    Counting the total dimension, we have exactly
    \[\frac{(m+1)m}{2}+(m+1)(n-1)+\frac{(n-1)(n-2)}{2}+(m+1)=\frac{(m+n+1)(m+n)}{2}-(n-1)\]
    degrees of freedom under the symmetric condition and the Dieudonn\'e conditions. 
    
    When $m=n-1$, we illustrate the picture with the example $m=2$, $n=3$. 
    We again find exactly $m+1$ pairs of free orbits, 
    those passing though $a_{2i,1}$, for $i=1,2,...,m+1$.
    In addition, there is one extra orbit (passing through $a_{2m+2,1}$)forming a cycle, 
    which does not contribute to unipotent dimension. 
    The same counting formula applies, giving 
    \[\frac{(m+1)m}{2}+(m+1)(n-1)+\frac{(n-1)(n-2)}{2}+(m+1)=\frac{(m+n+1)(m+n)}{2}-(n-1)\]
    degrees of freedom on the variables under the symmetric condition and the Dieudonn\'e conditions. 
    
    \begin{center}
        \begin{tikzpicture}
            \draw [step=1cm,gray,very thin] (0,0) grid (12,12);
            \filldraw[fill=gray, fill opacity=0.3, draw=black] (0,7) rectangle (5,12);
            \filldraw[fill=gray, fill opacity=0.3, draw=black] (6,1) rectangle (11,6);
            \filldraw[fill=gray, fill opacity=0.3, draw=black] (0,1) rectangle (5,6);
            \filldraw[fill=gray, fill opacity=0.3, draw=black] (6,7) rectangle (11,12);
            \node at (0.5,11.5){0};

            \node at (1.5,10.5){0};
            \node at (1.5,8.5){0};
            \node at (1.5,6.5){0};
            \node at (1.5,5.5){0};
            \node at (1.5,3.5){0};
            \node at (1.5,1.5){0};
            \node at (1.5,0.5){0};

            \node at (3.5,10.5){0};
            \node at (3.5,8.5){0};
            \node at (3.5,6.5){0};
            \node at (3.5,5.5){0};
            \node at (3.5,3.5){0};
            \node at (3.5,1.5){0};
            \node at (3.5,0.5){0};

            \node at (5.5,10.5){0};
            \node at (5.5,8.5){0};
            \node at (5.5,6.5){0};
            \node at (5.5,5.5){0};
            \node at (5.5,3.5){0};
            \node at (5.5,1.5){0};

            \node at (6.5,10.5){0};
            \node at (6.5,8.5){0};
            \node at (6.5,6.5){0};
            \node at (6.5,5.5){0};
            \node at (6.5,3.5){0};
            \node at (6.5,1.5){0};
            \node at (6.5,0.5){0};

            \node at (8.5,10.5){0};
            \node at (8.5,8.5){0};
            \node at (8.5,6.5){0};
            \node at (8.5,5.5){0};
            \node at (8.5,3.5){0};
            \node at (8.5,1.5){0};
            \node at (8.5,0.5){0};

            \node at (10.5,10.5){0};
            \node at (10.5,8.5){0};
            \node at (10.5,6.5){0};
            \node at (10.5,5.5){0};
            \node at (10.5,3.5){0};
            \node at (10.5,1.5){0};
            \node at (10.5,0.5){0};

            \node at (11.5,10.5){0};
            \node at (11.5,8.5){0};
            \node at (11.5,5.5){0};
            \node at (11.5,3.5){0};
            \node at (11.5,1.5){0};
            \node at (11.5,0.5){0};

            \node at (2.5,9.5){0};
            \node at (4.5,7.5){0};
            \node at (7.5,4.5){0};
            \node at (9.5,2.5){0};

            \filldraw[fill=orange, fill opacity=0.3, draw=black] (0,9) rectangle (1,10);
            \filldraw[fill=orange, fill opacity=0.3, draw=black] (0,7) rectangle (1,8);

            \filldraw[fill=orange, fill opacity=0.3, draw=black] (2,7) rectangle (3,8);

            \filldraw[fill=orange, fill opacity=0.3, draw=black] (2,11) rectangle (3,12);
            \filldraw[fill=orange, fill opacity=0.3, draw=black] (4,11) rectangle (5,12);

            \filldraw[fill=orange, fill opacity=0.3, draw=black] (4,9) rectangle (5,10);

            \filldraw[fill=orange, fill opacity=0.3, draw=black] (0,2) rectangle (1,3);
            \filldraw[fill=orange, fill opacity=0.3, draw=black] (2,2) rectangle (3,3);
            \filldraw[fill=orange, fill opacity=0.3, draw=black] (4,2) rectangle (5,3);

            \filldraw[fill=orange, fill opacity=0.3, draw=black] (0,4) rectangle (1,5);
            \filldraw[fill=orange, fill opacity=0.3, draw=black] (2,4) rectangle (3,5);
            \filldraw[fill=orange, fill opacity=0.3, draw=black] (4,4) rectangle (5,5);

            \filldraw[fill=orange, fill opacity=0.3, draw=black] (7,7) rectangle (8,8);
            \filldraw[fill=orange, fill opacity=0.3, draw=black] (7,9) rectangle (8,10);
            \filldraw[fill=orange, fill opacity=0.3, draw=black] (7,11) rectangle (8,12);

            \filldraw[fill=orange, fill opacity=0.3, draw=black] (9,7) rectangle (10,8);
            \filldraw[fill=orange, fill opacity=0.3, draw=black] (9,9) rectangle (10,10);
            \filldraw[fill=orange, fill opacity=0.3, draw=black] (9,11) rectangle (10,12);

            \filldraw[fill=orange, fill opacity=0.3, draw=black] (7,2) rectangle (8,3);
            
            \filldraw[fill=orange, fill opacity=0.3, draw=black] (9,4) rectangle (10,5);

            \draw[->] (-0.3,11.5) to[bend left] (0.5,10.7);
            \draw[->] (-0.3,9.5) to[bend left] (0.5,8.7);
            \draw[->] (-0.3,7.5) to[bend left] (0.5,6.7);
            \draw[->] (-0.3,6.5) to[bend left] (0.5,5.7);
            \draw[->] (-0.3,4.5) to[bend left] (0.5,3.7);
            \draw[->] (-0.3,2.5) to[bend left] (0.5,1.7);

            \draw[->] (0.7,12.5) to[bend left] (1.5,11.7);
            \draw[->] (0.7,10.5) to[bend left] (1.5,9.7);
            \draw[->] (0.7,8.5) to[bend left] (1.5,7.7);
            \draw[->] (0.7,5.5) to[bend left] (1.5,4.7);
            \draw[->] (0.7,3.5) to[bend left] (1.5,2.7);
            \draw[->] (0.7,1.5) to[bend left] (1.5,0.7);
            \draw[->] (0.7,0.5) to[bend left] (1.5,-0.3);

            \draw[->] (1.7,11.5) to[bend left] (2.5,10.7);
            \draw[->] (1.7,9.5) to[bend left] (2.5,8.7);
            \draw[->] (1.7,7.5) to[bend left] (2.5,6.7);
            \draw[->] (1.7,6.5) to[bend left] (2.5,5.7);
            \draw[->] (1.7,4.5) to[bend left] (2.5,3.7);
            \draw[->] (1.7,2.5) to[bend left] (2.5,1.7);

            \draw[->] (2.7,12.5) to[bend left] (3.5,11.7);
            \draw[->] (2.7,10.5) to[bend left] (3.5,9.7);
            \draw[->] (2.7,8.5) to[bend left] (3.5,7.7);
            \draw[->] (2.7,5.5) to[bend left] (3.5,4.7);
            \draw[->] (2.7,3.5) to[bend left] (3.5,2.7);
            \draw[->] (2.7,1.5) to[bend left] (3.5,0.7);
            \draw[->] (2.7,0.5) to[bend left] (3.5,-0.3);

            \draw[->] (3.7,11.5) to[bend left] (4.5,10.7);
            \draw[->] (3.7,9.5) to[bend left] (4.5,8.7);
            \draw[->] (3.7,7.5) to[bend left] (4.5,6.7);
            \draw[->] (3.7,6.5) to[bend left] (4.5,5.7);
            \draw[->] (3.7,4.5) to[bend left] (4.5,3.7);
            \draw[->] (3.7,2.5) to[bend left] (4.5,1.7);

            \draw[->] (4.7,12.5) to[bend left] (5.5,11.7);
            \draw[->] (4.7,10.5) to[bend left] (5.5,9.7);
            \draw[->] (4.7,8.5) to[bend left] (5.5,7.7);
            \draw[->] (4.7,5.5) to[bend left] (5.5,4.7);
            \draw[->] (4.7,3.5) to[bend left] (5.5,2.7);
            \draw[->] (4.7,1.5) to[bend left] (5.5,0.7);
            \draw[->] (4.7,0.5) to[bend left] (5.5,-0.3);

            \draw[->] (5.7,12.5) to[bend left] (6.5,11.7);
            \draw[->] (5.7,10.5) to[bend left] (6.5,9.7);
            \draw[->] (5.7,8.5) to[bend left] (6.5,7.7);
            \draw[->] (5.7,5.5) to[bend left] (6.5,4.7);
            \draw[->] (5.7,3.5) to[bend left] (6.5,2.7);
            \draw[->] (5.7,1.5) to[bend left] (6.5,0.7);
            \draw[->] (5.7,0.5) to[bend left] (6.5,-0.3);

            \draw[->] (6.7,11.5) to[bend left] (7.5,10.7);
            \draw[->] (6.7,9.5) to[bend left] (7.5,8.7);
            \draw[->] (6.7,7.5) to[bend left] (7.5,6.7);
            \draw[->] (6.7,6.5) to[bend left] (7.5,5.7);
            \draw[->] (6.7,4.5) to[bend left] (7.5,3.7);
            \draw[->] (6.7,2.5) to[bend left] (7.5,1.7);

            \draw[->] (7.7,12.5) to[bend left] (8.5,11.7);
            \draw[->] (7.7,10.5) to[bend left] (8.5,9.7);
            \draw[->] (7.7,8.5) to[bend left] (8.5,7.7);
            \draw[->] (7.7,5.5) to[bend left] (8.5,4.7);
            \draw[->] (7.7,3.5) to[bend left] (8.5,2.7);
            \draw[->] (7.7,1.5) to[bend left] (8.5,0.7);
            \draw[->] (7.7,0.5) to[bend left] (8.5,-0.3);

            \draw[->] (8.7,11.5) to[bend left] (9.5,10.7);
            \draw[->] (8.7,9.5) to[bend left] (9.5,8.7);
            \draw[->] (8.7,7.5) to[bend left] (9.5,6.7);
            \draw[->] (8.7,6.5) to[bend left] (9.5,5.7);
            \draw[->] (8.7,4.5) to[bend left] (9.5,3.7);
            \draw[->] (8.7,2.5) to[bend left] (9.5,1.7);

            \draw[->] (9.7,12.5) to[bend left] (10.5,11.7);
            \draw[->] (9.7,10.5) to[bend left] (10.5,9.7);
            \draw[->] (9.7,8.5) to[bend left] (10.5,7.7);
            \draw[->] (9.7,5.5) to[bend left] (10.5,4.7);
            \draw[->] (9.7,3.5) to[bend left] (10.5,2.7);
            \draw[->] (9.7,1.5) to[bend left] (10.5,0.7);
            \draw[->] (9.7,0.5) to[bend left] (10.5,-0.3);

            \draw[->] (10.7,11.5) to[bend left] (11.5,10.7);
            \draw[->] (10.7,9.5) to[bend left] (11.5,8.7);
            \draw[->] (10.7,7.5) to[bend left] (11.5,6.7);
            \draw[->] (10.7,6.5) to[bend left] (11.5,5.7);
            \draw[->] (10.7,4.5) to[bend left] (11.5,3.7);
            \draw[->] (10.7,2.5) to[bend left] (11.5,1.7);

            \draw[->] (11.7,11.5) to[bend left] (12.5,10.7);
            \draw[->] (11.7,9.5) to[bend left] (12.5,8.7);
            \draw[->] (11.7,7.5) to[bend left] (12.5,6.7);
            \draw[->] (11.7,6.5) to[bend left] (12.5,5.7);
            \draw[->] (11.7,4.5) to[bend left] (12.5,3.7);
            \draw[->] (11.7,2.5) to[bend left] (12.5,1.7);

            \node at (0.5,10.5){$a_{2,1}$};
        \end{tikzpicture}
    \end{center}

    When $n>m+1$, let us illustrate the idea using 
    $m=2$ and $n=5$. 
    From the picture below, we observe the following: 
    there are still  \[\frac{m(m+1)}{2}+(m+1)(n-1)+\frac{(n-1)(n-2)}{2}\] 
    symmetric pairs of orange boxes 
    contributing to the unipotent dimension. 
    However, there are more than $m+1$ pairs of free orbits. 

    In addition to the $m+1$ pairs of orbits passing through $a_{2i,1}$ for $i=1,2,...,m+1$, 
    we also have the orbits passing through $a_{2m+2i+1,1}$ for $i=1,...,n-m-1$. 
    These orbits are likewise free, 
    and together with the $m+1$ pairs in the case of $m\leqslant n-1$, 
    they account for all free orbits. 
    Counting both the orange boxes and the free orbits, the total number of degrees of freedom is:
    \[\frac{(m+1)m}{2}+(m+1)(n-1)+\frac{(n-1)(n-2)}{2}+n=\frac{(m+n+1)(m+n)}{2}-m.\]

        \begin{center}
        \begin{tikzpicture}
            \draw [step=1cm,gray,very thin] (0,0) grid (16,16);
            \filldraw[fill=gray, fill opacity=0.3, draw=black] (0,11) rectangle (5,16);
            \filldraw[fill=gray, fill opacity=0.3, draw=black] (6,1) rectangle (15,10);
            \filldraw[fill=gray, fill opacity=0.3, draw=black] (0,1) rectangle (5,10);
            \filldraw[fill=gray, fill opacity=0.3, draw=black] (6,11) rectangle (15,16);
            \node at (0.5,15.5){0};

            \node at (1.5,14.5){0};
            \node at (1.5,12.5){0};
            \node at (1.5,10.5){0};
            \node at (1.5,9.5){0};
            \node at (1.5,7.5){0};
            \node at (1.5,5.5){0};
            \node at (1.5,3.5){0};
            \node at (1.5,1.5){0};
            \node at (1.5,0.5){0};

            \node at (3.5,14.5){0};
            \node at (3.5,12.5){0};
            \node at (3.5,10.5){0};
            \node at (3.5,9.5){0};
            \node at (3.5,7.5){0};
            \node at (3.5,5.5){0};
            \node at (3.5,3.5){0};
            \node at (3.5,1.5){0};
            \node at (3.5,0.5){0};

            \node at (5.5,14.5){0};
            \node at (5.5,12.5){0};
            \node at (5.5,10.5){0};
            \node at (5.5,9.5){0};
            \node at (5.5,7.5){0};
            \node at (5.5,5.5){0};
            \node at (5.5,3.5){0};
            \node at (5.5,1.5){0};

            \node at (6.5,14.5){0};
            \node at (6.5,12.5){0};
            \node at (6.5,10.5){0};
            \node at (6.5,9.5){0};
            \node at (6.5,7.5){0};
            \node at (6.5,5.5){0};
            \node at (6.5,3.5){0};
            \node at (6.5,1.5){0};
            \node at (6.5,0.5){0};

            \node at (8.5,14.5){0};
            \node at (8.5,12.5){0};
            \node at (8.5,10.5){0};
            \node at (8.5,9.5){0};
            \node at (8.5,7.5){0};
            \node at (8.5,5.5){0};
            \node at (8.5,3.5){0};
            \node at (8.5,1.5){0};
            \node at (8.5,0.5){0};

            \node at (10.5,14.5){0};
            \node at (10.5,12.5){0};
            \node at (10.5,10.5){0};
            \node at (10.5,9.5){0};
            \node at (10.5,7.5){0};
            \node at (10.5,5.5){0};
            \node at (10.5,3.5){0};
            \node at (10.5,1.5){0};
            \node at (10.5,0.5){0};

            \node at (12.5,14.5){0};
            \node at (12.5,12.5){0};
            \node at (12.5,10.5){0};
            \node at (12.5,9.5){0};
            \node at (12.5,7.5){0};
            \node at (12.5,5.5){0};
            \node at (12.5,3.5){0};
            \node at (12.5,1.5){0};
            \node at (12.5,0.5){0};

            \node at (14.5,14.5){0};
            \node at (14.5,12.5){0};
            \node at (14.5,10.5){0};
            \node at (14.5,9.5){0};
            \node at (14.5,7.5){0};
            \node at (14.5,5.5){0};
            \node at (14.5,3.5){0};
            \node at (14.5,1.5){0};
            \node at (14.5,0.5){0};

            \node at (15.5,14.5){0};
            \node at (15.5,12.5){0};
            \node at (15.5,9.5){0};
            \node at (15.5,7.5){0};
            \node at (15.5,5.5){0};
            \node at (15.5,3.5){0};
            \node at (15.5,1.5){0};
            \node at (15.5,0.5){0};

            \node at (2.5,13.5){0};
            \node at (4.5,11.5){0};
            \node at (7.5,8.5){0};
            \node at (9.5,6.5){0};
            \node at (11.5,4.5){0};
            \node at (13.5,2.5){0};

            \filldraw[fill=orange, fill opacity=0.3, draw=black] (0,13) rectangle (1,14);
            \filldraw[fill=orange, fill opacity=0.3, draw=black] (0,11) rectangle (1,12);

            \filldraw[fill=orange, fill opacity=0.3, draw=black] (2,11) rectangle (3,12);

            \filldraw[fill=orange, fill opacity=0.3, draw=black] (2,15) rectangle (3,16);
            \filldraw[fill=orange, fill opacity=0.3, draw=black] (4,15) rectangle (5,16);

            \filldraw[fill=orange, fill opacity=0.3, draw=black] (4,13) rectangle (5,14);

            \filldraw[fill=orange, fill opacity=0.3, draw=black] (0,2) rectangle (1,3);
            \filldraw[fill=orange, fill opacity=0.3, draw=black] (2,2) rectangle (3,3);
            \filldraw[fill=orange, fill opacity=0.3, draw=black] (4,2) rectangle (5,3);

            \filldraw[fill=orange, fill opacity=0.3, draw=black] (0,4) rectangle (1,5);
            \filldraw[fill=orange, fill opacity=0.3, draw=black] (2,4) rectangle (3,5);
            \filldraw[fill=orange, fill opacity=0.3, draw=black] (4,4) rectangle (5,5);

            \filldraw[fill=orange, fill opacity=0.3, draw=black] (0,6) rectangle (1,7);
            \filldraw[fill=orange, fill opacity=0.3, draw=black] (2,6) rectangle (3,7);
            \filldraw[fill=orange, fill opacity=0.3, draw=black] (4,6) rectangle (5,7);

            \filldraw[fill=orange, fill opacity=0.3, draw=black] (0,8) rectangle (1,9);
            \filldraw[fill=orange, fill opacity=0.3, draw=black] (2,8) rectangle (3,9);
            \filldraw[fill=orange, fill opacity=0.3, draw=black] (4,8) rectangle (5,9);

            \filldraw[fill=orange, fill opacity=0.3, draw=black] (7,11) rectangle (8,12);
            \filldraw[fill=orange, fill opacity=0.3, draw=black] (7,13) rectangle (8,14);
            \filldraw[fill=orange, fill opacity=0.3, draw=black] (7,15) rectangle (8,16);

            \filldraw[fill=orange, fill opacity=0.3, draw=black] (9,11) rectangle (10,12);
            \filldraw[fill=orange, fill opacity=0.3, draw=black] (9,13) rectangle (10,14);
            \filldraw[fill=orange, fill opacity=0.3, draw=black] (9,15) rectangle (10,16);

            \filldraw[fill=orange, fill opacity=0.3, draw=black] (11,11) rectangle (12,12);
            \filldraw[fill=orange, fill opacity=0.3, draw=black] (11,13) rectangle (12,14);
            \filldraw[fill=orange, fill opacity=0.3, draw=black] (11,15) rectangle (12,16);

            \filldraw[fill=orange, fill opacity=0.3, draw=black] (13,11) rectangle (14,12);
            \filldraw[fill=orange, fill opacity=0.3, draw=black] (13,13) rectangle (14,14);
            \filldraw[fill=orange, fill opacity=0.3, draw=black] (13,15) rectangle (14,16);

            \filldraw[fill=orange, fill opacity=0.3, draw=black] (11,2) rectangle (12,3);

            \filldraw[fill=orange, fill opacity=0.3, draw=black] (9,2) rectangle (10,3);
            \filldraw[fill=orange, fill opacity=0.3, draw=black] (9,4) rectangle (10,5);
            
            \filldraw[fill=orange, fill opacity=0.3, draw=black] (7,2) rectangle (8,3);
            \filldraw[fill=orange, fill opacity=0.3, draw=black] (7,4) rectangle (8,5);
            \filldraw[fill=orange, fill opacity=0.3, draw=black] (7,6) rectangle (8,7);

            \filldraw[fill=orange, fill opacity=0.3, draw=black] (13,4) rectangle (14,5);
            
            \filldraw[fill=orange, fill opacity=0.3, draw=black] (13,6) rectangle (14,7);
            \filldraw[fill=orange, fill opacity=0.3, draw=black] (11,6) rectangle (12,7);
            
            \filldraw[fill=orange, fill opacity=0.3, draw=black] (13,8) rectangle (14,9);
            \filldraw[fill=orange, fill opacity=0.3, draw=black] (11,8) rectangle (12,9);
            \filldraw[fill=orange, fill opacity=0.3, draw=black] (9,8) rectangle (10,9);

            \draw[->] (-0.3,15.5) to[bend left] (0.5,14.7);
            \draw[->] (-0.3,13.5) to[bend left] (0.5,12.7);
            \draw[->] (-0.3,11.5) to[bend left] (0.5,10.7);
            \draw[->] (-0.3,10.5) to[bend left] (0.5,9.7);
            \draw[->] (-0.3,8.5) to[bend left] (0.5,7.7);
            \draw[->] (-0.3,6.5) to[bend left] (0.5,5.7);
            \draw[->] (-0.3,4.5) to[bend left] (0.5,3.7);
            \draw[->] (-0.3,2.5) to[bend left] (0.5,1.7);

            \draw[->] (0.7,16.5) to[bend left] (1.5,15.7);
            \draw[->] (0.7,14.5) to[bend left] (1.5,13.7);
            \draw[->] (0.7,12.5) to[bend left] (1.5,11.7);
            \draw[->] (0.7,9.5) to[bend left] (1.5,8.7);
            \draw[->] (0.7,7.5) to[bend left] (1.5,6.7);
            \draw[->] (0.7,5.5) to[bend left] (1.5,4.7);
            \draw[->] (0.7,3.5) to[bend left] (1.5,2.7);
            \draw[->] (0.7,1.5) to[bend left] (1.5,0.7);
            \draw[->] (0.7,0.5) to[bend left] (1.5,-0.3);

            \draw[->] (1.7,15.5) to[bend left] (2.5,14.7);
            \draw[->] (1.7,13.5) to[bend left] (2.5,12.7);
            \draw[->] (1.7,11.5) to[bend left] (2.5,10.7);
            \draw[->] (1.7,10.5) to[bend left] (2.5,9.7);
            \draw[->] (1.7,8.5) to[bend left] (2.5,7.7);
            \draw[->] (1.7,6.5) to[bend left] (2.5,5.7);
            \draw[->] (1.7,4.5) to[bend left] (2.5,3.7);
            \draw[->] (1.7,2.5) to[bend left] (2.5,1.7);

            \draw[->] (2.7,16.5) to[bend left] (3.5,15.7);
            \draw[->] (2.7,14.5) to[bend left] (3.5,13.7);
            \draw[->] (2.7,12.5) to[bend left] (3.5,11.7);
            \draw[->] (2.7,9.5) to[bend left] (3.5,8.7);
            \draw[->] (2.7,7.5) to[bend left] (3.5,6.7);
            \draw[->] (2.7,5.5) to[bend left] (3.5,4.7);
            \draw[->] (2.7,3.5) to[bend left] (3.5,2.7);
            \draw[->] (2.7,1.5) to[bend left] (3.5,0.7);
            \draw[->] (2.7,0.5) to[bend left] (3.5,-0.3);

            \draw[->] (3.7,15.5) to[bend left] (4.5,14.7);
            \draw[->] (3.7,13.5) to[bend left] (4.5,12.7);
            \draw[->] (3.7,11.5) to[bend left] (4.5,10.7);
            \draw[->] (3.7,10.5) to[bend left] (4.5,9.7);
            \draw[->] (3.7,8.5) to[bend left] (4.5,7.7);
            \draw[->] (3.7,6.5) to[bend left] (4.5,5.7);
            \draw[->] (3.7,4.5) to[bend left] (4.5,3.7);
            \draw[->] (3.7,2.5) to[bend left] (4.5,1.7);

            \draw[->] (4.7,16.5) to[bend left] (5.5,15.7);
            \draw[->] (4.7,14.5) to[bend left] (5.5,13.7);
            \draw[->] (4.7,12.5) to[bend left] (5.5,11.7);
            \draw[->] (4.7,9.5) to[bend left] (5.5,8.7);
            \draw[->] (4.7,7.5) to[bend left] (5.5,6.7);
            \draw[->] (4.7,5.5) to[bend left] (5.5,4.7);
            \draw[->] (4.7,3.5) to[bend left] (5.5,2.7);
            \draw[->] (4.7,1.5) to[bend left] (5.5,0.7);
            \draw[->] (4.7,0.5) to[bend left] (5.5,-0.3);

            \draw[->] (5.7,16.5) to[bend left] (6.5,15.7);
            \draw[->] (5.7,14.5) to[bend left] (6.5,13.7);
            \draw[->] (5.7,12.5) to[bend left] (6.5,11.7);
            \draw[->] (5.7,9.5) to[bend left] (6.5,8.7);
            \draw[->] (5.7,7.5) to[bend left] (6.5,6.7);
            \draw[->] (5.7,5.5) to[bend left] (6.5,4.7);
            \draw[->] (5.7,3.5) to[bend left] (6.5,2.7);
            \draw[->] (5.7,1.5) to[bend left] (6.5,0.7);
            \draw[->] (5.7,0.5) to[bend left] (6.5,-0.3);

            \draw[->] (6.7,15.5) to[bend left] (7.5,14.7);
            \draw[->] (6.7,13.5) to[bend left] (7.5,12.7);
            \draw[->] (6.7,11.5) to[bend left] (7.5,10.7);
            \draw[->] (6.7,10.5) to[bend left] (7.5,9.7);
            \draw[->] (6.7,8.5) to[bend left] (7.5,7.7);
            \draw[->] (6.7,6.5) to[bend left] (7.5,5.7);
            \draw[->] (6.7,4.5) to[bend left] (7.5,3.7);
            \draw[->] (6.7,2.5) to[bend left] (7.5,1.7);

            \draw[->] (7.7,16.5) to[bend left] (8.5,15.7);
            \draw[->] (7.7,14.5) to[bend left] (8.5,13.7);
            \draw[->] (7.7,12.5) to[bend left] (8.5,11.7);
            \draw[->] (7.7,9.5) to[bend left] (8.5,8.7);
            \draw[->] (7.7,7.5) to[bend left] (8.5,6.7);
            \draw[->] (7.7,5.5) to[bend left] (8.5,4.7);
            \draw[->] (7.7,3.5) to[bend left] (8.5,2.7);
            \draw[->] (7.7,1.5) to[bend left] (8.5,0.7);
            \draw[->] (7.7,0.5) to[bend left] (8.5,-0.3);

            \draw[->] (8.7,15.5) to[bend left] (9.5,14.7);
            \draw[->] (8.7,13.5) to[bend left] (9.5,12.7);
            \draw[->] (8.7,11.5) to[bend left] (9.5,10.7);
            \draw[->] (8.7,10.5) to[bend left] (9.5,9.7);
            \draw[->] (8.7,8.5) to[bend left] (9.5,7.7);
            \draw[->] (8.7,6.5) to[bend left] (9.5,5.7);
            \draw[->] (8.7,4.5) to[bend left] (9.5,3.7);
            \draw[->] (8.7,2.5) to[bend left] (9.5,1.7);

            \draw[->] (9.7,16.5) to[bend left] (10.5,15.7);
            \draw[->] (9.7,14.5) to[bend left] (10.5,13.7);
            \draw[->] (9.7,12.5) to[bend left] (10.5,11.7);
            \draw[->] (9.7,9.5) to[bend left] (10.5,8.7);
            \draw[->] (9.7,7.5) to[bend left] (10.5,6.7);
            \draw[->] (9.7,5.5) to[bend left] (10.5,4.7);
            \draw[->] (9.7,3.5) to[bend left] (10.5,2.7);
            \draw[->] (9.7,1.5) to[bend left] (10.5,0.7);
            \draw[->] (9.7,0.5) to[bend left] (10.5,-0.3);

            \draw[->] (10.7,15.5) to[bend left] (11.5,14.7);
            \draw[->] (10.7,13.5) to[bend left] (11.5,12.7);
            \draw[->] (10.7,11.5) to[bend left] (11.5,10.7);
            \draw[->] (10.7,10.5) to[bend left] (11.5,9.7);
            \draw[->] (10.7,8.5) to[bend left] (11.5,7.7);
            \draw[->] (10.7,6.5) to[bend left] (11.5,5.7);
            \draw[->] (10.7,4.5) to[bend left] (11.5,3.7);
            \draw[->] (10.7,2.5) to[bend left] (11.5,1.7);

            \draw[->] (11.7,16.5) to[bend left] (12.5,15.7);
            \draw[->] (11.7,14.5) to[bend left] (12.5,13.7);
            \draw[->] (11.7,12.5) to[bend left] (12.5,11.7);
            \draw[->] (11.7,9.5) to[bend left] (12.5,8.7);
            \draw[->] (11.7,7.5) to[bend left] (12.5,6.7);
            \draw[->] (11.7,5.5) to[bend left] (12.5,4.7);
            \draw[->] (11.7,3.5) to[bend left] (12.5,2.7);
            \draw[->] (11.7,1.5) to[bend left] (12.5,0.7);
            \draw[->] (11.7,0.5) to[bend left] (12.5,-0.3);

            \draw[->] (12.7,15.5) to[bend left] (13.5,14.7);
            \draw[->] (12.7,13.5) to[bend left] (13.5,12.7);
            \draw[->] (12.7,11.5) to[bend left] (13.5,10.7);
            \draw[->] (12.7,10.5) to[bend left] (13.5,9.7);
            \draw[->] (12.7,8.5) to[bend left] (13.5,7.7);
            \draw[->] (12.7,6.5) to[bend left] (13.5,5.7);
            \draw[->] (12.7,4.5) to[bend left] (13.5,3.7);
            \draw[->] (12.7,2.5) to[bend left] (13.5,1.7);

            \draw[->] (13.7,16.5) to[bend left] (14.5,15.7);
            \draw[->] (13.7,14.5) to[bend left] (14.5,13.7);
            \draw[->] (13.7,12.5) to[bend left] (14.5,11.7);
            \draw[->] (13.7,9.5) to[bend left] (14.5,8.7);
            \draw[->] (13.7,7.5) to[bend left] (14.5,6.7);
            \draw[->] (13.7,5.5) to[bend left] (14.5,4.7);
            \draw[->] (13.7,3.5) to[bend left] (14.5,2.7);
            \draw[->] (13.7,1.5) to[bend left] (14.5,0.7);
            \draw[->] (13.7,0.5) to[bend left] (14.5,-0.3);

            \draw[->] (14.7,15.5) to[bend left] (15.5,14.7);
            \draw[->] (14.7,13.5) to[bend left] (15.5,12.7);
            \draw[->] (14.7,11.5) to[bend left] (15.5,10.7);
            \draw[->] (14.7,10.5) to[bend left] (15.5,9.7);
            \draw[->] (14.7,8.5) to[bend left] (15.5,7.7);
            \draw[->] (14.7,6.5) to[bend left] (15.5,5.7);
            \draw[->] (14.7,4.5) to[bend left] (15.5,3.7);
            \draw[->] (14.7,2.5) to[bend left] (15.5,1.7);

            \draw[->] (15.7,15.5) to[bend left] (16.5,14.7);
            \draw[->] (15.7,13.5) to[bend left] (16.5,12.7);
            \draw[->] (15.7,11.5) to[bend left] (16.5,10.7);
            \draw[->] (15.7,10.5) to[bend left] (16.5,9.7);
            \draw[->] (15.7,8.5) to[bend left] (16.5,7.7);
            \draw[->] (15.7,6.5) to[bend left] (16.5,5.7);
            \draw[->] (15.7,4.5) to[bend left] (16.5,3.7);
            \draw[->] (15.7,2.5) to[bend left] (16.5,1.7);

            \node at (0.5,14.5){$a_{2,1}$};
        \end{tikzpicture}
    \end{center}

\end{proof}

\begin{corollary}\label{HomGnsym}
    The dimensions of \[\Hom(G_{n-2,1},G_{n-2,1}^\vee)^{\sym}\;\text{and}\; \Hom(G_{0,n-1},G_{0,n-1}^\vee)^{\sym}\] 
    are both equal to $\frac{n(n-1)}{2}$. 
\end{corollary}

\begin{proposition}\label{HomFVGmn}
The group $\Hom([FV],G_{m,n})$ is a direct product of $m+n+1$ copies of $\G_a$. 
\end{proposition}

\begin{proof}
    It suffices to compute the solutions to the equation $F-V=0$ in the covariant Dieudonné module of $G_{m,n}$ . 
    The covariant Dieudonné module of $G_{m,n}$ has the form: 
    \[e_1\xrightarrow{F}e_2\xleftarrow{V}e_3\xrightarrow{F}...\xleftarrow{V}e_{2m+1}\xrightarrow{F}e_{2m+2}\xrightarrow{F}e_{2m+3}\xleftarrow{V}...\xleftarrow{V}e_{2m+2n+2}\xleftarrow{V}e_1\]  
    A general element is written as \(\sum_{i=1}^{2m+2n+2}a_ie_i\), then we have \[\begin{aligned}
        F(\sum_{i=1}^{2m+2n+2}a_ie_i)=(\sum_{i=1}^{m}a_{2i-1}^\sigma e_{2i})+a_{2m+1}^\sigma e_{2m+2}+(\sum_{i=1}^n a_{2m+2i}^\sigma e_{2m+2i+1}),\\
        V(\sum_{i=1}^{2m+2n+2}a_ie_i)=(\sum_{i=1}^{m}a_{2i+1}^{\sigma^{-1}} e_{2i})+(\sum_{i=1}^n a_{2m+2i+2}^{\sigma^{-1}} e_{2m+2i+1})+a_{1}^{\sigma^{-1}} e_{2m+2n+2}.
    \end{aligned}\]
    The condition $F=V$ implies 
    \[0=a_1=a_3^{\sigma^{-2}}=...=a_{2m+1}^{\sigma^{-2m}}=0,\quad a_{2m+2}=a_{2m+4}^{\sigma^{-2}}=...=a_{2m+2n+2}^{\sigma^{-2n}},\]
    so the last chain of equalities contributes exactly one copy of $\G_a$.  
    The remaining free coefficients give $m+n$ additional copies of $\G_a$, so the total solution space is isomorphic to $\G_a^{\,m+n+1}$.  
\end{proof}

\begin{corollary}\label{HomFVGn}
The groups $\Hom([FV],G_{n-2,1})$ and $\Hom([FV],G_{0,n-1})$ are both direct products of $n$ copies of $\G_a$. 
\end{corollary}

\begin{proposition}
    Assume $A$ is a supersingular abelian variety, and $a(A)=g-1$. 
    Assume that $A[p]$ is isomorphic to 
    \[[FV]^{g-1-m-n}\oplus G_{m,n},\]
    then $\dim U[p]=\frac{g(g-1)}{2}-\min\{m,n-1\}. $
\end{proposition}

\begin{proof}
    This follows directly from Proposition \ref{HomGmnsym} and \ref{HomFVGmn}. 
    Assume $A[p]\cong  [FV]^{g-1-m-n}\oplus G_{m,n}$,
    so that $A[p]^\vee\cong [FV]^{g-1-m-n}\oplus G_{m,n}^\vee$. 
    Then  
    \begin{equation}
        \begin{aligned}
        &\Hom(A[p],A^\vee[p])^{\sym} \cong \Hom(G_{m,n},G_{m,n}^\vee)^{\sym} \\
        &\oplus (\Hom([FV],G_{m,n}^\vee)\times \Hom(G_{m,n},[FV]))^{\sym,(g-1-m-n)}\\
        &\oplus \Hom([FV]^{g-1-m-n},[FV]^{g-1-m-n})^{\sym}.
        \end{aligned}
    \end{equation}
The symmetric part $(\Hom([FV],G_{m,n}^\vee)\times \Hom(G_{m,n},[FV]))^{\sym,(g-1-m-n)}$ is isomorphic to 
$\Hom([FV],G_{m,n}^\vee)$, 
since duality gives an isomorphism between the two factors. 
Furthermore, the symmetric part of \[\Hom([FV]^{g-1-m-n},[FV]^{g-1-m-n})\]
has dimension $\frac{(g-1-m-n)(g-2-m-n)}{2}$. 
By Lemmas \ref{HomGnsym} and \ref{HomFVGn}, the total dimension sums up to
\[\begin{aligned}
    &\frac{(m+n)(m+n+1)}{2}-\min\{m,n-1\}\\
    +&(g-1-m-n)(m+n+1)+\frac{(g-1-m-n)(g-2-m-n)}{2}\\
    =&\frac{g(g-1)}{2}-\min\{m,n-1\}.
\end{aligned} \]
\end{proof}

\begin{corollary}\label{almostsuperspecial}
    Assume $A$ is a supersingular abelian variety, and $a(A)=g-1$. 
    Assume for some $2\leqslant n\leqslant g$ 
    such that \[A[p]\cong [FV]^{g-n}\oplus G_{n-2,1} \quad \text{or} \quad [FV]^{g-n}\oplus G_{0,n-1},\]
    Then \[\Hom(A[p],A^\vee[p])^{\sym}\] 
    has dimension $\frac{g(g-1)}{2}$. 
\end{corollary}

\begin{corollary}
    Let $A$ be a supersingular abelian threefold over $k$ with $a$-number $2$. 
    Then $U_A\cong \G_a^3$.  
\end{corollary}

\begin{proof}
    This is a direct consequence of Corollary \ref{Abe3g=2} and Corollary\ref{almostsuperspecial}.  
\end{proof}

\begin{proposition}\label{NPPabe3table}
    Let $A$ be an abelian $3$-fold. The following table lists all possible isogeny classes of $U$.
    \begin{small}
    \begin{center}{\rm
    \begin{tabular}{ c|c|c|c|c|c } \label{Table2}
    Newton polygon &    $a$      & $A[p]$           & $\dim(U[p])$& $\dim(U)$ & isogeny class of $U$\\ 
    supersingular, &    $3$      &  $[FV][FV][FV]$  &  $3$       & $3$      & $\G_a\times \G_a\times \G_a$  \\ 
    supersingular, &    $2$      & $[FV][FFVV]$     &  $3$       & $3$     & $\G_a\times \G_a\times \G_a$\\    
                   &             & $[FVFFVV]$       &  $3$       & $3$     & $\G_a\times \G_a\times \G_a$\\    
                   &             & $[VVFFVF]$       &  $3$       & $3$     & $\G_a\times \G_a\times \G_a$\\    
    supersingular, &    $1$      & $[FFFVVV]$       &  $2$       & $3$     & $\G_a\times \W_2$\\            
    $1/3$ type,    &    $2$      &  $[FFV][VVF]$    &  $2$       & $2$     & $\G_a\times \G_a$  \\ 
    $1/3$ type,    &    $1$      & $[FFFVVV]$       &  $2$       & $2$     & $\G_a\times \G_a$\\                 
    almost supersingular,& $2$   & $[F][V][FV][FV]$ &  $1$       & $1$     & $\G_a$  \\ 
    almost supersingular & $1$   & $[F][V][FFVV]$   &  $1$       & $1$     & $\G_a$\\           
     almost ordinary   & $1$   & $[F][F][V][V][FV]$ &  $0$       & $0$     & $0$  \\ 
     ordinary      &    $0$   & $[F][F][F][V][V][V]$&  $0$       & $0$     & $0$\\           
    \end{tabular}}
\end{center}
    \end{small} 
\end{proposition}

\begin{proof}
    The supersingular cases follow directly from Corollary \ref{supergeneraldimUp} and Corollary \ref{almostsuperspecial}. 
    For the $1/3$ type cases, we determine $A[p]$ according to its $a$-number. 
    If $a=2$, $A[p]$ has to be $[FFV][VVF]$ (as other possibilities give higher dimension of $U[p]$). 
    If $a=1$, the only possible case is $[FFFVVV]$. 
    The remaining cases (almost ordinary and ordinary) are immediate from the classification of 
    \( A[p] \)-types and the corresponding unipotent dimensions. 
\end{proof}

\chapter{Explicit Hom Group: Calculations} \label{Chap5}

We have classified the isogeny class of $U$ for abelian $3$-folds in Proposition \ref{NPPabe3table}. 
The only non-trivial cases arise for supergeneral abelian $3$-folds, where $U$ is isogenous to $\G_a\times \W_2$. 
Determining the precise structure of $U$ in these cases remains an open question. 
In this chapter, we present an explicit computation of the symmetric Hom group 
$\Hom(A[p^2],A^\vee[p^2])^{\sym}$ for a supergeneral
principally polarized abelian $3$-fold $A$. 
The result shows that the connected component $\Hom^{\circ}(A[p^2],A^\vee[p^2])^{\sym}$ 
distribute in a well-behaved manner 
across the moduli space of principally polarized abelian varieties. 

Recall from \cite{LiOo} that
the moduli space $\SSS_{3,1}$ of principally polarized supersingular abelian $3$-folds 
is described by a surjective map 
\[\bigsqcup_{\eta} \PP_{3,\eta}' \rightarrow \SSS_{3,1}\]
where $\PP_{3,\eta}'$ is the moduli of rigid polarized flag type quotients ($\mathtt{PFTQ}$) with respect to 
$(E^3, \eta)$, where $E^3$ is a triple product of a fixed supersingular elliptic curve, 
and the disjoint union ranges over polarizations $\eta$ with $\Ker(\eta)=\Ker(F^3)$. 
All components $\PP_{3,\eta}'$ are smooth of dimension $2$, and isomorphic to each other. 
Moreover, $\PP_{3,\eta}'$ is a $\mathbb{P}^1$-bundle over the Fermat curve of degree $p+1$, 
which we denote by $\mathcal{V}$. 
There exists a canonical section $t:\mathcal{V}\rightarrow \PP_{3,\eta}'$. 
\[\begin{aligned}
    &\PP_{3,\eta}'\xlongrightarrow{\mathbb{P}^1} \mathcal{V}=\mathcal{Z}(X^{p+1}+Y^{p+1}+Z^{p+1})\subset \mathbb{P}^2,\\
     &\PP_{3,\eta}'\supset \mathcal{T}\xlongleftarrow{t} \mathcal{V},
\end{aligned}\]
where $\mathcal{T}=t(\mathcal{V})$.

The finite automorphism group $G_\eta=\mathrm{Aut}(E^3, \eta)/\{\pm 1\}$ acts on $\PP_{3,\eta}'$, 
and the quotient $\PP_{3,\eta}'/G_\eta$ is birationally equivalent to an irreducible component of $\mathcal{S}_{3,1}$. 
The curve $\mathcal{T}$ is contracted to a superspecial point in $\mathcal{S}_{3,1}$. 
Away from this curve, the map $\PP_{3,\eta}'/G_\eta\rightarrow \mathcal{S}_{3,1}$ is finite-to-one. 

Points on $\mathcal{T}$ correspond exactly to superspecial points. Away from this section, 
the points whose image in $\mathcal{V}$ lie in $\mathcal{V}(\F_{p^2})$ correspond to 
abelian threefolds with $a=2$, 
all other points correspond to supergeneral abelian $3$-folds (\cite[Example 9.4]{LiOo}). 

Now, let $b,c\in k$. Define the following invariants: 
\[D(i,j)=\det\left(\begin{matrix}
    b^{\sigma^i}&c^{\sigma^i}\\
    b^{\sigma^{j}}&c^{\sigma^{j}}
\end{matrix}\right)=b^{\sigma^i}c^{\sigma^{j}}-c^{\sigma^i}b^{\sigma^{j}}.\]

For example, 
\[D(1,-1)=\det\left(\begin{matrix}
    b^{\sigma}&c^\sigma\\
    b^{\sigma^{-1}}&c^{\sigma^{-1}}
\end{matrix}\right)=b^{\sigma}c^{\sigma^{-1}}-c^{\sigma}b^{\sigma^{-1}}.\]

For each pair of $(b,c)$, suppose $D(1,-1)\neq 0$, we define a unipotent group $\W_2^{(b,c)}$ as a subgroup of $\G_a\times \W_2$, 
by choosing elements $(x,z)\in \G_a\times \W_2$ 
such that $x$ and $\overline{z}$ (the reduction of $z$ mod $p$) satisfy: 

\[\begin{aligned}
    &\frac{D(3,-1)-D(3,1)-D(1,-1)}{D(1,-1)}x+\frac{D(3,1)}{D(1,-1)}\overline{z}+\frac{D(3,-1)}{D(1,-1)}\overline{z}^p\\
    =&\frac{D(4,2)+D(4,0)}{D(4,2)}x^{p^3}+\frac{D(4,0)}{D(4,2)}\overline{z}^{p^3}+\frac{D(2,0)}{D(4,2)} \overline{z}^{p^4}.
\end{aligned}
\]

Then we have the following result. 

\begin{theorem}\label{distribution}
    Let $x\in \PP_{3,\eta}'\backslash \mathcal{T}(k)$ be a point, 
    whose image in $\mathcal{V}$ is the point $(1,b,c)\in \mathcal{V}(k)\backslash \mathcal{V}(\F_{p^2})$. 
    Let $A$ be the principally polarized abelian $3$-fold associated to the point $x$. 
    We then have an abstract isomorphism
    \[\Hom^{\circ}(A[p^2],A^\vee[p^2])^{\sym}\cong \W_2^{(b,c)}(k)\oplus \G_a(k).\]
\end{theorem}

Here, $\Hom^{\circ}(A[p^n],A^\vee[p^n])^{\sym}$ denotes 
the connected component of the total Hom group 
$\Hom(A[p^n],A^\vee[p^n])^{\sym}$. 
This shows that, on each fiber of the projection $\PP_{3,\eta}'\rightarrow \mathcal{V}$, 
the connected commutative unipotent group $\Hom^{\circ}(A[p^n],A^\vee[p^n])^{\sym}$ stays invariant
away from the superspecial point. 

\section{Explicit $\Hom$ Groups}

The contravariant Dieudonn\'e module of a supergeneral principally polarized abelian $3$-fold can be discribed as follows. 
First, choose a point $(1,b,c)$ in \[\mathcal{V}=\mathcal{Z}(X^{p+1}+Y^{p+1}+Z^{p+1})\subset P^2\]
where without loss of generality, we assume the first coordinate is nonzero and normalize it to be $1$. 
Assume further that $(1,b,c)\notin \mathcal{V}(\F_{p^2})$. 

\begin{lemma}
    Under the condition that $(1,b,c)\notin \mathcal{V}(\F_{p^2})$, 
    we have \[D(1,-1)\neq 0.\]
\end{lemma}

\noindent\emph{Proof:} First, note that $b,c\neq 0$. For instance, if $b=0$, then the defining equation 
implies $c^{p+1}+1=0$, which leads to $c^{p^2}=(\frac{-1}{c})^p=(-1)^{p+1}c=c$, and hence $c\in \F_{p^2}$, 
contradicting our assumption. 
Next, suppose $D=0$. Then we would have $((\frac{b}{c})^{\frac{1}{p}})^{p^2-1}=1$, 
which implies $k=\frac{b}{c}\in \F_{p^2}\backslash\{0\}$, and hence $k^{p+1}\in \F_p^\times$.
The defining equation becomes $1+(kc)^{p+1}+c^{p+1}=0$, implying $c^{p+1}=\frac{-1}{1+k^{p+1}}\in \F_p$, 
which again implies $c\in \F_{p^2}$, contradicting our assumption. $\square$

Let $M_0$ be the free $W$-module generated by $\{e_1,e_2,e_3,f_1,f_2,f_3\}$ 
with the Frobenius and Verschiebung actions defined by $Fe_i=Ve_i=f_i$, $Ff_i=Vf_i=pe_i$. 
Given any $l\in k$, lift $(b,c,l)$ arbitrarily to $W$, we define a Dieudonn\'e submodule $M_{b,c,l}\subset M_0$ 
to be the free $W$-module generated by the basis $\{\epsilon_1,\epsilon_2,\epsilon_3,\gamma_1,\gamma_2,\gamma_3\}$, 
where: 
\[\begin{aligned}
    &\epsilon_1=e_1+be_2+ce_3+lf_1\\
    &\epsilon_2=pe_2\\
    &\epsilon_3=pe_3\\
    &\gamma_1=pf_1\\
    &\gamma_2=f_1+b^{\sigma}f_2+c^{\sigma}f_3\\
    &\gamma_3=f_1+b^{\sigma^{-1}}f_2+c^{\sigma^{-1}}f_3.
\end{aligned}\]

By the theory Li and Oort developed in \cite{LiOo}, the module $M_{b,c,l}$ corresponds to the (contravariant) 
Dieudonn\'e module of principally polarized supergeneral abelian threefolds associated to the 
point $(1,l)\in \mathbb{P}^1(k)$ in the $\mathbb{P}^1$-fiber of the point 
$(1,b,c)\in \mathcal{V}(k)\backslash \mathcal{V}(\F_p^2)$. 
The module $M_{b,c,l}$ doesn't depend on the choice of the lift, 
but the basis does depend on the choice of the lift. 
We abuse the notation here, we also use $(b,c,l)$ to denote their chosen lifts to $W$ 
(or $W_2$, depending on the context). 
Let $A_{b,c,l}$ denote this abelian threefold.  

\begin{theorem}\label{explicitHom}
    For each $(b,c)\in \mathcal{V}(k)\backslash \mathcal{V}(\F_{p^2})$, and any $l\in k$, 
    there is an abstract isomorphism of abelian groups: 
    \[ \Hom^{\circ}(A_{b,c,l}[p^2],A_{b,c,l}^{\vee}[p^2])^{\sym}\cong \W_2^{(b,c)}(k)\oplus\G_a(k).\]
\end{theorem}

\begin{lemma}
    Under the basis $\{\epsilon_1,\epsilon_2,\epsilon_3,\gamma_1,\gamma_2,\gamma_3\}$, 
    the matrix of $[\mathcal{F}_M]$ and $[\mathcal{V}_M]$ are: 
    \[[\mathcal{F}_M]=\left(\begin{matrix}
        pl&0&0&p^2&p&p\\
        -lb&0&0&-pb&b^{\sigma^2}-b& 0\\
        -lc&0&0&-pc&c^{\sigma^2}-c&0\\
        -l^2&\frac{c^{\sigma}-c^{\sigma^{-1}}}{D}&\frac{b^{\sigma}-b^{\sigma^{-1}}}{-D}&-pl&-l&-l\\
        1& \frac{pc^{\sigma^{-1}}}{D}&\frac{pb^{\sigma^{-1}}}{-D}&0&0&0\\
        0&\frac{-pc^{\sigma}}{D}& \frac{-pb^{\sigma}}{-D}&0&0&0
    \end{matrix}\right),\]
    \[[\mathcal{V}_M]=\left(\begin{matrix}
        pl&0&0&p^2&p&p\\
        -lb&0&0&-pb&0& b^{\sigma^{-2}}-b\\
        -lc&0&0&-pc&0&c^{\sigma^{-2}}-c\\
        -l^2&\frac{c^{\sigma}-c^{\sigma^{-1}}}{D}&\frac{b^{\sigma}-b^{\sigma^{-1}}}{-D}&-pl&-l&-l\\
        0& \frac{pc^{\sigma^{-1}}}{D}&\frac{pb^{\sigma^{-1}}}{-D}&0&0&0\\
        1&\frac{-pc^{\sigma}}{D}& \frac{-pb^{\sigma}}{-D}&0&0&0
    \end{matrix}\right),\] 
where $D=D(1,-1)=b^{\sigma}c^{\sigma^{-1}}-c^{\sigma}b^{\sigma^{-1}}$. 
\end{lemma}

\begin{proof} 
    This is a straightforward calculation which we omit. 
\end{proof}

Since $A_{b,c,l}$ is principally polarzied, its covariant and contravariant Dieudonné module are isomorphic. 
This implies that we must solve the following equations in $\W_2$: 
\[[f][\mathcal{F}_M]=[\mathcal{V}_M]^{t\sigma}[f]^\sigma,\qquad [f]=-[f]^t.\]

Suppose that $(x,z)\in \W_2^{(b,c)}$, and $a_{1},a_{2}$ are elements in $\G_a(k)$ that satisfy 
\[(b-b^{\sigma^2})(a_{1}-a_{1}^{p^2})+(c-c^{\sigma^2})(a_{2}-a_{2}^{p^2})=(l^\sigma-l)x+(l^\sigma-l^{\sigma^2})x^{\sigma^2}.\]
Let $g$ and $h$ be the following elements determined by $(x,z,a_{1},a_{2})$:
\[g=lba_{1}+lca_{2}+l^2x+[lba_{1}+lca_{2}]^\sigma+(l^2x^\sigma)^\sigma,\]
\[h=lx^{\sigma^{-1}}+[(b^{\sigma^{-2}}-b)a_{1}+(c^{\sigma^{-2}}-c)a_{2}]^\sigma-(lx)^\sigma.\]
We claim that the following matrix $[f]$ 
is a solution of the equations $[f][\mathcal{F}_M]=[\mathcal{V}_M]^{t\sigma}[f]^\sigma$, and $[f]=-[f]^t$:
\begin{small}
\[[f]=\left(\begin{matrix}
    0&pa_{1}&pa_{2}&px&-z^\sigma+pg&z\\
    -pa_{1}&0&0&0&p[-(sx)^\sigma+(\frac{cz^{\sigma^2}+c^{\sigma^2}z^{\sigma}}{D^\sigma})]&p[(-sx)^{\sigma^{-1}}+(\frac{c^{\sigma^{-2}}z+cz^{\sigma^{-1}}}{D^{\sigma^{-1}}})]\\
    -pa_{2}&0&0&0&p[-(tx)^\sigma-(\frac{bz^{\sigma^2}+b^{\sigma^2}z^\sigma}{D^\sigma})]&p[(-tx)^{\sigma^{-1}}-(\frac{b^{\sigma^{-2}}z+bz^{\sigma^{-1}}}{D^{\sigma^{-1}}})]\\
    -px&0&0&0&0&0\\
    z^{\sigma}-pg&&&0&0&ph\\
    -z&&&0&-ph&0
\end{matrix}\right).\]
\end{small} 
Here, we have \[ s=\frac{c^{\sigma}-c^{\sigma^{-1}}}{D},\quad t=\frac{b^{\sigma}-b^{\sigma^{-1}}}{-D}.\]
There are some entries left blank, but these terms are uniquely determined by anti-symmetricity. 
Again, the verification is an direct calculation, which we omit. 

By counting dimension, this tells us that $\Hom^{\circ}(A_{b,c,l}[p^2],A_{b,c,l}^{\vee}[p^2])^{\sym}$ 
is isomorphic to the subgroup \[(a_{1},a_{2},x,z)\in (\G_a\times\G_a\times \G_a\times \W_2)(k)\] 
such that $(x,z)\in \W_2^{(b,c)}$ and 
\[(b-b^{\sigma^2})(a_{1}-a_{1}^{p^2})+(c-c^{\sigma^2})(a_{2}-a_{2}^{p^2})=(l^\sigma-l)x+(l^\sigma-l^{\sigma^2})x^{\sigma^2}.\]

\noindent{\bf Proof of Theorem \ref{explicitHom}:}
We claim that the above algebraic group is abstractly isomorphic to the group $\G_a\times \W_2^{(b,c)}$. 
To prove this, we invoke several key lemmas. 

\begin{lemma}\label{sectionlemma}
    Let $H\subset G$ be connected commutative unipotent algebraic groups.
Then the exact sequence
    \[0\rightarrow H \rightarrow G \rightarrow G/H\rightarrow 0\]
    admits an algebraic section (not necessarily a homomorphism) $G/H\rightarrow G$. 
\end{lemma}

\begin{proof}
    This is \cite[Proposition 2.3.4]{Ros}.
\end{proof}

This lemma implies that there is a factor system $f: G/H\times G/H \rightarrow H$ such that 
\[f(a,b)+f(a+b,c)-f(b,c)-f(a,b+c),\]
so that we can equip the space $H\times (G/H)$ with a group structure given by
$(h_1,a)+(h_2,b)=(h_1+h_2+f(a,b),a+b)$
such that $G$ is isomorphic to this product. 

We say that two factor systems $f_1,f_2$ are equivalent to each other, 
if there is a map $g: G/H\rightarrow H$ such that $f_1(a,b)-f_2(a,b)=g(a+b)-g(a)-g(b)$. 

\begin{lemma}\label{Serrelemma}
    Any commutative factor system $f:\G_a\times \G_a\rightarrow \G_a$ 
    is equivalent to a factor system of the form 
    \[\sum_{i=0}^n k_if_0^{p^i}, \quad k_i\in k\] where \[f_0(a,b)=\frac{1}{p}((a+b)^p-a^p-b^p).\] 
\end{lemma}
\begin{proof}
    This is Lemme 3 of \cite{Lar}. 
\end{proof}

\begin{lemma}\label{connected}
    Let $\alpha\in k\backslash{\F_{p^2}}$. Then the curve 
    \[(x^{p^2}-x)-\alpha(y^{p^2}-y)=0\]
    in $\G_a\times \G_a$ is connected. 
\end{lemma}
\begin{proof}
    The $k$-algebra 
$k[x,y]/((x^{p^2}-x)-\alpha(y^{p^2}-y))$ is an integral domain by Artin-Schreier theory, 
\end{proof}

Now consider the 2-dimensional subgroup $G$ of $(a_{1},a_{2},x)\in \G_a\times \G_a\times \G_a$ cut out by the equation: 
\[(b-b^{\sigma^2})(a_{1}-a_{1}^{p^2})+(c-c^{\sigma^2})(a_{2}-a_{2}^{p^2})=(l^\sigma-l)x+(l^\sigma-l^{\sigma^2})x^{\sigma^2}.\]

$G$ admits a projection to $\G_a$ that sends $(a_{1},a_{2},x)$ to $x$. 
The kernel of this projection is the curve $C$, defined by 
\[(b-b^{\sigma^2})(a_{1}-a_{1}^{p^2})+(c-c^{\sigma^2})(a_{2}-a_{2}^{p^2})=0.\]
The ratio $\frac{b-b^{\sigma^2}}{c-c^{\sigma^2}}$ does not lie in $\F_{p^2}$, or otherwise $A_{b,c,l}$'s $a$-number 
would be strictly larger than $1$.  
By Lemma \ref{connected}, $C$ is connected, so $C$ is abstractly isomorphic to $\G_a$. 
By Lemma \ref{sectionlemma}, $G_2$ is isomorphic to a group structure on $C\times \G_a$, 
where the second factor is the projection to $x$. 
According to Lemma \ref{Serrelemma}, we know every factor system on $\G_a\times \G_a\rightarrow \G_a$. 
Those factor systems that give an algebraic group killed by $p$ are all equivalent to the trivial factor system. 
Consequently, there is an isomorphism $G_2\cong C\times G_a$, where the second component is the projection to $x$-coordinate. 
As a result, $\Hom^{\circ}(A[p^2],A^{\vee}[p^2])^{\sym}$ is isomorphic to the subgroup of $G_2\times \W_2\cong C\times \G_a\times \W_2$ 
such that the last two components $(x,z)\in \W_2^{(b,c)}$. 
So $\Hom^{\circ}(A_{b,c,l}[p^n],A_{b,c,l}^{\vee}[p^n])^{\sym}\cong (\G_a\times \W_2^{(b,c)})(k)$, as desired. \qed

\section{Conjecture and Example}

We propose the following conjecture concerning the isogeny class of $U$ for an arbitrary supergeneral abelian variety of dimension $g$, 
as a natural generalization of our result in the case $g=3$. 
However, at present, we are unable to provide a justification for this conjecture, nor have we identified any counterexamples.

\begin{conj}\label{supergeneralconj}
    Let $A/k$ be a supergeneral abelian variety of dimension $g$. 
    Then $U$ is isogenous to \[\prod_{i=1}^{g-1} \W_{i}.\]
\end{conj}

We have a simplest example, where we can verify this conjecture.  
Consider the Dieudonné module $M:=W_{\sigma}[F,V]/(F^g-V^g)$, 
it is isoclinic of slope $1/2$ and it is generated by $1$ single element. 
It is the (covariant) Dieudonné module of a supersingular abelian variety $A$ 
whose $a$-number equals $1$. 

\begin{proposition}
    The unipotent group $U$ of $A$ is isogenous to 
    \[ \prod_{n=1}^{g-1}\W_n.\]
\end{proposition}
\begin{proof} The Dieudonné module $M$ is a free $W$-module of rank $2g$ with basis 
$$e_i=F^{i-1}, \ \  i=1,\ldots, g, \quad e_j=V^{2g+1-j},\ \  j=g+1,\ldots,2g.$$
We have $Fe_i=e_{i+1}$ for $i=1,\ldots, g$, and $Fe_i=pe_{i+1}$ for $i=g+1,\ldots, 2g$
(with the convention that subscripts are considered as integers modulo $2g$.) 

The Dieudonn\'e module $M$ has a principal polarisation,
that is, a perfect, alternating bilinear form
$$\psi\colon M\otimes_W M\to W,$$ 
such that $\psi(F(x),y)=\psi(x,V(y))^\sigma$,
defined as follows. Let $\theta\in W$
be an invertible element such that $\sigma^g(\theta)=-\theta$. Then
$$\psi(e_i,e_j)=0, \ \ i-j\neq g, \quad \psi(e_i,e_{i+g})=\sigma^{i-1}(\theta), \ \ i=1,\ldots, g,$$
is a principal polarisation. Write $B$ for the matrix of 
$\psi$ in the basis $e_1,\ldots,e_{2g}$, so that $B_{ij}=\psi(e_i,e_j)$ for $i,j=1,\ldots,2g$.

Since $M$ is generated by $e_1$ as a $W_{\sigma}[F,V]$-module, there is at most one endomorphism
$f\in\mathrm{End}_{W_{\sigma}[F,V]}(M/p^n)$ such that 
$$f(e_1)=\sum_{i=1}^{2g} a_i e_i,$$
for given $a_i\in W_n=W/p^n$ for $i=1,\ldots,2g$, and there is exactly one such $f$
if and only if $F^g f(e_1)=V^g f(e_1)$. We have
$$F^g f(e_1)=\sum_{i=1}^g a_i^{\sigma^g}p^{i-1}e_{i+g} +
\sum_{i=g+1}^{2g}a_i^{\sigma^g}p^{2g+1-i}e_{i-g},$$
$$V^g f(e_1)=\sum_{i=1}^g a_i^{\sigma^{-g}}p^{i-1}e_{i+g} +
\sum_{i=g+1}^{2g}a_i^{\sigma^{-g}}p^{2g+1-i}e_{i-g}.$$
Thus $F^g f(e_1)=V^g f(e_1)$ is equivalent to 
the condition that the image of $a_i$ in $W_{n+1-i}$ is invariant under $\sigma^{2g}$,
for $i=1,\ldots,2g$. Thus the connected component of $0$ of the group $k$-scheme
$\mathrm{End}_{W_{\sigma}[F,V]}(M/p^n)$ is isomorphic to 
\begin{equation}
\prod_{i=1}^g {\mathbb W}_{i-1}\times \prod_{i=g+1}^{2g}{\mathbb W}_{2g+1-i},
\label{W}\end{equation}
where the $i$-th factor is a group subscheme of $\W_n$ with coordinate $a_i$.

Let $A$ be the $(2g\times 2g)$-matrix of $f$ in the basis $\{e_1,\ldots, e_{2g}\}$,
that is, $A_{ij}$ is the $e_i$-th coordinate of $f(e_j)$.
To $f$ we associate the map $\tilde f\in\Hom(M/p^n,(M/p^n)^\vee)$ that sends
$x$ to the linear form on $M/p^n$ whose value on $y$ is $\psi(Ax,y)$.
The dual map $\tilde f^\vee\in\Hom(M/p^n,(M/p^n)^\vee)$ sends $x$ to the linear form
whose value on $y$ is $\psi(Ay,x)=-\psi(x,Ay)$. Thus $\tilde f^\vee=-f^\vee$ if and only if
$A^t=BAB^{-1}$, or equivalently, $A=B^{-1}A^tB$. Since $f$ is an endomorphism of $W_{\sigma}[F,V]$-modules, this condition is equivalent to 
$A(e_1)=B^{-1}A^tB(e_1)$. We have $B(e_1)=-\theta e_{g+1}$. Next, 
we use that $F^i(e_1)=e_{1+i}$ for $i=1,\ldots,g$ and $V^j(e_1)=e_{2g+1-j}$ 
for $j=1,\ldots, g$ to calculate
$$A^t(e_{g+1})=a_{g+1}e_1+a_g^\sigma e_2+\ldots+a_1^{\sigma^g}e_{g+1}+
a_{2g}^{\sigma^{-(g-1)}} e_{g+2}+\ldots+a_{g+2}^{\sigma^{-1}}e_{2g}.$$
From this we obtain 
$$B^{-1}A^tB(e_1)=a_1^{\sigma^g}e_1+
\sum_{i=2}^g \frac{\theta}{\theta^{\sigma^{i-1}}}a_{2g+2-i}^{\sigma^{i-g-1}}e_i
-\sum_{i=g+1}^{2g}\frac{\theta}{\theta^{\sigma^{i-g-1}}}a_{2g+2-i}^{\sigma^{i-g-1}}e_i.$$
Thus the involution on $\mathrm{End}_{W_\sigma[F,V]}(M/p^n)$ that sends $A$ to $B^{-1}A^tB$ (which is indeed
an involution because $B^t=-B^{-1}$)
sends the $i$-factor of (\ref{W}) isomorphically onto the $2g-i$-th factor, and acts on
the $g+1$-th factor $\W_g$ as $-1$. We conclude that $\mathrm{End}_{W_{\sigma}[F,V]}(M/p^n)^+$ is isomorphic to 
$\prod_{i=1}^g {\mathbb W}_{i-1}$.    
\end{proof}

Although we are currently unable to prove Conjecture \ref{supergeneralconj} in full, 
there is evidence suggesting that the $p$-exponent of $U$ can always be bounded above by $g-1$. 
Indeed, for non-supersingular abelian varieties, the exponent of the unipotent part of the full Hom group 
$\Hom(A[p^n],A^\vee[p^n])$ is known to be bounded by $g-1$ (see \cite{LNV}). 
In the case of supersingular abelian varieties, the key to establishing this bound 
lies in Proposition \ref{almostsuperspecial} and a Lemma due to Livia Grammatica: if an isogeny $f:A\rightarrow B$ has kernel isomorphic to a product of $\alpha_p$'s, 
then the $p$-exponent of $U_A$ and $U_B$ differ at most by $1$. 
A proof of this lemma, which makes uses of Nygaard filtration, will appear in a forthcoming joint 
paper \cite{Skor}.   

\part{Further Topics}

\chapter{Partial results and applications} \label{Chap6}

\section{Formal Brauer Group}\label{formalBrauer}
\medskip

In this section, we briefly review the theory of formal groups, Cartier modules, 
and the formal Brauer group $\widehat{\Br(X)}$. 

By a formal group, we mean a group object in the category of formal schemes over $k$. 
Every formal group over $k$ represents a sheaf on the big fppf site over $k$. 
Since $k$ is perfect, every formal group $G$ of finite presentation over $k$ can be written as a product $G^\circ\times \pi_0(G)$, 
where $G^{\circ}$ is connected and $\pi_0(G)$ is \'etale \cite[II.7, p. 34]{De}. 
Every connected formal group $G$ of finite presentation
fits into a unique exact sequence of connected groups: 
\[0\rightarrow G_{\red} \rightarrow G \rightarrow G_{\inf} \rightarrow 0,\]
where $G_{\red}$ is a smooth formal group (so that 
the underlying formal scheme is $Spf(k[[x_1,...,x_n]])$) 
and $G_{\inf}=G/G_{\red}$ is an infinitesimal group scheme (see \cite[Section II.10, p. 43]{De}). 

To study formal groups, one associates to each formal group $G$ the \emph{Cartier module}, 
or \emph{the module of typical curves}, $\TC(G):=\Hom(\widehat{\W},G)$ of $G$. 
Here, $\widehat{\W}$ is the formal Witt group over $k$ that associates to a $k$-algebra $R$ 
the group of sequences $(a_0,a_1,a_2,...)$ of nilpotent elements in $R$ such that $a_i=0$ for all  sufficiently large $i$. 
$\TC(G)$ is a $W$-module with semilinear operators $F$ and $V$, induced by the action of $V$ and $F$ on $\widehat{\W}$. 
They satisfy well-known semilinear formula $Fa=a^{\sigma}F$ and $Va=a^{\sigma^{-1}}V$, for any $a\in W$. 
Thus the Cartier module is a (left) module over the Dieudonn\'e ring $W_{\sigma}[F,V]$.  

It is worth mentioning that the Cartier module of a $p$-divisible formal group is isomorphic 
to the dual of the contravariant Dieudonn\'e module. 
Explicitly, for any finite group scheme $G$, 
we let $M(G)$ be the contravariant Dieudonn\'e module of $G$ (see \cite[Chapter III]{De}).
Let $\DM(G):=M(G^{\vee})$ be the covariant Dieudonn\'e module functor, 
where $G^\vee$ is the Cartier dual of $G$. 
For any connected $p$-divisible formal group $G:=\varinjlim G_n$
we define $\DM(G):=\varprojlim \DM(G_n)$.
Then there is an isomorphism $\TC(G)\cong \DM(G)$ compatible with the action of $F$ and $V$. 
This is stated in \cite[Section I.3, p.~102]{AM} and \cite[Remark 3.3]{Br}, 
see \cite{Mo} for an explicit proof. 

 
The Cartier module can be viewed as a covariant functor from the category of finite-dimensional connected
commutative smooth formal groups to the category of  finitely generated $W_{\sigma}[F,V]$-modules. 
The essential image of $\TC$ is the category of finitely generated $W_{\sigma}[F,V]$-modules $M$ 
such that $V$ is injective 
on $M$ and $M/VM$ is a finite dimensional $k$-vector space. 
The dimension of the formal group equals $\dim_k M/VM$. 

Unlike the contravariant theory in \cite{Man}, the covariant theory $\TC$ is insensitive to the infinitesimal part $G_{\inf}$. 
In particular, if $G$ is a connected finite-dimensional formal group, 
then $\TC(G)=\TC(G_{\red})$. The information about $G_{\inf}$ is reflected in the first derived functor of $\TC$: 
we have $R^1\TC(G)=\DM(G_{\inf})$, the covariant Dieudonn\'e module of $G_{\inf}$. 
We have $R^i\TC(G)=0$ for $i>1$ (\cite[Remark 11.10]{BO} and \cite[Theorem 4.2(c)]{Oo}). 

Artin and Mazur \cite{AM} studied the representability of 
various deformation functors for a proper smooth map $X \rightarrow S$. 
For a sheaf $E$ of abelian groups on the big \'etale site of $X$
they defined a functor $\widetilde{E}_S$ that assigns
to a pair $(X',S')$, where $S'$ is a $S$-scheme and $X'\rightarrow X\times_S S'$ is \'etale, 
the abelian group $ker(E(X')\rightarrow E(X'\times_{S'} S'_{\red}))$. 
This gives a sheaf on the relative \'etale site $(X/S)_{\et}$ (see, for example, \cite[Section 2.3]{BL}). 
Using the projection $\pi: (X/S)_{\et}\rightarrow S_{\Et}$, 
one defines $\Phi^i(X/S,E):=R^i\pi_*\widetilde{E}_S$. 
They showed that if $R^{i-1}f_*E$ is formally smooth, 
then $\Phi^i(X/S,E)$ is pro-representable by a formal group (\cite[Proposition 1.8]{AM}). 
In particular, when $E=\G_m$, if the Picard scheme $\Pic_{X/S}$ is smooth, then
the functor $\Phi^2(X/S,\G_m)$ is pro-representable by a formal group over $S$, called
the formal Brauer group $\widehat{\Br(X/S)}$.

For any sheaf $\mathcal{F}$ of abelian groups on the fppf site over $k$, 
we can extend the definition of TC to $\mathcal{F}$ by denoting 
$\TC(\mathcal{F}):=\Hom(\widehat{\W},\mathcal{F})$. 
If $\mathcal{F}$ is representable by a group scheme $G$, then this notion agrees 
with the Cartier module of the formal completion of $G$ along the zero section. 

The Artin-Mazur formal group has a flat variant, let us denote by $\Phi_{\fl}^i(X,\G_m)$, 
which is defined as the fppf-sheafification of the presheaf 
\[A\rightarrow ker(\HH^i(A,\G_m)\rightarrow \HH^i(A_{\red},\G_m))\]
on the category of Artinian $k$-algebra $(\art_{k})^{\op}$ eqquiped with fppf-topology. 
Our notation is consistent with the notation in \cite{BO}; in \cite{Ray}, 
it is denoted by $\widetilde{R^if_*}(\widehat{\G_m})$; 
in \cite{Ek2}, it is simply denoted by $R^if_*(\widehat{\G_m})$. 
It is proved in \cite[10.6]{BO} that $\Phi_{\fl}^i(X,\G_m)$ 
can be equivalently defined as $\widehat{R^if_*(\G_m)}$, the formal completion of $R^if_*(\G_m)$, the sheaf on $(\art_{k})^{\op}$ defined by 
\[A\rightarrow \Ker(R^if_*(\G_m)(A)\rightarrow R^if_*(\G_m)(A_{red})).\]
It is then a well-known theorem that $\Phi_{\fl}^i(X,\G_m)$ are all pro-representable. 

\begin{theorem}[Raynaud, Bragg]
    Let $X$ be a proper smooth variety over $k$. 
    The flat Artin-Mazur formal groups $\Phi_{\fl}^i(X,\G_m)$ are all pro-representable. 
\end{theorem}

\begin{proof}
    See \cite[Theorem 10.8]{BO}; \cite[Theorem 4.1.2]{Ray}.
\end{proof} 

We finish this section by recalling a relation between the flat Artin-Mazur formal groups 
and the Hodge--Witt cohomology groups $\HH^i(X,WO_X)$. 

\begin{theorem}\label{Eke81}
There is an exact sequence 
    \begin{equation*}
        0 \rightarrow R^1\TC(\Phi^{i-1}_{\fl}(X,\G_m)) \rightarrow \HH^i(X,WO_X) \rightarrow \TC(\Phi^{i}_{\fl}(X,\G_m)) \rightarrow 0,
    \end{equation*}
    where $R^1\TC(\Phi^{i-1}_{\fl}(X,\G_m))$ is isomorphic to the $V^{\infty}$-torsion subgroup of $\HH^i(X,WO_X)$, 
    and $\TC(\Phi^{i-1}_{\fl}(X,\G_m))$ is isomorphic to the quotient of $\HH^i(X,WO_X)$ by its $V^{\infty}$-torsion subgroup. 
\end{theorem}

\begin{proof}
    See \cite[Proposition III.8.1]{Ek2}.
\end{proof}

\section{Cartier Module of the Higher Direct Image} 

\begin{proposition}
    Let $n$ be an integer such that $p^n$ annihilates
$\HH^2(X,WO_X)[p^\infty]$. 

{\rm (1)} we have the following exact sequences: 
\begin{equation*}
    \begin{aligned}
    &0 \rightarrow \TC(\widehat{R^2f_*(\mu_{p^n})}) \rightarrow \HH^2(X,WO_{X})[p^{\infty}/V^\infty]\\
    & \rightarrow \DM(\widehat{\Pic}_{X/k}/p^n) \rightarrow \DM(\widehat{R^2f_*\mu_{p^n}}_{\inf}) \rightarrow \DM((\Phi^2_{\fl}(X,\G_m)[p^n])_{\inf}) \rightarrow 0.
    \end{aligned}
\end{equation*}

{\rm (2)} when $\Pic_{X/k}$ is smooth, we have the following two isomorphism:
\[  \TC(\widehat{R^2f_*(\mu_{p^n})}) \cong \HH^2(X,WO_{X})[p^{\infty}],\qquad (\widehat{R^2f_*\mu_{p^n}})_{\inf}\cong (\Phi^2_{\fl}(X,\G_m)[p^n])_{\inf}\]
where $\TC$ is the Cartier module functor (see Section \ref{2}).
\end{proposition}

\noindent{\bf Proof}: 
We have the following exact sequenses of fppf sheaves on $\art_k^{\op}$: 
\begin{equation}\label{Kummer}
    0\rightarrow \widehat{\Pic}_{X/k}/p^n \rightarrow \widehat{R^2f_*(\mu_{p^n})} \rightarrow \Phi^2_{\fl}(X,\G_m)[p^n] \rightarrow 0,
\end{equation}
\begin{equation}\label{pn}
    0\rightarrow \Phi^2_{\fl}(X,\G_m)[p^n] \rightarrow \Phi^2_{\fl}(X,\G_m) \xrightarrow{p^n} \Phi^2_{\fl}(X,\G_m).
\end{equation}
Here, the first sequence is obatined by applying the completion functor 
$F\rightarrow \widehat{F}$ (which is an exact functor, see \cite[Lemma 10.4]{BO}) to the short sequence $0\rightarrow \Pic_{X/k}/p^n\rightarrow  R^2f_*(\mu_{p^n})\rightarrow R^2f_*(\G_m)[p^n]\rightarrow 0$. 
Theorem \ref{Eke81} gives an isomorphism
$R^1\TC(\Phi^1_{\fl}(X,\G_m))\cong\HH^2(X,WO_{X})[V^\infty]$, and this 
is the covariant Dieudonné module of $\widehat{\Pic}_{X/k,\inf}$.
Writing $\HH^2(X,WO_{X})/V^\infty$ for the quotient
$\HH^2(X,WO_{X})/\HH^2(X,WO_{X})[V^{\infty}]$, we have an isomorphism
$\TC(\Phi^2_{\fl}(X,\G_m))\cong\HH^2(X,WO_{X})/V^\infty$.
Thus, applying the covariant functor $\TC$ to (\ref{pn}) we obtain an exact sequence
\[
    0\rightarrow \TC(\Phi^2_{\fl}(X,\G_m)[p^n]) \rightarrow \HH^2(X,WO_{X})/[V^{\infty}] \xrightarrow{p^n} \HH^2(X,WO_{X})/[V^{\infty}].
\]
For a sufficiently large $n$ such that $p^n$ kills the $p^{\infty}$-torsion part of $\HH^2(X,WO_{X})$, 
the above sequence gives an isomorphism
\[ 
    \TC(\Phi^2_{\fl}(X,\G_m)[p^n])\cong \HH^2(X,WO_{X})[p^{\infty}/V^n].
\]
Applying $\TC$ to (\ref{Kummer}) and noticing that $\TC(\widehat{\Pic}_{X/k}/p^n)=0$, we obtain an exact sequence
\begin{equation}
    \begin{aligned}
    &0 \rightarrow \TC(\widehat{R^2f_*(\mu_{p^n})}) \rightarrow \HH^2(X,WO_{X})[p^{\infty}/V^\infty]\\
    & \rightarrow \DM(\widehat{\Pic}_{X/k}/p^n) \rightarrow \DM((\widehat{R^2f_*\mu_{p^n}})_{\inf}) \rightarrow \DM(\Phi^2_{\fl}(X,\G_m)[p^n]_{\inf}) \rightarrow 0.
    \end{aligned}
\end{equation}
When $\Pic_{X/k}$ is smooth, $\widehat{\Pic}_{X/k}/p^n$ is trivial, so $\DM(\widehat{\Pic}_{X/k}/p^n)=0$ 
and we have the isomorphisms: 
\[\begin{aligned}
    \TC(\Phi^2_{\fl}(X,\G_m)) \cong \HH^2(X,WO_{X})[p^{\infty}/V^\infty],\\
\DM(\widehat{R^2f_*\mu_{p^n}}_{\inf}) \cong \DM((\Phi^2_{\fl}(X,\G_m)[p^n])_{\inf})
\end{aligned} \] as desired. 
This finishes the proof. $\square$

\medskip

\section{Comparing Flat with Crystalline Cohomology} \label{5}
In this section, we study the relation between the fppf cohomology and the crystalline cohomology. 
Consider the natural map 
\begin{equation}\label{dlog}
    d\log\colon \HH^i_{\fppf}(X,\mu_{p^n})\rightarrow \HH^i_{\cris}(X/W_n)
\end{equation}
induced by the logarithmic differential $d\log:\G_m/\G_m^{p^n}[-1]\rightarrow W_n\Omega_{X/k}^{\bullet}$, 
where $W_n\Omega_{X/k}^{\bullet}$ is the de Rham--Witt complex of level $n$. 
Ogus proved in \cite{Og} that $d\log\colon\HH^2_{\fppf}(X,\mu_p)\to \HH^2_{\dR}(X/k)$ is injective
when the Fr\"olicher spectral sequence degenerates. 
The following theorem can be seen as a generalization of this result.

\begin{theorem} \label{Inj}
    Let $X$ be a smooth proper variety over $k$.
    Suppose that the Fr\"olicher spectral sequence 
    degenerates at the $E_1$ page and that $\HH^1_{\rm{cris}}(X/W)$ is torsion-free. 
    Then the map $d\log$ is injective in degree $2$ for all $n$: 
    \begin{equation*}
        \HH^2_{\rm{fppf}}(X,\mu_{p^n})\hookrightarrow \HH^2_{\rm{cris}}(X/W_n).
    \end{equation*}
\end{theorem}

This theorem applies to arbitrary abelian varieties and K3 surfaces. Indeed,
by \cite[Proposition 5.1]{Od} and \cite[Proposition 2.2]{DI}, the Fr\"olicher spectral sequence degenerates for abelian varieties and K3 surfaces. It is well-known that
in these cases all crystalline cohomology groups are torsion-free, see, e.g., \cite[Exercise 1.8]{Li}.


One consequence of Theorem \ref{Inj} is that multiplication by $m$ on an abelian variety $A$ acts on 
$\HH^i_{\rm{fppf}}(A,\mu_{p^n})$ as $m^i$ when $i\leqslant 2$. 

The aim of this section is to prove Theorem \ref{Inj}.

For each $m$ and $n$, we have the short exact sequence of fppf sheaves:
\[ 0 \rightarrow \mu_{p^m} \rightarrow \mu_{p^{n+m}} \rightarrow \mu_{p^n} \rightarrow 0\]
Applying the derived functor $R^i\epsilon_*$ we have 
\[0\rightarrow R^1\epsilon_*(\mu_{p^m}) \rightarrow R^1\epsilon_*(\mu_{p^{n+m}}) \rightarrow R^1\epsilon_*(\mu_{p^n}) \rightarrow 0\]
from
which we get \[0\rightarrow \G_m/\G_m^{p^m} \xrightarrow[]{p^n} \G_m/\G_m^{p^{n+m}} \rightarrow \G_m/\G_m^{p^n} \rightarrow 0.\] 

We would like to construct a similar sequence for $W_n\Omega_{X/k}^\bullet$. 

Recall that Illusie defined the canonical filtration (see Definition \ref{CanFildef}) on the de Rham-Witt complex, 
and we showed that the multiplication by $p$ map:  $p: W_{n+1}\Omega_X^\bullet\rightarrow W_{n+1}\Omega_X^\bullet$ 
has kernel $\Fil^n W_{n+1}\Omega_X^{\bullet}$ (see Proposition \ref{kerP}). Moreover, for each $i$ such that $0\leqslant i\leqslant n+1$, we have 
\[\Ker(p^i: W_{n+1}\Omega^\bullet\rightarrow W_{n+1}\Omega^\bullet)=\Fil^{n+1-i}W_{n+1}\Omega_X^\bullet.\]
We conclude that multiplication by $p^n$ map on $W_{n+m}\Omega_X^\bullet$ induces an exact sequence of complexes
\[0\rightarrow W_m\Omega_X^\bullet \rightarrow W_{m+n}\Omega_X^\bullet \rightarrow W_{m+n}\Omega_X^\bullet /p^n \rightarrow 0.\]

By Lemma \ref{quasiisop}, for each $n\geqslant 0$, multiplication by $p$ induces a quasi-isomorphism
 $\gr^nW\Omega_X^\bullet \rightarrow \gr^{n+1}W\Omega_X^\bullet$. The following statement generalizes 
\cite[Corollaire I.3.16]{Il1} that says that the canonical projection 
$W\Omega^{\bullet}_{X/k}/p\rightarrow \Omega^{\bullet}_{X/k}$ is a quasi-isomorphism . 

\begin{proposition}\label{quasiiso}
    The natural projection map induces a quasi-isomorphism $W_{m+n}\Omega_X^\bullet/p^n \rightarrow  W_n\Omega_X^\bullet$. 
\end{proposition}
\begin{proof}
    Our proof is similar to the proof of \cite[Corollaire I.3.15]{Il1}. 
We use induction on $m$. For $m=0$, this is trivial. 
Suppose that for $m-1$ the claim is true. Consider the following commutative diagram
with exact rows:
\begin{center}
    \begin{tikzcd}
        0 \arrow{r} & \Fil^{n+m-1}W_{n+m}\Omega_X^\bullet \arrow{r} \arrow[d,"p^n=0"] & W_{n+m}\Omega_X^\bullet \arrow{r} \arrow[d,"p^n"] & W_{n+m-1}\Omega_X^\bullet  \arrow{r} \arrow[d,"p^n"] & 0\\
        0 \arrow{r} & \Fil^{n+m-1}W_{n+m}\Omega_X^\bullet \arrow{r} & W_{n+m}\Omega_X^\bullet \arrow{r} & W_{n+m-1}\Omega_X^\bullet \arrow{r} & 0
    \end{tikzcd}
\end{center}
Applying the snake lemma we obtain an exact sequence 
\[\begin{aligned}
    0 \rightarrow& \Fil^{n+m-1}W_{n+m}\Omega_X^\bullet \rightarrow \Fil^{m}W_{n+m}\Omega_X^\bullet \rightarrow \Fil^{m-1}W_{m+n-1}\Omega_X^\bullet \rightarrow \\
    &\Fil^{n+m-1}W_{n+m}\Omega_X^\bullet \rightarrow W_{m+n}\Omega_X^\bullet/p^n \rightarrow W_{m+n-1}\Omega_X^\bullet/p^n \rightarrow 0
\end{aligned}\]
We have obvious exact sequences
\[0 \rightarrow \Fil^{n+m-1}W_{n+m}\Omega_X^\bullet \rightarrow \Fil^{m}W_{n+m}\Omega_X^\bullet \rightarrow \Fil^{m}W_{m+n-1}\Omega_X^\bullet \rightarrow 0,\]
\[0 \rightarrow \Fil^{m}W_{m+n-1}\Omega_X^\bullet \rightarrow \Fil^{m-1}W_{m+n-1}\Omega_X^\bullet \rightarrow \Fil^{m-1}W_m\Omega^\bullet_X \rightarrow 0,\]
and by repeatly using Lemma \ref{quasiisop} we obtain a quasi-isomorphism
\[\Fil^{m-1}W_m\Omega^\bullet_X \rightarrow \Fil^{n+m-1}W_{n+m}\Omega_X^\bullet.\]
Thus $W_{m+n}\Omega_X^\bullet/p^n \rightarrow W_{m+n-1}\Omega_X^\bullet/p^n $ 
is a quasi-isomorphism.
\end{proof}

\medskip

So we have the following commutative diagram of complexes:
\begin{center}
    \begin{tikzcd}
        0 \arrow{r}  & \G_m/\G_m^{p^m}[-1] \arrow{r} \arrow[d,"d\log"] & \G_m/\G_m^{p^{n+m}}[-1] \arrow{r} \arrow[d,"d\log"] & \G_m/\G_m^{p^{n}}[-1] \arrow{r} \arrow{d} & 0\\
        0 \arrow{r} & W_m\Omega_X^\bullet \arrow[r,"p^n"]                   & W_{n+m}\Omega_X^\bullet \arrow{r}                 &  W_{n}\Omega_X^\bullet \arrow{r}          & 0
    \end{tikzcd}
\end{center}
This proves the following statement.

\begin{proposition}\label{LES}
    For each $m$ and $n$, we have a functorial commutative diagram, where the horizontal sequences are exact: 
    \begin{small}
    \begin{center}\label{LESD}
        \begin{tikzcd}
            ...\arrow{r} & \HH^{i-1}_{\fppf}(X,\mu_{p^n}) \arrow{r} \arrow[d,"d\log"]  & \HH^i_{\fppf}(X,\mu_{p^m}) \arrow{r} \arrow[d,"d\log"] & \HH^i_{\fppf}(X,\mu_{p^{n+m}}) \arrow{r} \arrow[d,"d\log"] & \HH^i_{\fppf}(X,\mu_{p^n}) \arrow{r} \arrow[d,"d\log"] &  ...\\
            ...\arrow{r} & \HH^{i-1}_{\cris}(X/W_n) \arrow{r} & \HH^i_{\cris}(X/W_m) \arrow{r} & \HH^i_{\cris}(X/W_{n+m}) \arrow{r} & \HH^i_{\cris}(X/W_n) \arrow{r} &...
        \end{tikzcd}
    \end{center}
    \end{small}

\end{proposition}

\section{A Criterion of Injectivity}

A natural question is when the map 
$d\log\colon \HH^i_{\fppf}(X,\mu_{p^n})\to \HH^i_{\cris}(X/W_n)$ is injective.
We will prove that it is injective when $i=2$ for abelian varieties and K3 surfaces.
However, even for K3 surfaces or abelian surfaces, this is not injective when $i=3$ and
$X$ is supersingular.

Our core argument is the following

\begin{proposition} \label{Injarg}
    Let $Z\Omega^1_X:=ker(\Omega_X^1\xrightarrow[]{d} \Omega^2_X)$ be the sheaf of closed $1$-forms. Assume that 
    
{\rm  (1)} $\HH^{i-2}(X,Z\Omega_X^1)\rightarrow \HH^{i-2}(X,\Omega_X^1)$ is surjective;

{\rm  (2)} the natural map $\HH^{i-1}(X,Z\Omega^1_X) \rightarrow \HH^i_{\dR}(X/k)$ 
    induced by $Z\Omega^1_X[-1] \rightarrow \Omega_X^\bullet$ is injective.

\noindent Then the map $d\log: \HH^i_{\fppf}(X,\mu_{p}) \rightarrow \HH^i_{\dR}(X/k)$ is injective. 
    If, moreover, $\HH^i_{\cris}(X/W)$ is torsion-free, 
    then $d\log: \HH^i_{\fppf}(X,\mu_{p^n}) \rightarrow \HH^i_{\cris}(X/W_n)$ is injective for all $n$. 
\end{proposition}
\begin{proof}
     
When $n=1$, we have the following exact sequence of etale sheaves (cf.~\cite[Proposition 0.2.1.18]{Il1}):
\[0\rightarrow \G_m/\G_m^{p} \xrightarrow[]{d\log} Z\Omega_X^1 \xrightarrow[]{W^*-C} \Omega_{X^{(p)}}^1 \rightarrow 0.\]
Here $C$ is the Cartier operator and $W$ is the map in the following diagram: 
\begin{center}
    \begin{tikzcd}\label{Cartier}
        X \arrow{d} & X^{(p)} \arrow[l,"W"'] \arrow{d} & X \arrow[l,"F_{X/k}"'] \arrow[ld,""]\\
        \spec(k) & \spec(k) \arrow[l,"\sigma"]
    \end{tikzcd}
\end{center}

Our first claim is that the map $W^*-C: \HH^{i-2}(X,Z\Omega_X^1) \rightarrow \HH^{i-2}(X^{(p)},\Omega_{X^{(p)}}^1)$ is surjective. 
Indeed, as maps between $k$-vector spaces, $W^*$ is Frobenius-linear and surjective, 
and $C$ is a $k$-linear map, and
the difference between two such maps is always surjective by semi-linear algebra 
(see, \cite[Lemma II.5.3]{Il1}).

The square in the above diagram is Cartesian, $W$ is an isomorphism of schemes (but 
not of $k$-schemes), and $F_{X/k}$ is the relative Frobenius. 
$W^*$ stands for the composition \[Z\Omega_X^1 \xrightarrow[]{\cong} W^*(Z\Omega_X^1) \hookrightarrow W^*\Omega_X^1 \xrightarrow[]{\cong} \Omega_{X^{(p)}}^1.\]
Under the assumption, using the long exact sequence associated to \[0\rightarrow Z\Omega^1_X \rightarrow \Omega_X^1 \rightarrow B\Omega^2_X \rightarrow 0,\] 
we have a surjection $\HH^{i-2}(X,Z\Omega^1_X)\rightarrow \HH^{i-2}(X,\Omega^1_X)$. 
Composing with the Frobenius linear isomorphism $\HH^{i-2}(X,\Omega^1_X) \rightarrow \HH^{i-2}(X^{(p)},\Omega^1_{X^{(p)}})$, 
we thus have that $W^*:\HH^{i-2}(X,Z\Omega^1_X) \rightarrow \HH^{i-2}(X^{(p)},\Omega_{X^{(p)}}^1)$ is a $k$-Frobenius-linear surjection. 

The Cartier operator $C$ is $O_{X^{(p)}}$-linear.
We view $Z\Omega_X^1$ as an $O_{X^{(p)}}$-module because the differential in the de Rham complex is $O_{X^{(p)}}$-linear. 
So on global sections, $C: \HH^{i-2}(X,Z\Omega_X^1) \rightarrow \HH^{i-2}(X^{(p)},\Omega^1_{X^{(p)}})$ is $k$-linear. 
Now the first claim follows from \cite[Lemma II.5.3]{Il1}. 

By the long exact sequence, the first claim implies $\HH^{i-1}(X,\G_m/\G_m^{p}) \rightarrow \HH^{i-1}(X,Z\Omega_X^1)$ is injective. 
Composing with $\HH^{i-1}(X,Z\Omega^1) \rightarrow \HH^{i}_{\dR}(X/k)$ is still injective, this proves the case $n=1$. 

For general $n$, we use the diagram \ref{LESD} and induction on $n$. 
Proposition \ref{quasiiso} gives a quasi-isomorphism $W_{m+n}\Omega^{\bullet}/p^n\rightarrow W_n\Omega^{\bullet}$. 
Taking the projective limit over $m$, we obtain a quasi-isomorphism $W\Omega^{\bullet}/p^n\rightarrow W_n\Omega^{\bullet}$.
The long exact sequence associated to the distinguished triangle
\[0\rightarrow W\Omega^{\bullet} \xrightarrow{p^n} W\Omega^{\bullet} \rightarrow W_n\Omega^\bullet \rightarrow 0\]
gives \[...\rightarrow \HH^{i-1}_{\cris}(X/W_n)\rightarrow \HH^i_{\cris}(X/W)\xrightarrow{p^n}\HH^i_{\cris}(X/W) \rightarrow...\]
When $\HH^i_{\cris}(X/W)$ is torsion-free, $p^n$ is injective on it, 
and the connecting map 
$\HH^{i-1}_{\cris}(X/W_n) \rightarrow \HH^i_{\cris}(X/W)$ is $0$. 
By composing with the projection to $\HH^i_{\dR}(X/k)$, the connecting map $\HH^{i-1}_{\cris}(X/W_n)\rightarrow \HH^i_{\dR}(X/k)$ is also trivial. 
Suppose that for both $1$ and $n$ the $i$-th $d\log$ maps are injective, 
then the connecting map $\HH^{i-1}_{\fppf}(X,\mu_{p^n})\rightarrow \HH^i_{\fppf}(X,\mu_{p^1})$ is also $0$. 
We can use the 5-lemma to conclude.
\end{proof} 

\begin{proposition}\label{HdRdegene}
    Assume that $X$ is a proper smooth variety over $k$ such that the Fr\"olicher spectral sequence degenerates at $E_1$. Then the following statements hold:
    
{\rm (1)} the map $\HH^0(X,Z\Omega^1_X) \rightarrow \HH^0(X,\Omega_X^1)$ is an isomorphism;

{\rm (2)} the map $\HH^1(X,Z\Omega^1_X) \rightarrow \HH^2_{\dR}(X/k)$ is  injective.
\end{proposition}
\begin{proof}
     We use \v{C}ech cohomology. Suppose $\{U_t\}_{t\in T}$ is an affine open cover of $X$ and 
$C_i(\Omega_X^j)$ is the group of degree $i$ cochains on $\Omega_X^j$, that is, 
\[C_i(\Omega_X^j)=\prod_{t_0,t_1,...,t_i\in T}\Omega_X^j(U_{t_0}\cap U_{t_1}\cap U_{t_2} \cap ... \cap U_{t_i})\]
Then the \v{C}ech complex 
\[0 \rightarrow C_0(\Omega_X^j) \xrightarrow{\delta} C_1(\Omega_X^j) \xrightarrow{\delta} ...\]
computes the cohomology of $\Omega^j$, 
and the total complex of the following diagram computes $\HH^n_{\dR}(X/k)$: 
\begin{center}
    \begin{tikzcd}
        C_0(O_X) \arrow[r,"\delta"] \arrow[rd,"d"] & C_1(O_X)         \arrow[r,"\delta"]\arrow[rd,"d"]  & C_2(O_X)        \arrow[r,"\delta"]\arrow[rd,"d"]  & C_3(O_X)              \\
                                                   & C_0(\Omega_X^1)  \arrow[r,"\delta"]\arrow[rd,"d"]  & C_1(\Omega_X^1)  \arrow[r,"\delta"]\arrow[rd,"d"] & C_2(\Omega_X^1)       \\
                                                   &                                                    & C_0(\Omega_X^2)  \arrow[r,"\delta"]\arrow[rd,"d"] & C_1(\Omega_X^2)       \\
                                                   &                                                    &                                                   & C_0(\Omega_X^3)                                               
    \end{tikzcd}
\end{center}

The degeneration of the Fr\"olicher spectral sequence implies that
the cohomology 
of the total complex of the bottom three rows maps isomorphically to the $Fil^1$ part of $\HH^*_{\dR}(X/k)$. 

For (1), the degeneration implies that $\HH^0(X,\Omega^1_X)$ is isomorphic to the $Fil^1$ part of $\HH^1_{\dR}(X/k)$, which is $\HH^0(X,Z\Omega_X^1)$ by definition. So (1) holds. 

For (2), suppose that $\{h_{ij}\}$ represents a cohomology class in $\HH^1(X,Z\Omega_X^1)$ that becomes $0$ in $\HH^2_{\dR}(X/k)$. 
Clearly, this cohomology class lies in $\Fil^1(\HH^2_{\dR}(X/k))$, so the fact that
it becomes $0$ implies that 
there exists an element $\{\omega_i\}\in C_0(\Omega_X^1)$
such that $h_{ij}=\omega_i-\omega_j$ and $d\omega_i=0$.
It follows that $\{h_{ij}\}\in \HH^1(Z\Omega_X^1)$ is identically $0$. So (2) holds.
\end{proof}

\begin{corollary}\label{InjCol}
    Suppose that the Fr\"olicher spectral sequence for $X$ degenerates at $E_1$, 
    and $\HH^1_{\cris}(X/W)$ is torsion-free. Then the dlog map for $n=2$ is injective for all $m$: 
    \[\HH^2_{\fppf}(X,\mu_{p^m})\hookrightarrow \HH^2_{\cris}(X/W_m).\]
\end{corollary}
\begin{proof}
    This is a direct consequence of Propositions \ref{Injarg} and \ref{HdRdegene}. 
\end{proof}

\section{Parts of the Slope Spectral Sequence}

The goal of this section is to study parts of the slope spectral sequence.
As a corollary we show that all
flat Artin-Mazur formal groups $\Phi_{\fl}^i(A,\G_m)$ are smooth for all abelian 3-folds 
(and most of abelian 4-folds), without invoking Ekedahl's diagonal $t$-structure in \cite{Ek3}. 
We first recall Ekedahl's duality filtration in \cite{Ek1}. 

Recall that in \cite{Ek1} there is a duality functor $D:D(R)^{\op} \rightarrow D(R)$, formally defined as 
\[D(M):=R\varprojlim_{n}((R_{n}\otimes ^{\Ld}_{R}M)^\vee). \]
For a proper smooth variety $X/k$, Ekedahl proves that the following duality morphism is an isomorphism (see \cite{Ek1}):
\[D(R\Gamma(X,W\Omega^\bullet_{X/k}))[-N](-N)\rightarrow R\Gamma(X,W\Omega^\bullet_{X/k}).\]
For a coherent $R$-module $M$, Ekedahl defines $D^i(M):=\HH^i(D(M))$.
Then there is a duality filtration on $M$ 
\[0\subseteq T^2(M)\subseteq T^1(M)\subseteq T^0(M)=M \] 
coming from the duality spectral sequence: \[D^i(D^j(M))\Rightarrow \HH^{i-j}(D(D(M))),\quad (\text{we have}\;D(D(M))\cong M).\]
It is a fact that $T^2M\subset M[p^\infty]$ and the quotient $M[p^\infty]/T^2M$ 
is a direct sum of finite length $W$ modules with $F$ semisimple. 
Ekedahl obtaines the following description of graded factors $A^i(M):=T^i(M)/T^{i+1}M$.

\begin{definition}    
{\rm (1)} $A^2_f(M)^i:=\text{\rm {c\oe ur}}A^2(M)^i$ in degree $i$ has finite length, with $F$ nilpotent.

{\rm (2)} $A^2_d(M)^i:=$ the domino part of $A^2(M)$ in degree $i$.

{\rm (3)} $A^1_{ss}(M)^i:=$ the $p$-torsion part of $A^1(M)$ in degree $i$ has finite length with $F$ bijective.

{\rm (4)} $A^1_n(M)^i:=$ the free quotient $A^1(M)^i/A^1(M)^i[p^\infty]$ is a positive slope Dieudonné module.

{\rm (5)} $A^0(M)^i:=$ is free with $F$ bijective. 
\end{definition}

Ekedahl then proves the following result (Theorem IV.3.3 \cite{Ek1}): 

\begin{theorem} Let $M$ be a coherent $R$-module. Then the following properties hold.

    {\rm (1)} $A^0(M)$ and $D^0(M)$ are torsion-free, finitely generated $W$-modules with $F$ bijective, and 
    \[D_0(M)=\Hom_W(A_0(M),W),\quad (F,d,V)=(F^{-1,*}\sigma_*,0,pF^{-1}).\]
    This means that $[F(f)](x)=f(F^{-1}x)^\sigma$. 

    {\rm (2)} $A^1_n(M)$ and $D^1(M)/D^(M)[p^\infty]$ are torsion-free, 
    finitely generated $W$-modules with $F$ topologically nilpotent. 
    $A^1_{ss}(M)$ and $D^1(M)[p^\infty]$ are finite length $W$-modules with $F$ bijective. 
We have $A^1(M)=A^1_f(M)\oplus A^1_{ss}(M)$ as $R$-modules and 
    \[D^1(M)=\Hom_W(A^1_n(M),W)(1)\oplus \Hom_W(A^1_{ss}(M),K/W),\]
    \[(F,d,V)=(V^*\sigma_*,0,F^*\sigma^{-1}_*)\oplus (F^{-1,*}\sigma_*,0,pF^{-1}).\]

    {\rm (3)} $A^2_f(M)^i$ and $\text{\rm {c\oe ur}}(D^2(M))^i$ are both finite length $W$-modules 
    with $F$ nilpotent, and 
    \[D^2(M)=\Hom_W(\bigcup_n \mathrm{ker}(F^nd)\cap \mathrm{ker}(F^n),K/W)(1)\]
    \[(F,d,V)=(V^*\sigma_*,d^*,F^*\sigma^{-1}_{*})\]
    \[A^2_f(D^2(M))^i\cong D^2(A^2_f(M)^{-i-1})(-1), \quad
   A^2_d(D^2(M))^i\cong D^2(A^2_d(M)^{-i-2})(-2).\]

    {\rm (4)} We have $D^i(M)=0$ if $i\neq 0,1,2$. 
\end{theorem}

Ekedahl then proves the following duality theorem (Theorem IV.3.5 \cite{Ek1}): 

\begin{theorem}\label{Ekeduality}
    Let $X/k$ be a proper smooth variety of dimension $N$. Then the following statements hold.
    
{\rm (1)} $A^1(R^i\Gamma(X,W\Omega^\bullet_{X/k}))\cong D^1(R^{N+1-i}\Gamma(X,W\Omega^\bullet_{X/k}))(-N).$ 

{\rm (2)} There is a natural exact sequence of $R$-modules: 
\[\begin{aligned}
    0\rightarrow A^0(R^i\Gamma(X,W\Omega_{X/k}^\bullet))\rightarrow D^0(R^{N-i}\Gamma(X,W\Omega_{X/k}^\bullet))(-N) \\
    \xrightarrow{d_2} D^2(R^{N+1-i}\Gamma(X,W\Omega^\bullet_{X/k}))(-N) 
\rightarrow A^2(R^{i+1}\Gamma(X,W\Omega^\bullet_{X/k})) \rightarrow 0.
\end{aligned}\]

{\rm (3)} $A^1_n(R^i\Gamma(X,W\Omega_{X/k}^\bullet))^j\cong D^1(A^1_n(R^{N+1-i}\Gamma(X,W\Omega_{X/k}^\bullet))^{N-1-j})(-1)$;
\[A^1_{ss}(R^i\Gamma(X,W\Omega_{X/k}^\bullet))^j\cong D^1(A^1_{ss}(R^{N+1-i}\Gamma(X,W\Omega_{X/k}^\bullet))^{N-j})\]

{\rm (4)} 
$A^2_{f}(R^i\Gamma(X,W\Omega_{X/k}^\bullet))^j\cong D^2(A^2_{f}(R^{N+2-i}\Gamma(X,W\Omega_{X/k}^\bullet))^{N-1-j})(-1)$.

\end{theorem}

This theorem is basically saying that all parts in the slope spectral sequence of a smooth proper variety satisfy a
certain duality, namely: the slope $0$ part $A^0(\HH^{ij})$ is ``dual" to $A^0(\HH^{N-i,N-j})$ 
(not isomorphic to the dual! There is a finite cokernel); 
the positive slope part $A^1_n(\HH^{ij})$ is dual to $A^1_n(\HH^{N-1-i,N+1-j})$; 
the $F$-semi-simple part $A^1_{ss}(\HH^{ij})$ is dual to $A^1_{ss}(\HH^{N-i,N+1-j})$; 
the $F$-nilpotent torsion part $A^2_{f}(\HH^{ij})$ is dual to $A^2_f(\HH^{N-1-i,N+2-j})$;
the domino part $A^2_d(\HH^{ij})$ is dual to $A^2_{d}(\HH^{N-2-i,N+2-j})$. 
This is analogous to the Poincare duality of singular cohomology: 
the torsion part of degree $i$ is dual to the torsion part in degree $2N+1-i$. 

\medskip


    

We need some lemmas for maps between dominos. 

\begin{proposition}{\bf Maps from $U_0$ to $U_1$.}
    Let $U_0$ and $U_1$ be the dominos as defined in \ref{Coherence}. 
    Any map $f:U_0\rightarrow U_1$ is determined by the image of $1$, $f(1)$ is of the form $a+bV$. Moreover: 

    {\rm (1)} The map $f$ is injective if and only if  $a=0$, $b\neq 0$. 
    In this case, the cokernel of $f$ equals $k(0)$, $F=V=0$. 

    {\rm (2)} The map $f$ is surjective if and only if $a\neq 0$. 
    In this case, its kernel equals $k(-1)$ with $V=0$. 
    Moreover if $b\neq 0$, then $F$ bijective;
    if $b=0$, then $F=0$. 
\end{proposition}

\begin{proof}
    In $U_0$, $Fd(1)=0$, so $Fd(f(1))=f(Fd(1))=0$. 
The kernel of $Fd$ in $U_1$ is exactly $k\oplus kV$, this proves the first statement. 
When $f(1)=a+bV$, we have \[\begin{aligned}
    f(k_0+k_1V+...)&=k_0(a+bV)+k_1(a^{\sigma^{-1}}V+b^{\sigma^{-1}}V^2)+... \\
                   &=k_0a+(k_0b+k_1a^{\sigma^{-1}})V+(k_1b^{\sigma^{-1}}+k_2a^{\sigma^{-1}})V^2+...
\end{aligned}\]
\[\begin{aligned}
    f(k_0d+k_1dV+k_2dV^2+...)&=k_0d(f(1))+k_1dVf(1)+...\\
                             &=(k_0b+k_1a^{\sigma^{-1}})dV+(k_1b^{\sigma^{-1}}+k_2a^{\sigma^{-2}})dV^2+...
\end{aligned}\]

If $a=0$, $b\neq 0$, then $f$ is injective with cokernel $k(0)$, $F=V=0$. 

If $a\neq 0$, then $f$ is surjective with kernel $k(-1)$, generated by a vector of the form 
\[k_0d+k_0\frac{b}{a^{\sigma^{-1}}}dV+k_0\frac{b}{a^{\sigma^{-1}}}\frac{b^{\sigma^{-1}}}{a^{\sigma^{-2}}}dV^2+...\]

If moreover $b=0$, then $F=V=0$. If $b\neq 0$, then $F$ is bijective and $V=0$. 
\end{proof} 

\begin{proposition}{\bf Maps from $U_0$ to $U_{\ell}$.}
    Let $0\leqslant \ell$, let $U_0$ and $U_{\ell}$ be the dominos as defined in \ref{Coherence}. 
    Any map $f:U_0\rightarrow U_{\ell}$ is determined by the image of $1$, 
and $f(1)$ is of the form $a_0+a_1V+...+a_{\ell}V^{\ell}$.

    Moreover: 
    
    {\rm (1)} The map $f$ is injective if and only if $a_0=...=a_{\ell-1}=0$ and $a_{\ell}\neq 0$. 
    In this case, the cokernel $M$ concentrates on degree $0$, and $F=0$, $V$ is nilpotent ($V^{\ell}=0$). 

    {\rm (2)} The map $f$ is surjective if and only if $a_0\neq 0$. 
    In this case the kernel concentrates in degree $1$. 

\end{proposition}

\begin{proof}
    We have $Fd(1)=0$ so $Fd(f(1))=0$ as well. The kernel of $Fd$ in $U_{\ell}$ is the direct sum 
$k\oplus kV\oplus...\oplus kV^{\ell}$, so $f(1)$ must be of the form \[a_0+a_1V+...+a_{\ell}V^{\ell}.\] 

If $a_i$ is non-zero for some of $i=0,...,\ell-1$, say $a_0=...=a_{i-1}=0$ but $a_i\neq 0$, 
then for every $k_0,k_1,...,k_{\ell-i-1}$ we can suitably choose the rest of the coefficients such that $f(k_0d+k_1dV+k_2dV^2+...)=0$. 
This proves (1). 

If $a_0=0$ then $f(k[[V]])\subseteq Vk[[V]]$ so $f$ is not surjective. This proves (2).
\end{proof} 

\begin{proposition}{\bf Maps from $U_k$ to $U_{\ell}$.} \label{Injdim1}
    Let $k\leqslant \ell$, let $U_0$ and $U_{\ell}$ be the dominos as defined in \ref{Coherence}. 
    Any map $f:U_k\rightarrow U_{\ell}$ is determined by the image of $1$, $f(1)$ is of the form 
 $a_0+a_1V+...+a_{\ell-k}V^{\ell-k}$.

    Moreover: 
    
    {\rm (1)} The map $f$ is injective if and only if $a_0=...=a_{\ell-k-1}=0$ and $a_{\ell-k}\neq 0$. 
    In this case, the cokernel concentrates on degree $0$. 

    {\rm (2)} The map $f$ is surjective if and only if $a_0\neq 0$. 
    In this case the kernel concentrates on degree $1$. 
\end{proposition}

\begin{lemma}{\bf Maps from $U_{\ell}$ to an arbitrary domino.}\label{InjUD}
    Let $f:U_{\ell}\rightarrow D$ be a map from a $1$-dimensional domino $U_{\ell}$ to an arbitrary domino $D$ of dimension $n\geqslant 1$. 
    
    {\rm (1)} If $f$ is injective, then the cokernel $D/U_{\ell}$ can be written as: 
    \[0\rightarrow M\rightarrow D/U_{\ell}\rightarrow D'\rightarrow 0\]
    where $M$ is a finite dimensional $k$-vector space with $F=0$, $V$-nilpotent; $D'$ is a domino of dimension $n-1$.
\end{lemma}

\begin{proof}
     Let $U'\subset D$ be the following module: 
$U'(0) =V^{-\infty}f(U(0))$, $U'(1)=F^\infty f(U(1))$. 
Then we can show that $D/U'$ is a domino $D'$. 
Thus the cokernel sits in a short exact sequence \[0\rightarrow U'/U\rightarrow D/U\rightarrow D'\rightarrow 0.\]
The cokernel $U'/U$ is concentrated in degree $0$, so we are done.
\end{proof} 

\begin{proposition}{\bf Maps between dominos.}\label{MapsDD}
    Let $f: D_1\rightarrow D_2$ be a map between two dominos of the same dimension, such that its kernel and cokernel do not contain any dominos. 
    Then the kernel of $f$ concentrates in degree $1$, the cokernel of $f$ concentrates in degree $0$. 

    In other words, the following statements hold: 

    {\rm (1)} Let $f: D_1\rightarrow D_2$ be an injective map between two dominos of the same dimension, then its cokernel concentrates in degree $0$. 

    {\rm (2)} Let $f: D_1\rightarrow D_2$ be an surjective map between two dominos of the same dimension, then its kernel concentrates in degree $1$. 
\end{proposition}

\begin{proof}
     We use induction on the dimension of the dominos. When $\dim(D_1)=\dim(D_2)=1$, this is Proposition \ref{Injdim1}. 
Assume $d:=\dim(D_1)=\dim(D_2)\geq 2$. Choose a subdomino  $U\subseteq D_1$ of dimension $1$ such that $D_1/U$ is still a domino of dimension $d-1$. 
Consider the map $U\hookrightarrow D_1\rightarrow D_2$, 
by Proposition \ref{InjUD}, the cokernel $D_2/Im(U)$ is an extension of a domino $D_3$ by a finite length module $M[0]$. 
We thus get a map of dominos of dimension $d-1$:  \[D_1/U\xrightarrow{\widetilde{f}} D_2/Im(U)\rightarrow D_3.\]
We have the following commutative diagram 
\begin{center}
    \begin{tikzcd}
         0 \arrow{r} & U \arrow{r} \arrow{d} & D_1 \arrow{r} \arrow{d} & D_1/U \arrow{r} \arrow{d} &0\\
         0\arrow{r} & Im(U) \arrow{r} & D_2 \arrow{r} & D_2/Im(U) \arrow{r} &0
    \end{tikzcd}
\end{center}
Now we consider the snake lemma: the kernel of $D_1/U\rightarrow D_2/Im(U)$ 
is a submodule of the kernel $D_1/U\rightarrow D_2/Im(U)\rightarrow D_3$, 
by induction hypothesis it concentrates on degree $1$. 
The cokernel of $D_1/U\rightarrow D_2/Im(U)$  is an extension of cokernel of $D_1/U\rightarrow D_3$ by $M[0]$, 
so it concentrates on degree $0$. Now the snake lemma gives us what we wanted.
\end{proof} 

\begin{proposition}\label{Partsandhearts}
    
    Let $M$ be a coherent $R$-module, let $T^2M\subseteq M[p^\infty]\subseteq M$ be Ekedahl's duality filtration. 
    
    {\rm (1)} We have injective maps: $\text{\rm {c\oe ur}}(T^2M^j)\hookrightarrow \text{\rm {c\oe ur}}(M[p^\infty]^j)\hookrightarrow \text{\rm {c\oe ur}}(M^j)$. 

    {\rm (2)} $\text{\rm {c\oe ur}}(M[p^\infty])$ is identified with the torsion module of $\text{\rm {c\oe ur}}(M^j)$. 

    {\rm (3)}  $\text{\rm {c\oe ur}}(T^2M^j)$ is identified with the nilpotent part of  $\text{\rm {c\oe ur}}(M[p^\infty])$, 
    the quotient $\text{\rm {c\oe ur}}(M[p^\infty])/\text{\rm {c\oe ur}}(T^2M^j)$ (isomorphic to 
    $M[p^\infty]^j/T^2M^j\cong A^1_{ss}(M)^j$), is thus isomorphic to the semi-simple torsion part of $\text{\rm {c\oe ur}}(M^j)$. 
\end{proposition}

\begin{proof}
In general, if $N \subseteq M$ is a coherent $R$-submodule, by definiton we have 
$V^{-\infty}Z(N^j)=V^{-\infty}Z(M^j)\cap N^j$, and $F^\infty B(N^j)\subseteq F^\infty B(M^j)\cap N^j$. 
We now use $T\subseteq N\subseteq M$ to denote Ekedahl's duality filtration $T^2M\subseteq M[p^\infty]\subseteq M$. 
We have injective maps of dominos: 
\begin{center}
    \begin{tikzcd}
        T^j/V^{-\infty}Z^j(T) \arrow{r}\arrow{d} &F^{\infty}B^j(T) \arrow{d}\\
        T^j/V^{-\infty}Z^j(N) \arrow{r}\arrow{d} &F^{\infty}B^j(N) \arrow{d}\\
        T^j/V^{-\infty}Z^j(M) \arrow{r}          &F^{\infty}B^j(M) 
    \end{tikzcd}
\end{center}
By Proposition \ref{MapsDD}, we have equalities $F^{\infty}B^j(T)\cong F^{\infty}B^j(N)\cong F^{\infty}B^j(M)$. 
We have the following commutative diagram with exact rows
\begin{center}
    \begin{tikzcd}
        0\arrow{r} & F^{\infty}B^j(N) \arrow{r}\arrow[d,"\cong"] &V^{-\infty}Z^j(N) \arrow{r}\arrow{d} & \text{\rm {c\oe ur}}(N^j) \arrow{r}\arrow{d} &0\\
        0\arrow{r} & F^{\infty}B^j(M) \arrow{r}          &V^{-\infty}Z^j(M) \arrow{r}          & \text{\rm {c\oe ur}}(M^j) \arrow{r}           &0
    \end{tikzcd}
\end{center}
It induces an injective map $\text{\rm {c\oe ur}}(N^j)\rightarrow \text{\rm {c\oe ur}}(M^j)$ such that 
\[ \text{\rm {c\oe ur}}(M^j)/\text{\rm {c\oe ur}}(N^j)\cong V^{-\infty}Z^j(M)/V^{-\infty}Z^j(N).\] 
The latter is free because $V^{-\infty}Z(N^j)=V^{-\infty}Z(M^j)\cap N^j$. This proves (2). 

We have the following diagram: 
\begin{center}
    \begin{tikzcd}
        0\arrow{r} & F^{\infty}B^j(T) \arrow{r}\arrow[d,"\cong"] &V^{-\infty}Z^j(T) \arrow{r}\arrow{d} & \text{\rm {c\oe ur}}(T^j) \arrow{r}\arrow{d} &0\\
        0\arrow{r} & F^{\infty}B^j(N) \arrow{r}          &V^{-\infty}Z^j(N) \arrow{r}          & \text{\rm {c\oe ur}}(N^j) \arrow{r}           &0
    \end{tikzcd}
\end{center}
that induces an injective map $\text{\rm {c\oe ur}}(T^j)\rightarrow \text{\rm {c\oe ur}}(N^j)$, and it gives an isomorphism
\[\text{\rm {c\oe ur}}(N^j)/\text{\rm {c\oe ur}}(T^j)\cong V^{-\infty}Z^j(N)/V^{-\infty}Z^j(T).\]

We also have the following diagram 
\begin{center}
    \begin{tikzcd}
        0\arrow{r} &V^{-\infty}Z^j(T)  \arrow{r}\arrow{d} &T^j \arrow{r}\arrow{d} & T^j/V^{-\infty}Z^j(T) \ \arrow{r}\arrow{d} &0\\
        0\arrow{r} &V^{-\infty}Z^j(N)  \arrow{r}                 & N^j \arrow{r}          & N^j/V^{-\infty}Z^j(N) \ \arrow{r}           &0
    \end{tikzcd}
\end{center}
that gives a short exact sequence 
\[0\rightarrow V^{-\infty}Z^j(N)/V^{-\infty}Z^j(T)\rightarrow A^1_{ss}(M^j)\rightarrow (N^j/V^{-\infty}Z^j(N))/(T^j/V^{-\infty}Z^j(T))\rightarrow 0\]
The $F$ action on $A^1_{ss}(M^j)$ is bijective, and the $F$-action on \[(N^j/V^{-\infty}Z^j(N))/(T^j/V^{-\infty}Z^j(T))\] is nilpotent, so 
we must have  $(N^j/V^{-\infty}Z^j(N))/(T^j/V^{-\infty}Z^j(T))=0$ and  $V^{-\infty}Z^j(N)/V^{-\infty}Z^j(T)\cong A^1_{ss}(M^j)$. 
As a result we see that
$$\text{\rm {c\oe ur}}(N^j)/\text{\rm {c\oe ur}}(T^j)\cong V^{-\infty}Z^j(N)/V^{-\infty}Z^j(T)\cong A^1_{ss}(M^j)$$ is $F$-semisimple. 
This completes our proof.
\end{proof}

\begin{proposition}\label{A1ssH11}
    Let $X/k$ be a smooth proper variety of dimension $N$. 
    
    {\rm (1)} We have $A^1_{ss}(\HH^1(X,W\Omega^\bullet_{X/k}))^0=0$, and 
$$A^1_{ss}(\HH^1(X,W\Omega^\bullet_{X/k}))^1=\NS(X)[p^\infty]\otimes_{\Z_p} W,$$
where $F$ acts trivially on $\NS(X)[p^\infty]$. 

    {\rm (2)} We have $A^1_{ss}(\HH^N(X,W\Omega^\bullet_{X/k}))^{N-1}=\NS(X)[p^\infty]^\vee\otimes_{\Z_p}W.$ Here $\vee$ denotes the Pontryagin dual, and $F$ acts trivially on $\NS(X)[p^\infty]^\vee$. 
\end{proposition}

\begin{proof}
 We have    $A^1_{ss}(\HH^1(X,W\Omega^\bullet_{X/k}))^0=0$ since
$\HH^1(X,WO_X)$ is the Cartier module of $\widehat{\Pic}_{X/k}$, which is always free. 
Next, we have 
$A^1_{ss}(\HH^1(X,W\Omega^\bullet_{X/k}))^1=\NS(X)[p^\infty]\otimes_{\Z_p} W$.
Indeed,
by Proposition \ref{Partsandhearts} we have \[A^1_{ss}(\HH^1(X,W\Omega^\bullet_{X/k}))^1\cong \text{\rm {c\oe ur}}(\HH^1(X,W\Omega^1_{X/k}))_{ss}.\] 
By Proposition \ref{log} we have $\HH^1(X,W\Omega^1_{X/k,\log})\cong \HH^2(X,\Z_p(1))$.
We also have $\HH^2(X,\Z_p(1))[p^\infty]\cong \NS(X)[p^\infty]$. To see this, note that
by Kummer sequence we have a long exact sequence
\[\HH^0(X,\G_m)\xrightarrow{p^n} \HH^0(X,\G_m)\rightarrow \HH^1(X,\mu_{p^n})\rightarrow \HH^1(X,\G_m)\xrightarrow{p^n} \HH^1(X,\G_m),\]
hence $\HH^1(X,\mu_{p^n})\cong \Pic(X)[p^n]$. Passing to the projective limit we obtain an isomorphism
$\HH^1(X,\Z_p(1))\cong T_p(\Pic(X))$. 
We also have the exact sequence 
\[\HH^1(X,\Z_p(1))\xrightarrow{p^n}\HH^1(X,\Z_p(1))\rightarrow \HH^1(X,\mu_{p^n})\rightarrow \HH^2(X,\Z_p(1)) \xrightarrow{p^n} \HH^2(X,\Z_p(1))\]
Taking a large enough $n$ gives $\HH^2(X,\Z_p(1))[p^\infty]\cong \NS(X)[p^\infty]$. Part
(2) is a consequence of (1) and Theorem \ref{Ekeduality}(3). 
\end{proof}

\begin{remark} {\rm
    If $X/k$ is a surface, this gives another proof of Theorem \ref{surfaces}(1) that the finite group part of 
    $\HH^3(X,\Z_p(1))$ equals $\NS(X)[p^\infty]^\vee$. 
    This shows that Ekedahl's duality theorem \ref{Ekeduality} is consistent with Milne's flat duality theorem \ref{FlatDuality}. }
\end{remark}

We have the following proposition explaining each parts of the slope spectral sequence of an arbitrary 
proper smooth surface. For a proper smooth surface $S$, we have the following invariants: 
\begin{itemize}
    \item Betti numbers: $b_n:=\dim(\HH^n_{\et}(S,\Q_{\ell}))$.
    \item Hodge numbers: $h^{ij}$.
    \item Irregularity: $q=\frac{b_1}{2}$.
    \item The finite group scheme $H:=\mathrm{coker}(\Pic_{\red,X/k}^\circ\rightarrow \Pic_{X/k}^\tau)$. 
    \item We have a decomposition $H\cong H_{el}\oplus H_{ll}\oplus H_{le}$, where $H_{el}$ is $\NS(X)[p^\infty]$. 
\end{itemize}


\begin{theorem}\label{surfacesparts}
    Let $S/k$ be a proper smooth surface. Then the following statements hold.

{\rm (1)} $\HH^0(S,W\Omega^i_{S/k})$ are all free for $i=0,1,2$. 

{\rm (2)} $A^1_{ss}(\HH^1(S,W\Omega^\bullet_{S/k}))^0=0$, 

\indent\indent\indent $A^1_{ss}(\HH^1(S,W\Omega^\bullet_{S/k}))^1\cong(\NS(S)[p^\infty]\otimes_{\Z_p}W)\cong H_{el}\otimes_{\Z_p}W$.  

\indent\indent\indent $A^1_{ss}(\HH^1(S,W\Omega^\bullet_{S/k}))^2$ is $\DM(H_{le})^\vee$. 

\indent\indent\indent $A^2_{f}(\HH^1(S,W\Omega^\bullet_{S/k}))^i=0$ for all $i=0,1,2$. 
That is, $\HH^1(S,W\Omega^\bullet_{S/k})$ has no nilpotent torsion. 

{\rm (3)} $A^1_{ss}(\HH^2(S,W\Omega^\bullet_{S/k}))^0\cong\DM(H_{le})$ 

\indent\indent\indent $A^1_{ss}(\HH^2(S,W\Omega^\bullet_{S/k}))^1\cong\NS(S)[p^\infty]^\vee\otimes_{\Z_p}W\cong(H_{el}\otimes_{\Z_p}W)^\vee$,

\indent\indent\indent $A^1_{ss}(\HH^2(S,W\Omega^\bullet_{S/k}))^2=0$,

\indent\indent\indent $A^2_{f}(\HH^2(S,W\Omega^\bullet_{S/k}))^0\cong\DM(H_{ll})$. 

\indent\indent\indent $A^2_{f}(\HH^2(S,W\Omega^\bullet_{S/k}))^1\cong\DM(H_{ll})^\vee$. 

\end{theorem}

We are abusing notations here. The action $(F,V)$ on the dual $M^\vee$ has different definitions: if $M$ is a torsion with $F$ bijective, then $(F,V)_{M^\vee}=(F^{-1,*}\sigma_{*},pF^{-1})$; 
if $M$ is torsion with $F$ nilpotent, then $(F,V)_{M^\vee}=(V^{*}\sigma_{*},F^*\sigma^{-1}_*)$. 

\begin{proof}
     (1) is because the sheaves $W\Omega^i_{S/k}$ are all torsion-free (\cite{Il1}). 
The first and second statements of (2) are Proposition \ref{A1ssH11}. 
The third statement of (2), by Theorem \ref{Ekeduality}(3) is equivalent to the first statement of (3), 
which is a consequence of Theorem \ref{Eke81}. 
The fourth statement of (2), according to Theorem \ref{Ekeduality}(4) is equivalent to the vanishing of 
$A_{f}^2(\HH^{N+1}(S,W\Omega^\bullet_{S/k}))$, which is trivial. 
The second and third statement of (3) by duality is equivalent to the first and second statement of (2) respectively. 
The fourth and fifth statement of (3) are equivalent to each other and is a consequence of Theorem \ref{Eke81}.
\end{proof}

\begin{remark}{\rm
    Theorem \ref{surfacesparts} shows that every finite torsion part in the slope spectral sequence 
    of a surface $S/k$ comes from a certain finite torsion subgroup of $\Pic_{S/k}$, the Picard scheme. 
    Conversely, each finite torsion part of $\Pic_{S/k}$ (local-\'etale, or \'etale-local, or local-local parts) 
    is reflected in a certain finite torsion part of the slope spectral sequence of $S$. }
\end{remark}

We can compute the Hodge numbers of a surface with $T^{02}=0$ from its slope spectral sequences, generalizing a result of Suwa \cite{Su}. 
Suwa's calculation for surfaces satisfying $T^{02}=m^{02}=0$ extends to the case $m^{02}\neq 0$. 

\begin{theorem}\label{HodgeofSurfaceT02}
    Let $S$ be a proper smooth surface over an algebraically closed field $k$ of characteristic $p$. 
    Suppose that the unipotent number $T^{02}=0$, then the following statements hold.

    {\rm (1)}
    We have the following table of Hodge numbers for $S$:
    \begin{center}{\rm\label{Hodgetable}
        \begin{tabular}{ c|c|c|c } 
         $\dim$    &   $\Omega^0_{S/k}$           & $\Omega^1_{S/k}$              & $\Omega^2_{S/k}$ \\
         $\HH^2$   &   $m^{02}+\mathrm{rk}(_FH)$ &     $q+\mathrm{rk}(_VH)$        &  1   \\
         $\HH^1$   &    $q+\mathrm{rk}(_FH)$      &$b_2-2m^{02}+2\mathrm{rk}(_VH)$&  $q+\mathrm{rk}(_FH)$ \\           
         $\HH^0$   &    1                         &     $q+\mathrm{rk}(_VH)$        &   $m^{02}+\mathrm{rk}(_FH)$.\\ 
        \end{tabular}}
    \end{center}

    {\rm (2)}  We have the following table of de Rham numbers for $S$: 
    \[h_{\dR}^n=\begin{cases}
        1           &\text{when}\quad n=0,4\\
        b_1+\rk(_pH)&\text{when} \quad n=1,3\\
        b_2+2\rk(_pH)&\text{when} \quad n=2.
    \end{cases}\]
\end{theorem}

Here we follow the notation of Suwa \cite{Su}: $_FH$ means the kernel of the map $H\xrightarrow{F_{H}}H^{(p)}$, 
$_VH$ means the kernel of the map $H\xrightarrow{V_{H}}H^{(1/p)}$, $_pH$ means the kernel of the map $H\xrightarrow{p}H$. 

\begin{proof}
The short exact sequence \[0\rightarrow WO_X\xrightarrow{V}WO_X \rightarrow O_X\rightarrow 0\]
gives the long exact sequence
\[\begin{aligned}
     &0\rightarrow \HH^1(WO_X)\xrightarrow{V}\HH^1(WO_X)\rightarrow \HH^1(O_X)\rightarrow \HH^2(WO_X)\\
     &\xrightarrow{V} \HH^2(WO_X)\rightarrow \HH^2(O_X)\rightarrow 0.
\end{aligned}\]
We have \[\begin{aligned}
    &\dim_k(\HH^1(WO_X)/V)=q,\\
    &\dim_k( \HH^2(WO_X)[V])=\rk_FH,\\
    &\dim_k( \HH^2(WO_X)/V)=m^{02}+\rk_FH
\end{aligned}\]
so the first column (and the third column, by Serre duality) of (1) follows. 

To compute the cohomology of $\Omega^1_{X/k}$, we use the exact sequence: 
\[0\rightarrow WO_{X} \xrightarrow{(F,-Fd)}WO_X\oplus W\Omega^1_{X/k} \xrightarrow{dV+V}W\Omega^1_{X/k}\rightarrow \Omega^1_{X/k}\rightarrow 0.\]
Let $\mathrm{Im}$ be the image sheaf of $dV+V$. We have short exact sequences
\[0\rightarrow WO_{X} \xrightarrow{(F,-Fd)}WO_X\oplus W\Omega^1_{X/k} \rightarrow\mathrm{Im}\longrightarrow 0,\]
\[0\rightarrow \mathrm{Im} \rightarrow W\Omega^1_{X/k}\rightarrow \Omega^1_{X/k}\rightarrow 0.\]
We have long exact sequences: 
\[\begin{aligned}\label{Im1}
    0\rightarrow \HH^0(WO_X) \xrightarrow{(F,-Fd)} \HH^0(WO_X) \oplus \HH^0(W\Omega^1_{X/k}) \longrightarrow \HH^0(\mathrm{Im})\\
    \rightarrow \HH^1(WO_X) \xrightarrow{(F,-Fd)} \HH^1(WO_X) \oplus \HH^1(W\Omega^1_{X/k}) \longrightarrow \HH^1(\mathrm{Im})\\
    \rightarrow \HH^2(WO_X) \xrightarrow{(F,-Fd)} \HH^2(WO_X) \oplus \HH^2(W\Omega^1_{X/k}) \longrightarrow \HH^2(\mathrm{Im})\rightarrow 0,
\end{aligned}\]

\[\begin{aligned}\label{Im2}
    0\rightarrow \HH^0(\mathrm{Im}) \xrightarrow{} \HH^0(W\Omega^1_X)  \longrightarrow \HH^0(\Omega^1)\\
    \rightarrow \HH^1(\mathrm{Im}) \xrightarrow{} \HH^1(W\Omega^1_X)  \longrightarrow \HH^1(\Omega^1)\\
    \rightarrow \HH^2(\mathrm{Im}) \xrightarrow{} \HH^2(W\Omega^1_X)  \longrightarrow \HH^2(\Omega^1)\rightarrow 0. 
\end{aligned}\]

We can compute $\HH^2(\Omega^1)$ as $\HH^2(W\Omega^1)/(dV+V)$. If $T^{02}=0$, then $dV=0$, 
so $\HH^2(W\Omega^1)/(dV+V)\cong \HH^2(W\Omega^1)/V$, in this case we have 
\[0\rightarrow A_{ss}^1\oplus A_n^2 \rightarrow \HH^2(W\Omega^1)\rightarrow \HH^2(W\Omega^1)/[p^\infty]\rightarrow 0,\]
the last term is free and $V$ is injective on it. 
We can compute \[\begin{aligned}
    &\dim(\HH^2(W\Omega^1)/V)\\
    =&\dim(A_{ss}^1(\HH^2(W\Omega^1))/V)+\dim(A_n^2(\HH^2(W\Omega^1))/V)\\
     &+\dim(A^0(\HH^2(W\Omega^1))/V)+\dim(A^1_n(\HH^2(W\Omega^1))/V)\\
    =&\dim_{\F_p}\NS(S)[p]+\rk _VH_{ll}+m^{12}=q+\rk_VH. 
\end{aligned}\]

The rest of (1) follows from Serre's duality and the formula: $2\chi(O_X)-\chi(\Omega^1_{X/k})=c_2$. 
This formula is a consequence of the Hodge-de Rham spectral sequence and the long exact sequence: 
\[...\rightarrow \HH^n_{\cris}(S/W)\xrightarrow{p}\HH^n_{\cris}(S/W) \rightarrow \HH^n_{\dR}(S/k)\rightarrow...\]
The split degeneration of slope spectral sequence when $T^{02}=0$, and the above long exact sequence can give us (2).
\end{proof}

Ekedahl's duality also can be used to show that the flat Artin-Mazur formal groups 
of an arbitrary abelian 3-fold are all smooth. 

\begin{theorem}\label{flatAMFsmooth}
    Let $A/k$ be an abelian variety of dimension $g$. The following statements hold: 

{\rm  (1)}    when $g=3$, then the flat Artin-Mazur formal groups 
    $\Phi^i_{\fl}(A,\G_m)$ are all smooth for $i=1,2,3$; 

{\rm  (2)}    when $g=4$, the only possible case that at least one of the flat Artin-Mazur formal groups 
$\Phi^i_{\fl}(A,\G_m)$ is not smooth is the flat formal Brauer group $\Phi^2_{\fl}(A,\G_m)$ of an abelian $4$-fold whose slopes of 
$\HH^1_{\cris}(A/W)$ are \[\{0,0,\frac{1}{2},\frac{1}{2},\frac{1}{2},\frac{1}{2},1,1\}.\] 

{\rm  (3)}    when $A$ is supersingular,  $A_f^2$ and $A_{ss}^1$ are trivial everywhere. 
As a corollary, the flat Artin-Mazur formal groups 
$\Phi^i_{\fl}(A,\G_m)$ are all smooth for all $i$.  
\end{theorem}

\begin{proof}
    Let us first prove (1). By Ekedahl's result \ref{Eke81}, 
it suffices to show that $\HH^i(A,WO_A)$ has no $V$-torsion for $i=2,3$. 
The $V$-torsion of these groups are canonically isomorphic to the direct sum of the semi-simple torsion $A_{ss}^1$, 
and the nilpotent torsion $A_{f}^2$ at these 
positions by Proposition \ref{Partsandhearts}. Since the Picard scheme is smooth, and $A$ is $3$-dimensional, 
it suffices to show $\Phi^2_{\fl}(A,\G_m)$ is smooth. By Ekedahl's duality, the semi-simple torsion at $(i,j)$ is dual to the semi-simple torsion at 
$(N-i,N+1-j)$, the nilpotent torsion at $(i,j)$ is dual to the nilpotent torsion at $(N-1-i,N+2-j)$, 
so it suffices to show that $A_{ss}^1$ of $\HH^1(A,W\Omega^3_{A/k})$ is trivial, and 
$A_{f}^2$ of $\HH^2(A,W\Omega^2_{A/k})$ is also trivial. 
The groups $\HH^i(A,W\Omega^3_{A/k})$ are free (\cite[Proposition II.3.12]{Il1}), so the first statement holds. 
For $\HH^2(A,W\Omega^2_{A/k})$, notice that $\HH^1(A,W\Omega^3_{A/k})$ 
is trivial unless $A$ is ordinary or almost ordinary. 
In these two cases $A$ is Hodge-Witt, so of course $\HH^2(A,W\Omega^2_{A/k})$ does not have nilpotent torsion; 
in other cases, the nilpotent torsion $A_{f}^2$ of $\HH^2(A,W\Omega^2_{A/k})$ 
is the heart of a submodule of $\HH^2(A,W\Omega^2_{A/k})$, so it survives to the $\infty$-page, 
but the trivialness of $\HH^1(A,W\Omega^3_{A/k})$ makes sure that $E_\infty^{2,2}$ is free. So the result follows. 

Let us prove (2). It suffices to prove that $\Phi_{\fl}^2(A,\G_m)$ and $\Phi_{\fl}^3(A,\G_m)$ are smooth. 
By Ekedahl's result \ref{Eke81}, 
it suffices to show that $\HH^i(A,WO_A)$ doesn't have $V$-torsion for $i=3,4$.
For an abelian variety, there is no semi-simple torsion (\ref{AbeT}). 
It suffices to show that the nilpotent torison $A_f^2$ of 
$\HH^2(A,W\Omega_{A/k}^3)$ and $\HH^3(A,W\Omega_{A/k}^3)$ are trivial. 
The nilpotent torison of $\HH^2(A,W\Omega_{A/k}^3)$ is trivial for all $A$: 
again, $\HH^1(X,W\Omega^4_{A/k})$ is trivial unless $A$ is ordinary or almost ordinary (Hodge-Witt), 
but in other cases the nilpotent torsion $A_{f}^2$ of $\HH^2(A,W\Omega^3_{A/k})$ 
is the heart of a submodule of $\HH^2(A,W\Omega^3_{A/k})$, so it survives to the $\infty$-page, 
but the trivialness of $\HH^1(A,W\Omega^4_{A/k})$ makes sure that $E_\infty^{2,2}$ is free. 
For the nilpotent torsion of $\HH^3(A,W\Omega_{A/k}^3)$, 
we notice that $\HH^2(A,W\Omega^4_{X/k})$ is non-trivial 
only when the multiplicity of slope $0$ in $\HH^1_{\cris}(A/W)$ is at least $2$. 
The only cases when $A$ is not Hodge-Witt, is when the slopes of 
$\HH^1_{\cris}(A/W)$ are $\{0,0,\frac{1}{2},\frac{1}{2},\frac{1}{2},\frac{1}{2},1,1\}$. 
In other cases, either $A$ is Hodge-Witt, or  $\HH^2(A,W\Omega^4_{X/k})$ is trivial so 
the nilpotent torsion $A_{f}^2$ of $\HH^3(A,W\Omega^3_{A/k})$ 
is the heart of a submodule of $\HH^3(A,W\Omega^3_{A/k})$, so it survives to the $\infty$-page, 
but the trivialness of $\HH^2(A,W\Omega^4_{A/k})$ makes sure that $E_\infty^{3,3}$ is free.

Finally let us prove (3). By \ref{AbeT} there is no semi-simple torsion anywhere. 
Since $A$ is supersingular, $E_\infty^{i,j}$ is torsion whenever $j<i$. 
They are all $0$ since we can use induction to prove that they are torsion subobjects of $\HH^{i+j}_{\cris}(A/W)$. 
This further shows that the heart of $\HH^j(A,W\Omega^i_{A/k})$, for any $j\leqslant i+1$, is free. 
Therefore there is no nilpotent torsion $A_f^2$ whenever $j\leqslant i+1$. 
Now by Ekedahl duality, we conclude that there is no nilpotent torsion anywhere. 
\end{proof} 

\begin{corollary}
    Let $A/k$ be an abelian variety of dimension $g$ 
    over $k$ an algebraically closed field of characteristic $p$. 
    Then the \'etale Artin-Mazur formal group
    $\Phi^i(A,\G_m)$ are all pro-representable  on $\Art^{\op}_k$ and they agree with $\Phi^i_{\fl}(A,\G_m)$, for all the following cases: 
    
    {\rm (1)} $g=3$. 

    {\rm (2)} $g=4$ and the slopes of $\HH^1_{\cris}(A/W)$ are not $\{0,0,\frac{1}{2},\frac{1}{2},\frac{1}{2},\frac{1}{2},1,1\}$. 

    {\rm (3)} Supersingular.  
\end{corollary}

\begin{proof}
     This is a consequence of \ref{flatAMFsmooth} and Proposition 10.11 of \cite{BO}. 
\end{proof}

\begin{remark}{\rm
    Abelian varieties are Mazur-Ogus in the sense of \cite{Ek3}. 
    According to Theorem IV.1.2 and Theorem III.4.6 in \cite{Ek3}, 
    it should be true that for any abelian variety $A/k$, $A_f^2$ and $A_{ss}^1$ are all trivial everywhere. 
    From this we can conclude that the flat Artin-Mazur formal groups $\Phi^i_{\fl}(A,\G_m)$ are formally smooth for all $n$. 
    This further shows that on $\Art^{\op}_k$, the (\'etale) Artin-Mazur formal group 
    $\Phi^i(A,\G_m)$ are all pro-representable as well and they agree with $\Phi^i_{\fl}(A,\G_m)$. 
    I don't know how to prove this using elementary method without using Ekedahl's diagonal complexes.}
\end{remark}

\section{On the Tate Conjecture for Divisors}

In this section we explain how to view Tate conjecture for divisors (over a finite field)
 as an analogue of Hodge conjecture over $\Cm$. 
Let $k=\overline{\F}_p$, the algebraic closure of the finite field. 
Let $X/k$ be a proper smooth variety. Since $X$ must be defined over a finite field, 
there is a proper smooth variety $X_0/\F_q$ such that $X_0\otimes_{\F_q}k\cong X$. 
Then the Tate conjecture for divisors for $X_0$ states that the following map is surjective: 
\[(T_{\ell}): \Pic(X_0)\otimes \Q_\ell \longrightarrow \HH^2_{\et}(X,\Q_{\ell}(1))^{Gal(k/\F_q)}.\] 
The $p$-adic version for the crystalline cohomology states that the following map is surjective: 
\[(T_{p}): \Pic(X_0)\otimes \Q_{p} \longrightarrow \HH^2_{\cris}(X_0/W(\F_q))_{\Q}^{F=p}. \]

There is a nice viewpoint to view Tate conjecture as an analogue of Hodge conjecture, based on the following fact: 
the union of $W(\F_{p^n})$, the ring of Witt vectors of the finite field $\F_{p^n}$ in $W$ does not equal $W$. 
This is because, if we denote $\breve{W}$ the union, then $\breve{W}$ is not complete, but $W$ is complete. 
This is very similar to the fact that $\Q$ is not complete in $\Cm$. 
Inside $\HH^2_{\cris}(X/W)$, the union \[\breve{\HH}^2_{\cris}(X/\breve{W}):=\bigcup_{n}\HH^2_{\cris}(X_0\otimes_{\F_q}\F_{q^n}/W(\F_{q^n}))\]
is the subspace of $\HH^2_{\cris}(X/W)$ consisting of the elements that are "uniformly finite". 
This is the analogue of $\HH^2(X,\Q)\subset \HH^2(X,\Cm)$ in the Hodge conjecture situation.
To compare, the $F=p$ subspace of $\HH^2_{\cris}(X/W)$ is the analogue of the $(1,1)$-piece: 
$\HH^{1,1}(X,\Cm):=\HH^1(X,\Omega^1_{X/C})$ in the Hodge conjecture setting. 

Then the Tate conjecture for divisors maybe reformulated as: the following map is an isomorphism:
\[(T_p'): \Pic(X)\otimes \Q_p \longrightarrow \breve{\HH}^2_{\cris}(X/\breve{W})_{\Q} \cap \HH^2_{\cris}(X/W)_{\Q}^{F=p}.\]

Viewing this, the right form of Torelli theorem for K3 surfaces should be: let $X/k$ be a K3 surface over $k=\overline{\F}_p$, 
then the isomorphism class of $X$ depends only on the isomorphism class of (non-completed version) F-crystal with a perfect pairing:
\[(\breve{\HH}^2_{\cris}(X/\breve{W}),F,\langle-,-\rangle).\]
When $X/k$ is not supersingular, analogous to the complex case,
the isomorphism class of $(\HH^2_{\cris}(X/W),F,\langle-,-\rangle)$ only depends on the height of $X$. 

For Hodge conjecture, the divisor case is proved as the Lefschetz $(1,1)$-theorem. 
One proof of the Lefschetz $(1,1)$-theorem uses the exponential map of sheaves: 
\[0\rightarrow \Z \rightarrow \mathcal{O}_X \xrightarrow{\exp} \mathcal{O}_X^* \rightarrow 0. \]

In the Tate conjecture setting there is no exponential map. 
However, there is the following Artin-Hasse exponential map (\ref{ArtinHasse}), viewing both sides as sheaves on various sites over $k$:  
\begin{equation}\begin{aligned}
    \mathrm{E}:\widehat{\W}\times \W\rightarrow \G_m, \\ \label{ArtinHasse}
    \mathrm{E}((u,x))=\mathrm{E}(u\cdot x).
\end{aligned}\end{equation}
It is a question what does this map give us about Tate conjecture. 
For example, we have \[\widehat{\W}\cong \mathcal{H}om(\W,\G_a), \quad
R\Gamma(X,\widehat{\W})\cong R\Gamma(X,\mathcal{H}om(\W,\G_a)).\]
It would be ideal if we could find a short exact sequence of sheaves of the following form: 
\[0\rightarrow \breve{\W} \rightarrow \W \rightarrow \G_m \rightarrow 0, \]
where the cohomology groups of $\breve{\W}$ are exactly $\breve{\HH}^*_{\cris}(X/\breve{W})$. 
If this is done, we could prove Tate conjectures for divisors by analyzing the short exact sequence 
\[\HH^1(X,\G_m)\rightarrow \breve{\HH}^2_{\cris}(X/\breve{W})\rightarrow \HH^2(X,WO_X),\]
because if $x\in \breve{\HH}^2_{\cris}(X/\breve{W})\cap \HH^2_{\cris}(X/W)^{F=p}$, 
its image in $\HH^2(X,WO_X)$ has to be $0$ by Illusie's theory. So it must comes from $\Pic(X)$. 

\appendix

\begin{appendices}  
    \chapter{Properties of the De Rham-Witt complex}
    
    The main reference for this section is \cite{Il1}. 
    
    \begin{definition}{\bf V-pro complexes.}\label{Vpro}
        Let $X$ be a ringed topos. A V-pro complex on $X$ is a projective system of commutative differential $\Z$--graded algebras
        \[\{M_n^\bullet: M_{n+1}^\bullet \xrightarrow{R} M_n^\bullet\}\]
        equipped with a family of additive maps $V:M_n^\bullet\rightarrow M_{n+1}^\bullet$
        satisfying the following properties: 
    
        {\rm (0)} $RV=VR$. 
    
        {\rm (1)} $M_n^\bullet=0$ if $n\leqslant 0$, $M_1^0$ is an $\F_p$-algebra, and for $n\geqslant 1$, $M_n^0=W_n(M_1^0)$, 
        and $R,V$ are the usual restriction and shifting operators on Witt vectors: 
        $R:W_{n+1}(M_1^0)\rightarrow W_n(M_1^0)$, $V:W_{n}(M_1^0)\rightarrow W_{n+1}(M_1^0)$. 
    
        {\rm (2)} For each $n,i,j$, and $x\in M_n^i$, $y\in M_n^j$, we have $V(xdy)=Vxd(Vy)$. 
    
        {\rm (3)} For each $n\geqslant 1$, and each element $x\in M_1^0$, $y\in M_n^0$, 
        then we have $(Vy)d\underline{x}=V(\underline{x}^{p-1}y)d(V\underline{x})$, 
        where $\underline{x}=(x,0,...,0)\in M_n^0$ is the multiplicative lift.
    \end{definition}
    
    Note that condition (2) implies: 
    
    {\rm (2')}  For each $n$ and $a,x_1,...,x_i\in M_n^0$, we have \[V(adx_1dx_2...dx_i)=(Va)(dVx_1)...(dVx_i).\]
    
    Conditions (2) and (2') are equivalent when each $M_n^\bullet$ is generated by $M_n^0$ as a differential graded algebra. 
    Illusie explicitly constructed $W_n\Omega_A^\bullet$ inductively as a quotient of $\Omega^\bullet_{W_n(A)}$, 
    and proved the following:  
    \begin{theorem}
        The functor from V-pro complexes to $\F_p$-algebras: 
        \[V\text{-}\mathrm{pro}(X)\rightarrow \F_p\text{-}\mathrm{Alg}(X): \{M_n^\bullet\}\rightarrow M_1^0\]
        admits a right adjoint \[A\rightarrow \{W_n\Omega^\bullet_A\},\] 
        in the sense that, for any V-pro complexes $\{M_n^\bullet\}$, we have a natural isomorphism 
        \[\Hom_{V\text{-}\mathrm{pro}(X)}(\{W_{n} \Omega^{\bullet}_A\},\{M_{n}^{\bullet}\})=\Hom_{\F_p\text{-}\mathrm{Alg}(X)}(A,M_1^0).\]
    \end{theorem}
    Moreover, the construction of $W_n\Omega^\bullet_A$ is compatible with \'etale base change. 
    The resulting pro-complex of \'etale sheaf on $X$ formed by $A\rightarrow W_n\Omega^\bullet_A$ 
    defines the de Rham-Witt pro-complex $W_n\Omega^\bullet_{X}$. 
    By taking projective limit, we obtain $W\Omega_{A}^\bullet:=\varprojlim_{n}W_n\Omega_A^\bullet$. 
    
    Until now, we have not discussed the Frobenius operator $F$ yet. 
    The definition of $F$, 
    is highly dependent on the explicit computation of $W_{n}\Omega^\bullet_A$ for a smooth algebra $A$, 
    particularly, for those $A$ of the form $\F_p[t_1,...,t_m,t_{m+1}^{\pm1},...,t_n^{\pm 1}]$. 
    Illusie defined  the operator $F$ explicitly for such algebras, 
    and then showed that this definition generalizes to the general case by a limit argument. 
    
    \begin{theorem}
        Let $X$ be a ringed topos of $\F_p$-algebras. 
        Then the projective sysmtem of ring homomorphisms: $FR=RF: W_{n+1}O_X\rightarrow W_nO_X$ 
        extends uniquely to a projective system of graded algebra homomorphisms: 
        \[F:W_{n+1}\Omega^\bullet_{X}\rightarrow W_n\Omega^\bullet_{X},\]
        such that: 
    
        {\rm (1)} for each $n\geqslant 2$ and each $x\in O_X$, 
        let $\underline{x}_{\leqslant n}=(x,0,...,0)\in W_nO_X$ be the multiplicative lift, we have 
        \[Fd\underline{x}_{\leqslant n}=\underline{x}_{\leqslant n-1}^{p-1}d\underline{x}_{\leqslant n-1}.\]
    
        {\rm (2)} for each $n\geqslant 1$, $FdV=d:W_{n}\Omega^i_{X}\rightarrow W_n\Omega^{i+1}_X$. 
    \end{theorem}
    
    \begin{proposition}The following formulas hold: 
        
        {\rm (1)} $FV=VF=p: W_n\Omega^i_X\rightarrow W_n\Omega^i_X.$
    
        {\rm (2)} $dF=pFd: W_n\Omega^i_X\rightarrow W_{n-1}\Omega^{i+1}_X$, and
                  $Vd=pdV: W_n\Omega^i_X\rightarrow W_{n+1}\Omega^{i+1}_X$.
    
        {\rm (3)} $FdV=d: W_n\Omega^i_X\rightarrow W_n\Omega^{i+1}_X$.
    
        {\rm (4)} $xVy=V(Fx.y)$, for any $x\in W_n\Omega^i_X,\;y\in W_{n}\Omega^j_X$.
    
        These formulas also hold for the projective limit $W\Omega^\bullet_{X/k}$. 
    \end{proposition}

    From now on, assume $X$ is a smooth scheme over a perfect base $S$ of characteristic $p$. 
    We will need the following results on the structure of \( W\Omega_X^\bullet \).
    
    \medskip
    
    {\bf Canonical Filtration.}
    \begin{definition}\label{CanFildef}
    For each $n,r$, we define the canonical filtration on $W_r\Omega^\bullet_X$ and $W\Omega_{X}^\bullet$ to be
    \begin{equation}
    \label{CanFil}\Fil^nW_r\Omega^\bullet_X:=\begin{cases}
        W_r\Omega_X^\bullet &\text{if}\ n\leqslant 0  \;\text{or}\; r\leqslant 0\\
        \Ker(R^{r-n}: W_r\rightarrow W_n) &\text{if}\ 1\leqslant n < r\\
        0 &\text{if}\ n\geqslant r\\
    \end{cases}
    \end{equation}
    
    \[\Fil^nW\Omega^\bullet_X:=\begin{cases}
        W\Omega_X^\bullet &\text{if}\ n\leqslant 0 \\
        \Ker(R_n: W\Omega^\bullet_X\rightarrow W_n\Omega^\bullet_{X}) &\text{if}\ n\geqslant 1 \\
    \end{cases}\]
    
    We then naturally define the graded pieces: 
    \[\gr^nW_r\Omega^\bullet_{X}:=\Fil^nW_r\Omega^\bullet_{X}/\Fil^{n+1}W_r\Omega^\bullet_{X},\]
    \[\gr^nW\Omega^\bullet_{X}:=\Fil^nW\Omega^\bullet_{X}/\Fil^{n+1}W\Omega^\bullet_{X}.\]
    
    \end{definition}
    \begin{proposition}
        For each $n\geqslant 0$ and each $r$, $\Fil^nW_r\Omega^\bullet_X$ 
        is the differential graded ideal generated by $V^nW_{r-n}O_X$. 
        In other words, we have 
        \[\Fil^nW_r\Omega_X^i=V^nW_{r-n}\Omega_X^i+dV^nW_{r-n}\Omega_X^{i-1},\]
        \[\Fil^nW\Omega_X^i=V^nW\Omega_X^i+dV^nW\Omega_X^{i-1}.\]
    \end{proposition}
    
    \begin{proposition}
        For each $n$ and each $i$, the homomorphism $F:W_{n+1}\Omega^i_{X}\rightarrow W_n\Omega^i_{X}$ induces 
        (via quotient by $\Fil^{n}$) a homomorphism
        \[F: W_n\Omega^i_X\rightarrow W_n\Omega^i_{X}/dV^{n-1}\Omega^{i}_X.\]
        For $n=1$, this is the inverse Cartier operator
        \[C^{-1}:\Omega^i_{X/k}\rightarrow \Omega^i_{X/k}/d\Omega^{i-1}_{X/k}.\]
    \end{proposition}
    
    \begin{proposition}\label{kerP}
        For each $n$, the multiplication by $p$ map 
        $p: W_{n+1}\Omega_X^\bullet\rightarrow W_{n+1}\Omega_X^\bullet$ has kernel $\Fil^n W_{n+1}\Omega_X^{\bullet}$, 
        so it induces an injection
        $p: W_n\Omega_X^\bullet\rightarrow W_{n+1}\Omega_X^\bullet$. 
        Moreover, for each $i$ such that $0\leqslant i\leqslant n+1$, we have 
        \[\Ker(p^i: W_{n+1}\Omega^\bullet\rightarrow W_{n+1}\Omega^\bullet)=\Fil^{n+1-i}W_{n+1}\Omega_X^\bullet.\]
    \end{proposition}
    
    \medskip
    
    {\bf Graded Pieces of the Canonical Filtration.}
    
    \begin{proposition}\label{kerVndVn}
        For each $n$ and each $i$, we have the following morphism of exact sequences:
        \begin{center}
            \begin{tikzcd}
                0\arrow{r} & B_n\Omega^i_X \arrow{r}\arrow{d} & \Omega_X^i \arrow[r,"V^n"]\arrow{d} & W_{n+1}\Omega_X^i \arrow{d} \\
                0\arrow{r} & B_{n+1}\Omega^i_X \arrow{r} & \Omega_X^i \arrow[r,"V^n"] & W_{n+1}\Omega_X^i/dV^n\Omega_X^{i-1}
            \end{tikzcd}
        \end{center}
        \begin{center}
            \begin{tikzcd}
                0\arrow{r} & Z_{n+1}\Omega^{i-1}_X \arrow{r}\arrow{d} & \Omega_X^{i-1} \arrow[r,"dV^n"]\arrow{d} & W_{n+1}\Omega_X^i \arrow{d} \\
                0\arrow{r} & Z_{n}\Omega^{i-1}_X \arrow{r} & \Omega_X^{i-1} \arrow[r,"dV^n"] & W_{n+1}\Omega_X^i/V^n\Omega_X^{i}
            \end{tikzcd}
        \end{center}
        where the right vertical arrows are the natural projections. 
    \end{proposition}

    The proof of Theorem \ref{kerVndVn} relies heavily on the analysis in the next subsection. We will use these properties freely. 
    
    \begin{proof}
         By \'etale localization we may assume $S=\spec(k)$, $X=\spec(A)$, $A=k[T_1,...,T_r]$. 
    Suppose $B=W[T_1,...,T_r]$, let $\Omega_B^\bullet$ the de Rham complex of $B/W$, 
    and let $E^\bullet$ the complex of integral logarithmic forms used in the construction of $F$ 
    ($E^\bullet$ is the de Rham-Witt complex of $A$), we have for each $m$, a quasi-isomorphism 
    \[W_m\Omega_A^\bullet\cong E^\bullet/\Fil^m E^\bullet,\]
    and for each $i$ we have
    \[\Fil^mE^i=V^mE^i+dV^mE^{i-1}.\]
    Since $E$ is generated by $E^0$ (as a commutative differential graded algebra) 
    and $E^0=\sum_{j\geqslant 0}V^jB$, we can see that 
    \[\Fil^mE^i=\sum_{j\geqslant0} V^{m+j}\Omega_B^i+dV^{m+j}\Omega_B^{i-1}.\]
    
    Suppose $x\in \Omega_A^i$ is a section who has a lift $y\in \Omega_B^i$. 
    Suppose that $V^ny\in \Fil^{n+1}E^i$, 
    we can write $V^n$ as the following finite sum: 
    \[V^ny=\sum_{1\leqslant m\leqslant N}V^{n+m}z_m+dV^{n+m}t_m,\]
    where $z_m\in \Omega^i_B$ and $t_m\in \Omega^{i-1}_B$. Apply $F^{n+N}$ on both sides, 
    we obtain that\[p^nF^Ny=\sum_{1\leqslant m\leqslant N} p^{n+m}F^{N-m}z_m+ F^{N-m}dt_m\]
    Multiplying $p^{N-1}$ on both sides, we obtain
    \[p^{n+N-1}(F^Ny-\sum_{1\leqslant m\leqslant N}p^mF^{N-m}z)=d(\sum_{1\leqslant m\leqslant N}p^{m-1}F^{N-m}t_m).\]
    Now we apply Corollary \ref{BnZnlifting}(2), since $F^Ny-\sum_{1\leqslant m\leqslant N}p^mF^{N-m}z$ is a lift of $C^{-N}x$, 
    we conclude that $C^{-N}x\in B_{n+N}\Omega^i_A/B_n\Omega^i_A$. This implies that $x\in B_N\Omega^i_A$. 
    So we finished the proof of the first part. The second part of the proof is similar. 
    \end{proof}

    \begin{corollary}\label{Gradedpieces}
        For each $n\geq 0$ and each $i$, we have the following exact sequences of finite type locally free $O_X$-modules: 
        \[\begin{aligned}
            &0\rightarrow F^{n+1}_*\Omega^i_X/B_n^i \xrightarrow{V^n} \gr^nW\Omega^i_{X/k} \xrightarrow{\beta} F^{n+1}_*\Omega_X^{i-1}/Z_n^{i-1}\rightarrow 0\\
            0& \rightarrow F^{n+1}_*\Omega_X^{i-1}/Z_{n+1}^{i-1} \xrightarrow{dV^n} \gr^nW\Omega^i_{X/k} \xrightarrow{\beta'} F^{n+1}_*\Omega_X^{i}/B_{n+1}^i \rightarrow 0,
        \end{aligned}\]
        where the $O_X$-module of $\gr^nW\Omega^i_{X/k}$ is defined by \[F:O_X\cong W_{n+1}O_X/VW_{n}O_X\rightarrow W_{n+1}O_X/pW_{n+1}O_X,\]
        and $\beta$ maps the class of $V^nx+dV^ny$ to the class of $y$, $\beta'$ maps the class of $V^nx+dV^ny$ to the class of $x$. 
    \end{corollary}

    \medskip

    {\bf The kernel and cokernel of $F$ and $V$.}

    \begin{proposition}
        For each $n\geq 0$ and each $i$, there is a exact sequence 
        \[0\rightarrow \Omega^i_{X/k}/B_n\Omega^i_X \xrightarrow{V^n} W_{n+1}\Omega^i_{X/k} \xrightarrow{F} W_n\Omega^i_X\]
        \[0\rightarrow \Omega_X^{i-1}/Z_n\Omega_X^{i-1} \xrightarrow{dV^n} W_{n+1}\Omega_X^i \xrightarrow{V} W_{n+2}\Omega_X^i \]
    \end{proposition}

    \medskip 

    {\bf The canonical filtration and the p-adic filtration.}

    \begin{proposition}
        For each integer $n\geq 0$, 
    \end{proposition}

    \begin{lemma}\label{quasiisop}
        For each $n\geqslant 0$, multiplication by $p$ induces a quasi-isomorphism
        $\gr^nW\Omega_X^\bullet \rightarrow \gr^{n+1}W\Omega_X^\bullet$.
    \end{lemma}
    \begin{proof}
        This is \cite[Proposition I.3.13]{Il1}.
    \end{proof}  

    \medskip
    {\bf The fixed point of $F$}

    For each $n$, there is a map $d\log:O_X^* \rightarrow W_n\Omega^1_{X/k}$, 
    a homomorphism between sheaves of abelian groups defined by $x\rightarrow d\underline{x}/\underline{x}$
    where $\underline{x}=(x,0,0,...)$ is the multiplicative lift of $x\in W_nO_X$. 

    \begin{proposition}\label{GmWOmega1log}
        For each $n\geq 1$, the following sequence is exact: 
        \[0\rightarrow O_X^{*p^n} \rightarrow O_X^* \xrightarrow[]{d\log} W_n\Omega_{X/k}^1\]
    \end{proposition}

\begin{proof}
    Since $p^nW_n\Omega^1_{X/k}=0$, we have $d\log (O_X^{*p^n})=0$ in $W_n\Omega_{X/k}^1$. 
    Conversely, suppose $x\in O_X^*$ such that $d\log(x)=0$, we prove inductively that locally $x$ is of the form 
    $y^{p^n}$ for some local section $y$. 
    Suppose in some open neighborhood $x=y^{p^m}$ for some $0\leqslant m< n$. 
    Then $0=d\log(x)=p^md\log(y)$ and Proposition \ref{kerP} 
    shows that $d\log(y)\in \Fil^{n-m}W_n\Omega^1_{X/k}\subset \Fil^1W_n\Omega^1_{X/k}$.
    So the image of $d\log(y)\in \Omega^1_{X/k}$ equals $0$, this tells us that $y=y'^p$. 
    By induction, the claim is true. 
\end{proof}

    \chapter{Propositions for a Local Lifting}

    The proof of Proposition \ref{kerVndVn} is not trivial at all. 
    We need the following proposition in the following context. 
    Assume that $X\rightarrow S$ is smooth, where $S$ is a characteristic $p$ base which lifts to 
a $p$-adic complete and separate formal scheme $T$ such that 
    $O_T/pO_T\cong O_S$, and there is a endomorphism $F_T$ of $T$ that lifts the (absolute) Frobenius $F_S$. Assume also that
$X\rightarrow S$ lifts to smooth map between formal schemes $Y\rightarrow T$ such that there is an endomorphism 
    $F_Y$ that lifts the absolute Frobenius $F_X$. 
    When $S$ is perfect and when $X$ is affine, such lifts always exist. 
    We can take $T=W(S)$ and take $F_T$ to be the Frobenius of Witt vectors, 
    and since $X$ is affine, there is no obstruction to lift $X/S$ to a formally flat map $Y/T$, and we take $F_Y:Y\rightarrow T\times_{T,F_T}T$. 
    More generally, the lifts exist when both $X$ and $S$ are affine and $S$ is smooth over a perfect scheme.  
    
    We let $Y^{F/T}=Y^{F}$ be the formal scheme otained by fiber product 
    and we have the 'relative Frobenius' $F_{Y/T}:Y\rightarrow Y^{F}$ obtained by universal property. 
    Let $\Omega^1_{Y/T}$ be the continuous differential of $Y/T$, let $\Omega^\bullet_{Y/T},d$ be the corresponding de Rham complex. 
    
    \begin{proposition}
        Under the above conventions, 
    
        {\rm (1)} $F_{Y/T}$ is finite, and locally free.
    
        {\rm (2)} For each integer $i\geqslant 0$, the image of the canonical arrow 
        $F_{Y/T}^*\Omega^i_{Y^F/T}\rightarrow \Omega^i_{Y/T}$ (resp, 
        $\Omega^i_{Y^F/F}\rightarrow F_{Y/T,*}\Omega^i_{Y/T}$) is contained in $p^i\Omega^i_{Y/T}$
        (resp, $p^iF_{Y/T,*}\Omega^i_{Y/T}$). 
    \end{proposition} 
    
   Part (2) follows from the fact that when reducing modulo $p$ to $X/S$, 
the map $F_{X/S}^*\Omega^1_{X^{(p)}/S}\rightarrow \Omega^1_{X/S}$ is trivial. 
    This tells us that there exists a unique homomorphism of graded algebras 
    \[F: \Omega^i_{Y^F/T}\rightarrow \Omega^i_{Y/T}\] 
    such that $dF=pFd$, and that if $F_{Y/T}$ is the canonical map 
    $\Omega^i_{Y^F/T}\rightarrow \Omega^i_{Y/T}$ then $F_{Y/T}=p^iF$. 
    Let $\overline{F}: \Omega^\bullet_{X^{(p)}/S}\rightarrow \Omega^\bullet_{X/S}$ be the reduction of $F$. 
    Then (Mazur) gives the following observation: 
    
    \begin{proposition}
        We have for each $i$, $\overline{F}(\Omega^i_{X^{(p)}/S})\subseteq Z\Omega^i_{X/S}$, and the following composition map is the identity: 
        \[\Omega^i_{X^{(p)}/S}\xrightarrow{\overline{F}} Z\Omega^i_{X/S}\rightarrow \HH^i(\Omega^\bullet_{X/S},d)\xrightarrow{C_{X/S}} \Omega^i_{X^{(p)}/S}. \]
    \end{proposition}
    
    \begin{proof}
         In degree $0$, $\overline{F}$ is the relative Frobenius $F_{X/S}$. 
    In degree $1$, let $y\in O_Y$ be a local section whose reduction is $x\in O_X$. 
    By definition of the Cartier operator, $C^{-1}W^*dx$ is the class of $x^{p-1}dx\in \HH^1(\Omega^\bullet_{X/S},d)$. 
    By the definition of $\overline{F}$, the class of $\overline{F}(W^*dx)$ is the class of $F(W^*dy)=\frac{1}{p}(d(F_YW^*y))$ modulo $p$. 
    We have $F_YW^*y=y^p+pz$ for some $z\in O_Y$, 
    so the class of $F(W^*dy)$ equals the class of $y^{p-1}dy+dz$, modulo $p$ we get the class of $y^{p-1}dy\in \HH^1(\Omega^\bullet_{X/S},d)$, 
    so we proved the proposition. 
    \end{proof}
    
    Now let $Y^{F^n}=Y^{(n)}$ be the $n$-th twist and let $F^n: Y\rightarrow Y^{(n)}$ and 
    $F^n:\Omega^i_{Y^{(n)}/T}\rightarrow \Omega^i_{Y/T}$ be the $n$-fold 
    composition of $F$. Let $\overline{F^n}: \Omega^i_{X^{(n)}/S}\rightarrow \Omega^i_{X/S}$ be the reduction of $F^n$, 
    then we have: 
    
    \begin{corollary}\label{F=C-1}
        For each $i$, $\overline{F^n}(\Omega^i_{X^{(n)}/S})\subseteq Z_n\Omega^i_{X/S}$, 
        and the following composition map is identity: 
        \[\Omega^i_{X^{(n)}/S}\xrightarrow{\overline{F^n}} Z_n\Omega^i_{X/S}\rightarrow Z_n\Omega^i_{X/S}/B_n\Omega^i_{X/S} \xrightarrow{C_{X/S}^n} \Omega^i_{X^{(n)}/S}. \]
    
        In particular, we have 
        \[\overline{F^n}(\Omega^\bullet_{X^{(n)}/S})+B_n\Omega^\bullet_{X/S}=Z_n\Omega_{X/S}^\bullet\]
        \[\overline{F^n}(B_1\Omega^\bullet_{X^{(n)}/S})+B_n\Omega^\bullet_{X/S}=B_{n+1}\Omega_{X/S}^\bullet.\]
    \end{corollary}
    
    \begin{proposition}
        For each integers $n\geq 0$ and $i\geq 0$, let $d^{-1}p^n\Omega^{i+1}_{Y/T}$ 
        be the subsheaf of $\Omega^{i}_{Y/T}$ formed by sections $x$ such that $dx\in p^n\Omega^{i+1}_{Y/T}$. 
        We have: 
        \[d^{-1}p^n\Omega^{i+1}_{Y/T}=\sum_{0\leqslant k\leqslant n}p^kF^{n-k}(\Omega^i_{Y^{(n-k)}/T})+\sum_{0\leqslant k\leqslant n-1}F^k(d\Omega^{i-1}_{Y^{(k)}/T}).\]
    \end{proposition}

    \begin{proof}
         We use induction on $n$ to prove the theorem. When $n=0$, the proposition is clear. 
    Suppose the proposition is proved for $n-1$, let us prove $n$. 
    Let $x\in \Omega_{Y/T}^i$ such that $dx=p^ny$, where $y\in\Omega^{i+1}_{Y/T}$. 
    By induction hypothesis, after localization, we can write $x$ as the following form: 
    \[x=\sum_{0\leqslant k\leqslant n-1}p^kF^{n-k}x_k+\sum_{0\leqslant k\leqslant n-2}F^k(dy_k)\]
    where $x_k\in \Omega^i_{Y^{(n-k)}/T}$ and $y_k\in \Omega^{i-1}_{Y^{(k)}/T}$. We deduce, 
    \[dx=p^{n-1}\sum_{k=0}^{n-1}F^{n-k}dx_k.\]
    Assume \[\sum_{k=0}^{n-1}F^{n-k}dx_k=py,\] modulo $p$, we have \[\sum_{k=0}^{n-1}\overline{F}^{n-k}d\overline{x_k}=0.\]
    Since $\overline{F}^{n-2}d\overline{x}_1,...,d\overline{x}_{n-1}$ all lie in $B_{n-1}\Omega^i_{X/S}$, 
    so is $\overline{F^{n-1}}d\overline{x}_0$. This shows that $d\overline{x}_0=0$. 
    Repeating this analysis, we conclude that $d\overline{x}_0=...=d\overline{x}_{n-1}=0$. 
    We can thus (according to Corollary \ref{F=C-1}) write $x_k=F(u_k)+dv_k+pw_k$, as a result we have 
    \[\begin{aligned}
        x&=F^nu_0+pF^{n-1}(u_1+w_0)+...+p^{n-1}F(u_{n-1}+w_{n-2})+p^nw_{n-1}\\
        &+d(y_0+p^{n-1}v_{n-1})+Fd(y_1+p^{n-2}v_{n-2})+...+F^{n-2}d(y_{n-2}+pv_1)+F^{n-1}dv_0.
    \end{aligned}\]
    \end{proof}
    
    \begin{corollary}\label{BnZnlifting}
        Suppose integers $n\geqslant 0$ and $i\geqslant 0$, and $x\in \Omega^i_{X/S}$ is a local section. 
    
        {\rm (1)} For $x$ to be in $Z_n\Omega^i_{X/S}$, it is necessary and sufficient that there is a local section 
        $y\in \Omega^i_{Y/T}$ lifting $x$ such that $dy\in p^n\Omega^i_{Y/T}$. 
    
        {\rm (2)} For $x$ to be in $B_n\Omega^i_{X/S}$, it is necessary and sufficient that there is a local section 
        $y\in \Omega^{i-1}_{Y/T}$ such that $dy\in p^{n-1}\Omega^i_{Y/T}$, and $dy/p^{n-1}$ lifts $x$. 
    \end{corollary}
\end{appendices}
    

\addcontentsline{toc}{chapter}{Bibliography}
\bibliographystyle{alpha}
\bibliography{bibliography/bibliography}

\end{document}